\documentclass{article}
\usepackage{geometry}
\usepackage{float}
\usepackage{amssymb}
\usepackage{amsmath}
\usepackage{amsfonts}
\usepackage{amsthm}
\usepackage[T1]{fontenc}
\usepackage{lmodern,microtype,mathtools,mathrsfs,booktabs,graphicx,xcolor}
\usepackage[shortlabels]{enumitem}
\usepackage[colorlinks,linkcolor=blue,anchorcolor=blue,citecolor=red]{hyperref}
\usepackage[all]{xy}
\usepackage{perpage}
\usepackage{needspace}
\MakePerPage{footnote}

\newtheoremstyle{theorem}
{3pt}
{3pt}
{\itshape}
{}
{\bfseries}
{.}
{.5em}
{}
\theoremstyle{theorem}
\newtheorem{theorem}{Theorem}[section]
\newtheorem{lemma}[theorem]{Lemma}
\newtheorem{proposition}[theorem]{Proposition}
\newtheorem{corollary}[theorem]{Corollary}

\newtheorem{example}[theorem]{Example}

\newtheoremstyle{definition}
{3pt}{3pt}{\normalfont}{}{\bfseries}{.}{.5em}{}
\theoremstyle{definition}
\newtheorem{definition}[theorem]{Definition}
\theoremstyle{remark}
\newtheorem{remark}[theorem]{Remark}

\numberwithin{equation}{section}

\newcommand\D{{\rm d}}

\newcommand{\pps}[1]{\partial_{#1}}

\newcommand{\norm}[2]{\left\| #1\right\|_{#2}}
\newcommand{\Lnorm}[2]{{\left\| #1 \right\|}_{{L}^{#2}}}

\newcommand{\abs}[1]{\left| #1 \right|}
\newcommand{\set}[1]{\left\{ #1 \right\}}
\newcommand{\tnorm}[2]{\| #1\|_{#2}}

\newcommand{\tset}[1]{\{ #1 \}}

\newcommand{\ball}[2]{B_{#1}(#2)}
\newcommand{\cball}[2]{\bar{B}_{#1}(#2)}
\newcommand{\lsp}[3]{{#1}^{#2}(#3)}

\newcommand{\lssp}[4]{{#1}_{#2}^{#3}(#4)}

\newcommand{\id}{\mathrm{id}}
\DeclareMathOperator{\TP}{TP}

\DeclareMathOperator{\diag}{diag}

\DeclareMathOperator{\dist}{dist}
\DeclareMathOperator{\supp}{supp}
\DeclareMathOperator{\sgn}{sgn}

\DeclareMathOperator{\Span}{span}

\DeclareMathOperator{\clos}{clos}

\DeclareMathOperator{\Orb}{Orb}

\DeclareMathOperator{\ind}{ind}
\DeclareMathOperator{\interior}{int}

\DeclareMathOperator{\Diff}{Diff}
\DeclareMathOperator{\Cont}{Cont}
\DeclareMathOperator{\loc}{loc}

\newcommand {\F}{\mathbb F}

\newcommand {\R}{\mathbb R}

\newcommand {\N}{\mathbb N}
\newcommand {\Z}{\mathbb Z}

\newcommand{\rn}[1]{\R^{#1}}

\newcommand{\bs}[1]{\boldsymbol{#1}}

\newcommand{\ssp}[2]{_{#1}^{#2}}
\newcommand{\inv}{^{-1}}

\newcommand{\mca}[1]{\mathcal{#1}}

\DeclareFontFamily{U}{mathx}{\hyphenchar\font45}
\DeclareFontShape{U}{mathx}{m}{n}{<-> mathx10}{}
\DeclareSymbolFont{mathx}{U}{mathx}{m}{n}
\DeclareFontSubstitution{U}{mathx}{m}{n}
\DeclareMathAccent{\widebar}{\mathalpha}{mathx}{"73}

\makeatletter
\newcommand*{\bigcdot}{\mathpalette\bigcdot@{.6}}  
\newcommand*{\bigcdot@}[2]{%
  \mathbin{\vcenter{\hbox{\scalebox{#2}{$\m@th#1\bullet$}}}}%
}
\makeatother

\DeclareMathOperator{\Crit}{Crit}

\DeclareMathOperator{\CF}{CF}

\DeclareMathOperator{\grad}{grad}

\newcommand{\Fsig}{F_\Sigma}
\newcommand{\Bpar}{\mathbb B}

\title{Fast growth of periodic orbits, translated points and short Floer bars for contactomorphisms}
\date{}
\author{
    Wenmin Gong \thanks{Laboratory of Mathematics and Complex Systems (Ministry of Education), School of Mathematical Sciences, Beijing Normal University, Beijing, 100875, China. Email: wmgong@bnu.edu.cn}
    \and
    Zhiyuan Hu \thanks{Department of Mathematical Sciences, Tsinghua University, Beijing, 100084, China. Email:\\ huzy23@mails.edu.tsinghua.cn}
}

\begin{document}
\maketitle
\begin{abstract}
    We establish arbitrarily fast growth of distinct periodic orbits and geometrically distinct translated points for contactomorphisms of the unit cotangent bundle of a closed hyperbolic surface. The conclusions hold on residual subsets of a nonempty $C^1$-open set and follow from Chaidez's contact blender and localized contact perturbations. For a specified analytically admissible normalized Rabinowitz Floer family on a closed contact manifold with Liouville filling, we bound fixed-window barcode entropy by topological entropy. The proof compares long-bar counts with graph volume through contact tomography and filtered continuation. Under strong translated nondegeneracy, we show that super-exponential translated-pair growth in a compact interior window forces super-exponential growth of short bars born in the surrounding window.
\end{abstract}
\tableofcontents

\section{Introduction}

    Counting periodic points is a fundamental way to quantify recurrence in dynamical systems. For Axiom~A diffeomorphisms, exponential counting bounds and their relation to topological entropy make exponential growth a natural benchmark~\cite{Bowen1970Markov,Bowen1970Entropy}. Beyond uniform hyperbolicity, much faster growth can occur. Arnold's question about exceeding a prescribed sequence along a subsequence of iterates has motivated a substantial literature on this phenomenon~\cite[Problem~1989-2]{Arnold2005}.

    Contact dynamics brings this question into a geometric setting in which diffeomorphisms preserve a contact hyperplane distribution. Alongside periodic points, translated points describe returns to Reeb trajectories with conformal factor one. Introduced by Sandon~\cite{Sandon2012,Sandon2013}, these points connect contact rigidity and fixed-point questions with Rabinowitz Floer theory~\cite{AlbersMerry2013}. Their abundance under iteration therefore raises both dynamical and homological questions.

    We establish arbitrarily fast growth of distinct periodic orbits in a nonempty $C^1$-open set of contactomorphisms of the unit cotangent bundle of a closed hyperbolic surface. For each prescribed positive sequence, a residual subset has orbit counts exceeding that sequence by arbitrarily large factors along an unbounded sequence of least periods. We also prove arbitrarily fast growth of geometrically distinct translated points in every fixed shift window of positive radius for a suitable contact form. Under the Floer hypotheses stated below, an entropy bound controls the number of sufficiently long bars born in a fixed action window. With the additional nondegeneracy and translated-point growth hypotheses of the short-bar theorem, the abundance of translated points forces super-exponential growth of short bars.

    \medskip\noindent\textbf{Fast growth of periodic points and the blender mechanism}

    Let $(S,h)$ be a closed oriented hyperbolic surface of constant curvature $-1$, and write
    \[
        M=S^*S:=\{(q,p)\in T^*S:|p|_h=1\},\qquad
        \alpha_{\mathrm{can}}:=\lambda_{\mathrm{can}}|_{S^*S},\qquad
        \xi:=\ker\alpha_{\mathrm{can}},
    \]
    where $\lambda_{\mathrm{can}}$ is the Liouville one-form on $T^*S$. The manifold $S^*S$ is the \emph{unit cotangent bundle}. Metric duality $(q,v)\mapsto(q,h_q(v,\cdot))$ identifies it with the unit tangent bundle $T^1S$, under which the Reeb flow of $\alpha_{\mathrm{can}}$ becomes the unit-speed geodesic flow~\cite[Section~1.5]{Geiges2008}. We write $\Cont^1(M,\xi)$ for the space of \emph{coorientation-preserving} $C^1$ contactomorphisms with its $C^1$ topology. Coorientation preservation means that the pullback of a defining contact form is a positive multiple of that form. A subset is \emph{residual} if it contains a countable intersection of open dense subsets, and residual subsets are dense in the Baire spaces considered here~\cite{Oxtoby1980}. We use \emph{generic} in this residual sense.

    The \emph{least period} of a periodic point is the smallest positive integer $n$ for which $f^n(x)=x$. Let $O_n(f)$ denote the number of distinct periodic orbits of least period $n$, so that each orbit is counted once. Much of the earlier literature counts fixed points of iterates, whereas the conclusion below concerns distinct orbits with a specified least period.

    Artin and Mazur proved the density of maps whose isolated periodic-point counts grow at most exponentially~\cite[Theorem~2.2]{ArtinMazur1965}. However the growth among a residual subset can behave differently. Bonatti, D\'iaz and Fisher established arbitrarily fast growth inside homoclinic classes containing saddles of different indices for $C^1$-generic diffeomorphisms~\cite{BonattiDiazFisher2008}, while Asaoka, Shinohara and Turaev obtained higher-regularity results from partially hyperbolic dynamics with one-dimensional center and suitable heterodimensional transitions~\cite{AsaokaShinoharaTuraev2021}. Such growth also occurs in geometrically constrained classes. Asaoka establishes arbitrarily fast growth for real-analytic Hamiltonian maps of the two-torus and smooth area-preserving surface maps \cite{Asaoka2017}.

    Our first theorem realizes arbitrarily fast growth in the contact category, as illustrated by the unit cotangent bundle introduced above.

    \begin{theorem}\label{thm:intro-periodic}
        Let \(M=S^*S\) be the unit cotangent bundle of a closed hyperbolic surface \(S\), equipped with its canonical contact structure \(\xi\), as above. There exists a nonempty \(C^1\)-open set \(\mca{U}\subset\Cont^1(M,\xi)\) such that, for every sequence \((a_n)_{n\ge1}\) of positive real numbers, there is a residual subset \(\mca{R}_a\subset\mca{U}\) satisfying
        \[
            \limsup_{n\to\infty}\frac{O_n(f)}{a_n}=\infty
            \qquad\text{for every }f\in\mca{R}_a.
        \]
    \end{theorem}

    The residual set may depend on the comparison sequence. Theorem~\ref{thm:arbitrarily-fast-contact-growth} gives the same conclusion on any closed contact manifold admitting an Anosov Reeb flow with smooth strong stable and unstable foliations.

    The proof uses the \emph{blender mechanism}. Introduced by Bonatti and D\'iaz in the study of robustly transitive nonhyperbolic diffeomorphisms~\cite{BonattiDiaz1996}, a blender makes the closure of a stable manifold robustly meet a suitable open family of disks even when the dimensions involved would ordinarily prevent transversality. This intersection property brings invariant manifolds of saddles of different indices close to one another and plays a central role in constructions of robust heterodimensional cycles~\cite{BonattiDiaz2008,BonattiDiazViana2005}. Chaidez's contact blender~\cite{Chaidez2025} makes this geometry available for contactomorphisms. In our setting, a localized contact perturbation turns the robust proximity into a strong heteroclinic connection, while a transverse connection in the opposite direction completes a co-index-one cycle, whose saddles have unstable dimensions differing by one. Smooth maps carrying this cycle are called \emph{based maps}, and Proposition~\ref{prop:contact heterodim cycle} proves their $C^1$-density in the chosen open set. The precise blender and cycle conventions appear in Section~\ref{subsection:heterodimensional cycles and blenders} and Definition~\ref{def:based-maps}.

    \medskip\noindent\textbf{Translated points and their growth}

    Fix any smooth contact form $\alpha$ defining $\xi$. Its \emph{Reeb vector field} $R_\alpha$ is uniquely determined by $\alpha(R_\alpha)=1$ and $\iota_{R_\alpha}\D\alpha=0$, and we denote its flow by $R_s^\alpha$.  For $f\in\Cont^1(M,\xi)$, let $\kappa_f:M\to(0,\infty)$ be the continuous \emph{conformal factor} defined by $f^*\alpha=\kappa_f\alpha$. A \emph{translated-point-shift} pair for $f^n$ is a solution $(x,s)\in M\times\mathbb{R}$ of
    \begin{equation}\label{eq:TP}
        f^n(x)=R_s^\alpha(x),\qquad \kappa_{f^n}(x)=1.
    \end{equation}
    The underlying \emph{translated point} is $x$, while $s$ is its \emph{time shift}. A point on a closed Reeb orbit may admit several shifts. We therefore distinguish the number of translated-point-shift pairs from the number of their geometrically distinct underlying points. For $s_*>0$ define
    \begin{equation}\label{eq:intro-TP-count}
        \TP_n(f;s_*):= \#\bigl\{(x,s)\in M\times(-s_*,s_*): f^n(x)=R_s^\alpha(x),\ \kappa_{f^n}(x)=1\bigr\} \in\mathbb{N}_0\cup\{+\infty\}.
    \end{equation}
    The number of \emph{geometrically distinct translated points} in the same shift window is
    \begin{equation}\label{eq:intro-N-count}
        N_n(f;s_*):=\#\bigl\{x\in M: \kappa_{f^n}(x)=1,\ \exists s\in(-s_*,s_*) \text{ with } f^n(x)=R_s^\alpha(x)\bigr\} \in\mathbb{N}_0\cup\{+\infty\}.
    \end{equation}
    Thus $N_n$ counts each point once, irrespective of the number of its admissible shifts, whereas $\TP_n$ counts all pairs. Both counts refer to the $n$th iterate without a least-period condition, and the form $\alpha$ is fixed throughout. See also the translated-point conventions of \cite[Section~3]{AlbersMerry2013} and \cite{Sandon2012}.

    \begin{theorem}\label{thm:intro-translated}
        Let \(M=S^*S\) be the unit cotangent bundle of a closed hyperbolic surface \(S\), equipped with its canonical contact structure \(\xi\), as above. There exist a smooth function $u:M\to (0,\infty)$ and a nonempty $C^1$-open set $\mca{U}\subset\Cont^1(M,\xi)$ such that, for the contact form $\alpha=u\alpha_{\mathrm{can}}$, the following holds. Fix $s_*>0$. For every sequence $(a_n)_{n\ge1}$ of positive integers, there is a residual set $\mca{R}_a\subset\mca{U}$ satisfying
        \[
            \limsup_{n\to\infty}\frac{N_n(f;s_*)}{a_n}
            =\limsup_{n\to\infty}\frac{\TP_n(f;s_*)}{a_n}
            =+\infty
            \qquad(f\in\mca{R}_a).
        \]
        In particular, there is a residual subset $\mca{R}_\infty\subset\mca{U}$ such that
        \[
            \limsup_{n\to\infty}\frac{1}{n}
            \log\max\{1,N_n(f;s_*)\}
            =\limsup_{n\to\infty}\frac{1}{n}
            \log\max\{1,\TP_n(f;s_*)\}
            =+\infty
            \qquad(f\in\mca{R}_\infty).
        \]
    \end{theorem}

    Theorems~\ref{thm:intro-periodic} and~\ref{thm:intro-translated} use the same open set from Proposition~\ref{prop:persistent-contact-data}, in which the residual sets may depend on the sequence and shift window. The proof creates isolated translated pairs with nonzero local \emph{Brouwer degree}, the integer invariant detecting zeros of a local map \cite[Chapter~1]{Deimling1985}. Its homotopy invariance yields robust local counts. Lemma~\ref{lemma:uniform-shift-multiplicity} then bounds shift multiplicity and transfers growth to distinct points.

    For Reeb flows, C\^ot\'e's unpublished note \cite[Theorem~1]{CoteReeb} describes contact forms on every closed connected contact manifold of dimension at least three whose numbers of closed Reeb orbits of period at most $k$ exceed any prescribed increasing sequence. This cumulative flow count differs from the least-period count in Theorem~\ref{thm:intro-periodic} and the count in a fixed Reeb-shift window in Theorem~\ref{thm:intro-translated}.

    \medskip\noindent\textbf{Filtered Floer theory and short bars}

    The geometric count above can grow much faster than the part of a Floer complex detected at a fixed action scale. We make this distinction using \emph{barcodes}, multisets of intervals describing the birth and death of homology classes in an increasing filtration; see \cite{PolterovichShelukhin2016,ChazalDeSilvaGlisseOudot2016,UsherZhang2016,CineliGinzburgGurel2024Entropy,CineliGinzburgGurel2024Growth}. Our action filtration is a normalized Rabinowitz Floer filtration whose critical points are translated-point-shift pairs and whose critical action equals the Reeb shift. The construction uses the perturbed Rabinowitz functional of Albers--Frauenfelder and its contact interpretation by Albers--Merry \cite{AlbersFrauenfelder2010Leafwise,AlbersMerry2013}.

    Let $(W,\lambda)$ be a \emph{Liouville domain}, meaning that $W$ is compact, $\D\lambda$ is symplectic, and the vector field $Z$ defined by $\iota_Z\D\lambda=\lambda$ points outward along $Y=\partial W$. The boundary contact form is $\alpha=\lambda|_Y$; see Definition~\ref{def:liouville} and \cite{McDuffSalamon1998}. Write $\Cont_0^\infty(Y,\ker\alpha)$ for the smooth contactomorphisms joined to the identity by a smooth contact isotopy. For $f\in\Cont_0^\infty(Y,\ker\alpha)$, put
    \[
        g=\log\kappa_f,\qquad G_n=\sum_{j=0}^{n-1}g\circ f^j,
        \qquad \Gamma_n=\{(x,f^n(x),G_n(x)):x\in Y\}.
    \]
    Volume on this graph is measured with a fixed product Riemannian metric.

    We use an \emph{analytically admissible family} in the precise sense of Definition~\ref{def:analytically-admissible-family}. This means that the specified, separated-time Rabinowitz functionals have the regularizations, finite-window complexes and continuation theory required there. This analytic hypothesis concerns the filtered construction, while barcode finiteness is proved subsequently. The associated barcode $\mca{B}_n$ is defined in Definition~\ref{def:associated-finite-birth-barcode} from a completion of the action-filtered complex and, for degenerate members, by limits in the \emph{bottleneck distance}, which minimizes the largest endpoint displacement, allowing an unmatched bar at a cost of half its length \cite[Section~4.2]{ChazalDeSilvaGlisseOudot2016}. It includes bars $[b,d)$ with $b\in\R$ and $b<d\le+\infty$. All relevant free homotopy classes are included. For a compact interval $I\subset\R$ and $\varepsilon>0$, set
    \[
        b_\varepsilon^I(\mca{B}_n)
        =\#\{[b,d)\in\mca{B}_n:b\in I,\ d-b\ge\varepsilon\}.
    \]
    Bars are counted with multiplicity, and an infinite bar has infinite length. The \emph{fixed-window barcode entropy} is
    \[
        h_{\mathrm{TP}}^I(f;\alpha;\{K_{n,\cdot}\})
        =\lim_{\varepsilon\downarrow0}\limsup_{n\to\infty}
        \frac{1}{n}\log\max\{1,b_\varepsilon^I(\mca{B}_n)\}.
    \]
    This is the action-window version of the barcode-growth convention in \cite{CineliGinzburgGurel2024Entropy}, and the Hamiltonian choices $K_{n,\cdot}$ form part of the filtered construction. We write $h_{\mathrm{top}}(f)$ for topological entropy, defined by the exponential growth of orbit segments distinguishable at a fixed spatial scale, followed by the limit as that scale tends to zero \cite[Chapter~3]{KatokHasselblatt1995}.

    \begin{theorem}\label{thm:intro-barcode}
        Let $(W,\lambda)$ be a Liouville domain with boundary $(Y,\alpha)$, let $f\in\Cont_0^\infty(Y,\ker\alpha)$, and choose an analytically admissible normalized family as in Definition~\ref{def:analytically-admissible-family}. For every compact Reeb-shift interval $I$ and every $\varepsilon>0$, there is a constant $C_{I,\varepsilon}$ independent of $n$ such that
        \[
            b_\varepsilon^I(\mca{B}_n)
            \le C_{I,\varepsilon}\operatorname{Vol}_{\dim Y}(\Gamma_n).
        \]
        Consequently,
        \[
            h_{\mathrm{TP}}^I(f;\alpha;\{K_{n,\cdot}\})\le h_{\mathrm{top}}(f).
        \]
    \end{theorem}

    The estimate combines a smooth circle compactification of $\Gamma_n$ with a \emph{contact tomograph}, a finite-dimensional family that makes parameterized graph intersections transverse and bounds their mean number by graph volume. The tomograph adapts the Lagrangian method of~\cite[Section~5.2]{CineliGinzburgGurel2024Entropy}. Barcode stability then transfers this intersection bound to bars of a fixed positive minimum length, while the graph-volume estimate relates their exponential growth to topological entropy.

    A smooth map is called \emph{strongly translated-nondegenerate} if every translated-point-shift pair of every positive iterate is nondegenerate, meaning that the map
    \[
        (x,s)\longmapsto(f^n(x),R_s^\alpha(x),G_n(x))
    \]
    is transverse to $\Delta_Y\times\{0\}$ at every solution; see Definition~\ref{def:regular}, Proposition~\ref{prop:normalized-hessian-translated-pair-kernel}, and \cite{AlbersMerry2013}. For a bounded open interval $I$, define
    \[
        S_\varepsilon^I(n)
        =\#\{[b,d)\in\mca{B}_n:b\in I,\ 0<d-b<\varepsilon\}.
    \]
    We call these the \emph{short bars born in $I$}. When we instead allow either finite endpoint to lie in $I$, we write $S_{\varepsilon,\mathrm{end}}^I(n)$.

    \begin{theorem}
        \label{thm:intro-short}
        Assume the hypotheses of Theorem~\ref{thm:intro-barcode}, and suppose that $f$ is strongly translated-nondegenerate. Let $I$ be a bounded open interval. If some compact interval $J\subset I$ satisfies
        \[
            \limsup_{n\to\infty}\frac{1}{n}
            \log\max\{1,\#\TP(f^n;J)\}=+\infty,
        \]
        then, for every fixed $\varepsilon>0$,
        \[
            \limsup_{n\to\infty}\frac{1}{n}
            \log\max\{1,S_\varepsilon^I(n)\}=+\infty.
        \]
        If the growth assumption is known only in $I$, the same conclusion holds for $S_{\varepsilon,\mathrm{end}}^I(n)$.

        In the setting of Theorem~\ref{thm:intro-translated}, one can choose a $C^1$-residual set $\widetilde{\mca{R}}_\infty\subset\mca{U}$ so that the conclusion for $I=(-s_*,s_*)$ holds for every $f\in\widetilde{\mca{R}}_\infty$ satisfying the smoothness, isotopy, nondegeneracy and admissibility hypotheses above.
    \end{theorem}

    The proof counts barcode endpoints and uses the entropy bound to control the contribution of long bars. Since $J$ is compactly contained in $I$, sufficiently short bars with an endpoint in $J$ are born in $I$.

    \medskip\noindent\textbf{Outlines}

    Section~\ref{section:preliminary} establishes the contact, dynamical and Floer conventions used throughout the paper. The local perturbation and approximation lemmas of Section~\ref{section:Contact-lemmas} support the construction of a persistent blender configuration and a $C^1$-dense family of smooth based maps in Section~\ref{section:based-maps}. Starting from these maps, Section~\ref{section:Long flat orbit near the cycle} produces long periodic orbits with flat center germs. These orbits provide the common input for the periodic-orbit growth theorem in Section~\ref{section:superexponential} and the translated-point growth theorem in Section~\ref{section:translated-point-proliferation}.

    The Floer part begins in Section~\ref{section:contact-graphs-and-barcode-entropy}, where graph-volume estimates and the normalized filtered construction lead to the fixed-window barcode entropy bound under the stated Floer hypotheses. With this bound in place, Section~\ref{section:superexponential-short-bar-growth} derives the short-bar growth theorem and its prescribed-growth variant under their additional hypotheses. Appendix~\ref{app:hamiltonian-lifts-rabinowitz-functional} verifies the Hamiltonian normalization and specifies the analytic input for continuation, while Appendix~\ref{app:completed-filtrations-barcode-endpoints} establishes the completion and endpoint-counting results used in the Floer argument.
\subsection*{Acknowledgements}
    Both authors are very grateful to Prof. Jinxin Xue and Prof. Jinhua Zhang for their generous help and valuable comments. We also thank Xue's grant supports NSFC 11790273 and 12271285, and the New Cornerstone investigator program.
\subsection*{AI disclosure}
    We developed and refined the major framework of the work. With ChatGPT's assistance, we established a suitable form of the contact approximation lemma, simplified the strategy for Section~\ref{section:Long flat orbit near the cycle}, and constructed the submersion in Lemma~\ref{lemma:contact-tomograph-submersion} to obtain the desired contact tomograph. This tomograph was inspired by \c{C}ineli--Ginzburg--G\"{u}rel's Lagrangian tomograph. We also thank ChatGPT for its support in concrete computations, the preparation of figures, and manuscript review.
\section{Preliminaries}\label{section:preliminary}
\subsection{Contactomorphisms, Reeb flows and contact Hamiltonian flows}
    We recall the standard contact-geometric conventions; see \cite[Sections~1.1 and 2.3]{Geiges2008} or \cite{Geiges2006}. A \emph{contact form} on a $(2d+1)$-manifold $M$ is a smooth one-form $\alpha$ for which $\alpha\wedge(\D\alpha)^d$ is nowhere zero. The \emph{contact structure}, or the \emph{contact distribution}, is its maximally nonintegrable hyperplane distribution $\xi=\ker\alpha$. Its coorientation is specified by choosing the positive multiples of $\alpha$. A $d$-dimensional submanifold tangent to $\xi$ is called \emph{Legendrian}; see \cite[Section~2.3]{Geiges2008}. The Reeb vector field $R_\alpha$ is defined by 
    \[
        \alpha(R_\alpha)=1,\qquad \iota_{R_\alpha}\D\alpha=0.
    \]
    A $C^r$ diffeomorphism $f$ is a \emph{coorientation-preserving contactomorphism} if $f^*\alpha=\kappa_f\alpha$ for a positive function $\kappa_f$, called its \emph{conformal factor}. It is \emph{strict} with respect to $\alpha$ if $\kappa_f\equiv1$. We write $\Cont^r(M,\xi)$ for the space of coorientation-preserving $C^r$ contactomorphisms, and $\Cont^r(M,\alpha)$ for the space of strict $C^r$ contactomorphisms. For a $C^1$ contactomorphism, we initially assume only that the conformal factor is continuous and establish its additional regularity in Section~\ref{subsection:contact-approximation}.

    Given a smooth function $H$ on $M$, the associated \emph{contact Hamiltonian vector field} $X_H$ is defined (with our sign convention) by
    \begin{equation}\label{eq:contact Hamiltonian vf.}
        \alpha(X_H)=H,\qquad \iota_{X_H}\D\alpha=(R_\alpha H)\alpha-\D H
    \end{equation}
    where \(R_\alpha\) is the Reeb vector field of $\alpha$.

    In Darboux coordinates, the contact form is $\alpha_0=\D{z}-x\cdot\D{y}$ and its Reeb field is $R_{\alpha_0}=\partial_z$ \cite[Theorem~2.5.1]{Geiges2008}. Substituting this form in \eqref{eq:contact Hamiltonian vf.} gives
    \begin{equation}\label{eq:Hamiltonian eq. in xyz}
        \begin{cases}
            \dot{x} = (\pps{y}H+x\pps{z}H),\\
            \dot{y} = -\pps{x}H,\\
            \dot{z} = H-x\pps{x}H .
        \end{cases}
    \end{equation}
\subsection{Partial hyperbolicity}\label{subsection:partial-hyperbolicity}
    We follow the conventions for dominated splittings and invariant manifolds in \cite{BonattiDiazViana2005,HirschPughShub1977}. Fix a Riemannian metric on $M$, let $f\in\Diff^1(M)$, and let $\Lambda$ be a compact $f$-invariant set. A continuous $Df$-invariant splitting $T_\Lambda M=E_1\oplus\cdots\oplus E_k$ into nonzero bundles of constant dimensions is \emph{dominated} over $\Lambda$ if there exist $N\ge1$ and $0<\theta<1$ such that
    \[
        \|Df^N|_{E_i(x)}\|\cdot\|(Df^N|_{E_{i+1}(x)})\inv\|\le\theta\qquad(x\in\Lambda,\ 1\le i<k).
    \]
    We denote this ordering by $E_1\prec\cdots\prec E_k$.

    The set $\Lambda$ is \emph{partially hyperbolic} if it admits such a dominated splitting
    \[
        T_\Lambda M=E^{ss}\oplus E^c\oplus E^{uu},\qquad E^{ss}\prec E^c\prec E^{uu},
    \]
    whose extreme bundles are uniformly contracting and expanding, respectively. We write $d_s=\dim E^{ss}$, $d_c=\dim E^c$, and $d_u=\dim E^{uu}$. The corresponding strong invariant manifolds are denoted by $W^{ss},W^{uu}$, and the subscript $\loc$ indicates a local plaque. At a hyperbolic periodic point, $W^s,W^u$ denote the full stable and unstable manifolds, which may also contain the center direction.

    A basic source is an \emph{Anosov flow} \cite[Chapter~17]{KatokHasselblatt1995}. A flow $\varphi^t$ on a closed manifold, generated by a nowhere-vanishing $C^1$ vector field $X$, is \emph{Anosov} if it admits a continuous $D\varphi^t$-invariant splitting $TM=E^{ss}\oplus\R X\oplus E^{uu}$, with nonzero strong bundles and constants $C\ge1$, $\lambda>0$ such that $\|D\varphi^t|_{E^{ss}(x)}\|\le Ce^{-\lambda t}$ and $\|D\varphi^{-t}|_{E^{uu}(x)}\|\le Ce^{-\lambda t}$ for all $x\in M$ and $t\ge0$. Since $D\varphi^t_xX(x)=X(\varphi^t(x))$ and $X$ is bounded above and away from zero, the derivative along the flow direction is uniformly bounded for positive and negative times. Consequently, the time-one map is partially hyperbolic with $E^c=\R X$.
\subsection{Heterodimensional cycles and blenders}\label{subsection:heterodimensional cycles and blenders}
    Let $f$ be a partially hyperbolic diffeomorphism of $M^d$ with splitting $TM=E^{ss}\oplus E^c\oplus E^{uu}$. Throughout this subsection, write $\bs d=(d_s,d_c,d_u)$. Thus $d=d_c+d_s+d_u$, and we use unstable indices for saddles. We first recall the topological property of a blender and explain its role in producing heteroclinic connections. We then describe the geometric construction that ensures this property.

    \medskip\noindent{\bf Heterodimensional cycles}

    A heterodimensional cycle consists of two transitive hyperbolic sets $\Lambda_1,\Lambda_2$ of different indices whose invariant manifolds meet in both directions, namely $\lsp{W}{s}{\Lambda_1}\cap\lsp{W}{u}{\Lambda_2}\neq\varnothing$ and $\lsp{W}{s}{\Lambda_2}\cap\lsp{W}{u}{\Lambda_1}\neq\varnothing$. For cycles associated to periodic saddles, we use the following configuration of transverse and quasi-transverse intersections, cf.~\cite[p.~192]{BonattiDiazPujalsRocha2003}.
    \begin{definition}
        A \emph{heterodimensional cycle associated to a pair of saddles}, or \emph{saddle cycle}, consists of hyperbolic periodic points $P,Q$ on distinct orbits such that
        \begin{itemize}
            \item $\dim W^s(P)=d_c+d_s$ and $\dim W^u(P)=d_u$;
            \item $\dim W^s(Q)=d_s$ and $\dim W^u(Q)=d_c+d_u$;
            \item $W^s(P)$ and $W^u(Q)$ have a nonempty transverse intersection;
            \item $W^u(P)$ and $W^s(Q)$ have a nonempty quasi-transverse intersection $a$, meaning that $T_aW^u(P)\cap T_aW^s(Q)=\{0\}$ at a chosen heteroclinic point $a$.
        \end{itemize}
    \end{definition}

    \begin{figure}[H]
        \centering
        \includegraphics[scale=0.65]{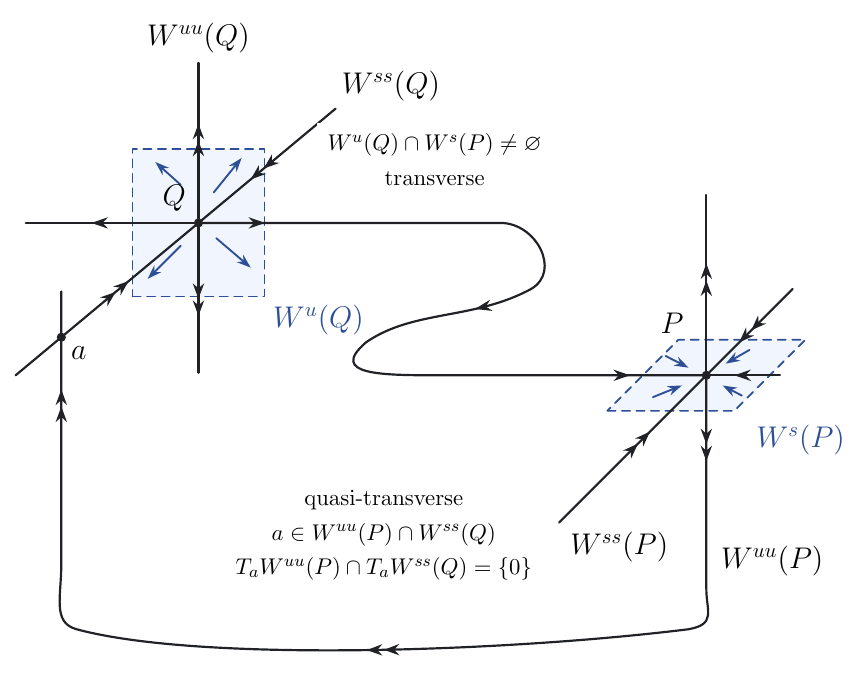}
        \caption[Heterodimensional saddle cycle]{\textbf{Heterodimensional saddle cycle.}
            In the schematic, the connection from $Q$ to $P$ is transverse, whereas $a\in W^{uu}(P)\cap W^{ss}(Q)$ is quasi-transverse.}
        \label{fig:heterodimensional-cycle}
    \end{figure}

    In this configuration the \emph{co-index}, namely the difference between the unstable indices, is $d_c$. In the remainder of the subsection we take 
    $$
    d_c=1
    $$
    and, after passing to an iterate when necessary, assume that $P,Q$ are fixed. Thus $P$ is center-contracting and $Q$ is center-expanding, with $W^u(P)=W^{uu}(P)$ and $W^s(Q)=W^{ss}(Q)$. Figure~\ref{fig:heterodimensional-cycle} distinguishes the transverse connection from $Q$ to $P$ from the quasi-transverse strong connection through $a$, whose lack of $C^1$-robustness motivates the blender mechanism.

    For $k\ge1$, let $D^k$ be a closed unit disk and let $\mca{D}^k(M)$ be the space of its $C^1$ embedded images, with the $C^1$ topology. Fix a hyperbolic saddle $Q$ with $\ind_u(Q)=d_u+1$.

    \begin{definition}[\cite{BonattiDiazViana2005}, Definition~6.11]\label{def:topological blender}
        We say that $f$ has a \emph{topological blender $\mathfrak{B}$ associated to $Q$} if there are a $C^1$-neighborhood $\mca{U}$ of $f$ and a nonempty $C^1$-open family $\mca{R}\subset\mca{D}^{d_u}(M)$ such that every $G\in\mca{U}$ and every disk $D\in\mca{R}$ satisfy $D\cap\overline{W^s(Q_G,G)}\ne\varnothing$, where $Q_G$ is the unique nearby hyperbolic periodic point of the same period (the \emph{hyperbolic continuation} of $Q$).
    \end{definition}
    The family $\mca{R}$ is the \emph{superposition region} of $\mathfrak{B}=(f,Q,\mca{U},\mca{R})$. Although $\dim W^s(Q)=d_s$, its closure robustly meets an open family of $d_u$-disks, behaving in this respect as though it had an additional dimension. When $W^{uu}(P)$ contains a disk in $\mca{R}$, points of $W^{ss}(Q)$ approach that disk arbitrarily closely, and this proximity persists under $C^1$ perturbations. The \emph{nonempty} strong connection $W^{uu}(P)\cap W^{ss}(Q)$ is necessarily non-transverse in $M$, since the dimensions of these manifolds sum to $d-1$. In Section~\ref{subsection:contact-blender-connection}, orbit selection and a local contact perturbation turn the robust proximity into such a connection.

    We call this a \emph{center-unstable} blender because the reference saddle has an expanding center. Bonatti--D\'iaz use the name \emph{cs-blender} in \cite{BonattiDiaz1996}, emphasizing its stable-manifold property. The geometric criterion below follows their construction.

    \medskip\noindent{\bf Conefields}

    Let $L^s\oplus L^c\oplus L^u$ be a continuous splitting, not necessarily $Df$-invariant, defined over a coordinate neighborhood, with dimensions $(d_s,1,d_u)$. Write $v=v^s+v^c+v^u$ according to this splitting. For $0<\varepsilon<1$, define its \emph{associated $\varepsilon$-conefields} by
    \[
        \mca{C}^{ss}(p)=\{v:\|v^c\|+\|v^u\|\le\varepsilon\|v^s\|\},\qquad
        \mca{C}^{cs}(p)=\{v:\|v^u\|\le\varepsilon(\|v^c\|+\|v^s\|)\}.
    \]
    The conefields $\mca{C}^{uu}(p)$ and $\mca{C}^{cu}$ are analogously defined. Put $\mca{C}^c=\mca{C}^{cs}\cap\mca{C}^{cu}$. When the splitting is the partially hyperbolic splitting, i.e., $L^s=E^{ss},L^c=E^c$ and $L^u=E^{uu}$, these conefields are called the \emph{partially hyperbolic $\varepsilon$-conefields}.

    Let $F=f^r$ be a return branch on a set $K$, with $K$ and $F(K)$ in the coordinate neighborhood. Following condition (H3) of \cite{BonattiDiaz1996}, we require forward invariance of the expanding cones and backward invariance of the stable cones. Explicitly, for every $p\in K$,
    \[
    \begin{aligned}
        DF|_p\bigl(\mca{C}^{cu}(p)\setminus\{0\}\bigr)&\subset\interior\mca{C}^{cu}(F(p)),\\
        DF|_p\bigl(\mca{C}^{uu}(p)\setminus\{0\}\bigr)&\subset\interior\mca{C}^{uu}(F(p)),\\
        D(F^{-1})|_{F(p)}\bigl(\mca{C}^{ss}(F(p))\setminus\{0\}\bigr)&\subset\interior\mca{C}^{ss}(p).
    \end{aligned}
    \]
    Here interiors are taken in the tangent fibers. We also require a constant $\lambda>1$ such that
    \[
        \|DF|_pv\|\ge\lambda\|v\|\quad(v\in\mca{C}^{cu}(p)),\qquad
        \|D(F^{-1})|_{F(p)}w\|\ge\lambda\|w\|\quad(w\in\mca{C}^{ss}(F(p))).
    \]
    These will be called the \emph{center-unstable cone conditions} for the branch $F$. The \emph{center-stable cone conditions} are obtained by reversing time and exchanging the strong directions. In particular, forward center-unstable and backward center-stable cone invariance can coexist.

    \medskip\noindent{\bf Local charts and component conditions}

\begin{figure}[H]
        \centering
        \includegraphics[scale=0.75]{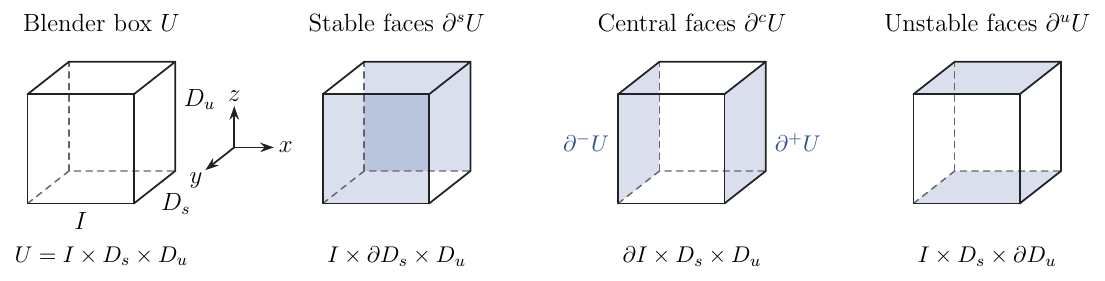}
        \caption[Blender box and its boundary faces]{\textbf{Blender box and its boundary faces.}
            The blender box $U=I\times D_s\times D_u$ has center, stable, and unstable directions labeled $x,y,z$, respectively. The last three cubes highlight the stable faces $\partial^sU$, central faces $\partial^cU$, and unstable faces $\partial^uU$.}
        \label{fig:blender-box-boundary}
    \end{figure}

    Choose coordinates $(x,y,z)\in\R\times\R^{d_s}\times\R^{d_u}$ in a neighborhood $U'$ of $Q$, where $x$ is the center coordinate, and consider the box
    \[
        U=I\times D_s\times D_u,\qquad
        I=[-2,2],\quad D_s=[-2,2]^{d_s},\quad D_u=[-2,2]^{d_u}.
    \]
    We allow boxes and disks with corners. The chart extends beyond $U$, and the conefields are those around its coordinate splitting $T\rn{d_s}\oplus T\R\oplus T\rn{d_u}$. Norms and lengths in this construction are measured with the product Euclidean metric of the chart. Define the stable, central, and unstable faces by
    \[
        \partial^sU=I\times\partial D_s\times D_u,\qquad
        \partial^cU=\partial I\times D_s\times D_u,\qquad
        \partial^uU=I\times D_s\times\partial D_u,
    \]
    and write $\partial^\pm U=\{\pm2\}\times D_s\times D_u$ for the right and left central faces, which are specified in Figure~\ref{fig:blender-box-boundary}.

    The primary branch has a compact domain $A\subset U$ such that $f(A)$ is a connected component of $f(U)\cap U$, contains $Q$, and satisfies
    \[
        f(A)\cap\bigl(\partial^sU\cup f(\partial^cU\cup\partial^uU)\bigr)=\varnothing.
    \]
    Together with the cone conditions, these conditions give the primary crossing geometry of \cite[Section~1]{BonattiDiaz1996}. The primary branch selects a local stable disk $W_0^s(Q)$ as the reference stable disk. This disk crosses $U$ in the stable direction, is tangent to $\mca{C}^{ss}$, and avoids $\partial^cU\cup\partial^uU$. We call such a box \emph{$Q$-compatible}.

    \begin{figure}[H]
        \centering
        \includegraphics[scale=0.5]{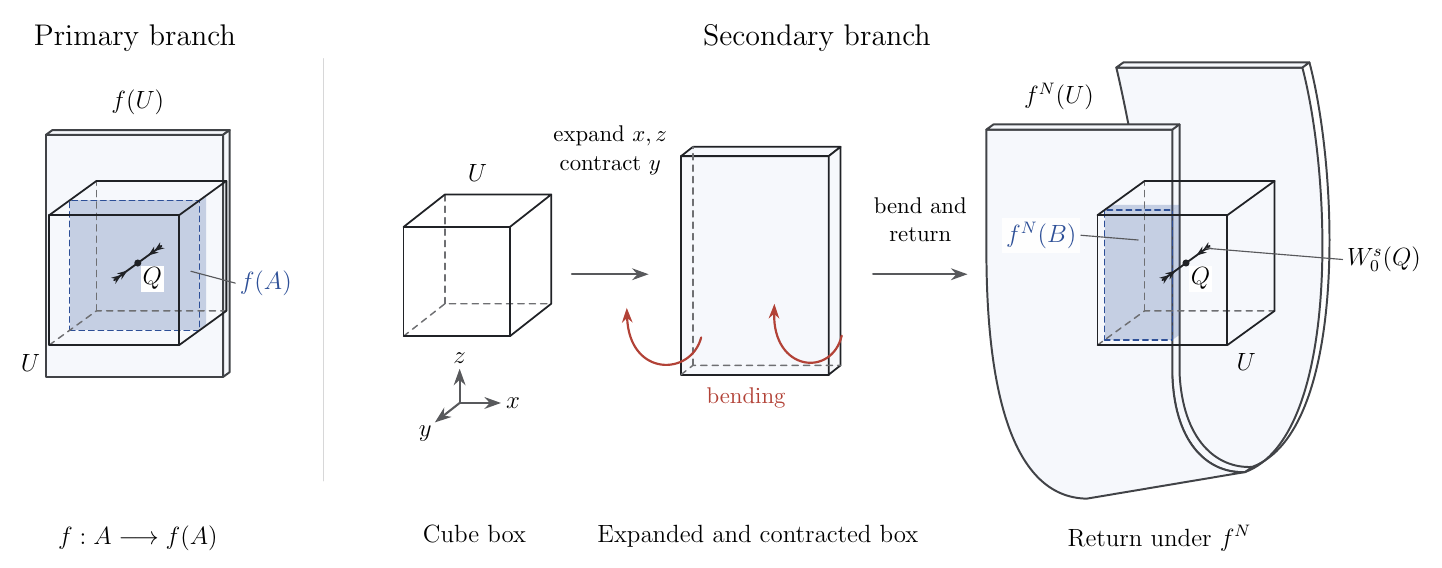}
        \caption[Primary and secondary blender branches]{\textbf{Primary and secondary blender branches.}
            The primary component $f(A)\subset f(U)\cap U$ contains $Q$, while the secondary component $f^N(B)\subset f^N(U)\cap U$ avoids $Q$ and meets $W_0^s(Q)$.}
        \label{fig:blender-branches}
    \end{figure}

    The secondary branch has a compact domain $B\subset U$ and a return time $N>1$ such that $f^N(B)$ is another connected component of $f^N(U)\cap U$, disjoint from $Q$ and meeting $W_0^s(Q)$. Its crossing conditions are
    \[
        f^N(B)\cap\bigl(\partial^sU\cup\partial^+U\cup f^N(\partial^uU)\bigr)=\varnothing.
    \]
    These boundary avoidance conditions prevent escape through the stable faces or the right central face, as well as clipping at an image of the unstable boundary. We call this a \emph{$W^s(Q)$-returning} branch. Here $N$ is a prescribed branch time, which need not be the first-return time, and the orbit may leave $U$ between its initial and final visits. The geometric construction assumes these crossing conditions in addition to the existence of a return.
    
    The dynamics of the branches are illustrated in Figure~\ref{fig:blender-branches}.

    \medskip\noindent{\bf The superposition region}

    A \emph{vertical disk through $U$} is a $d_u$-dimensional disk given as a graph over $D_u$, tangent to $\mca{C}^{uu}$, with boundary on $\partial^uU$ and disjoint from $\partial^sU\cup\partial^cU$. Let $\mca{R}'$ denote this family of crossing disks. We choose the cone size $\varepsilon$ small enough that the vertical disks disjoint from $W_0^s(Q)$ fall into two homotopy classes, represented by disks near $\partial^-U$ and $\partial^+U$. These homotopies pass through vertical disks avoiding $W_0^s(Q)$, with their boundaries remaining on $\partial^uU$. Such a disk is \emph{to the left} or \emph{to the right} according to its class, as in \cite[Section~1.1]{BonattiDiaz1996}.
    
    \begin{figure}[H]
        \centering
        \includegraphics[scale=0.6]{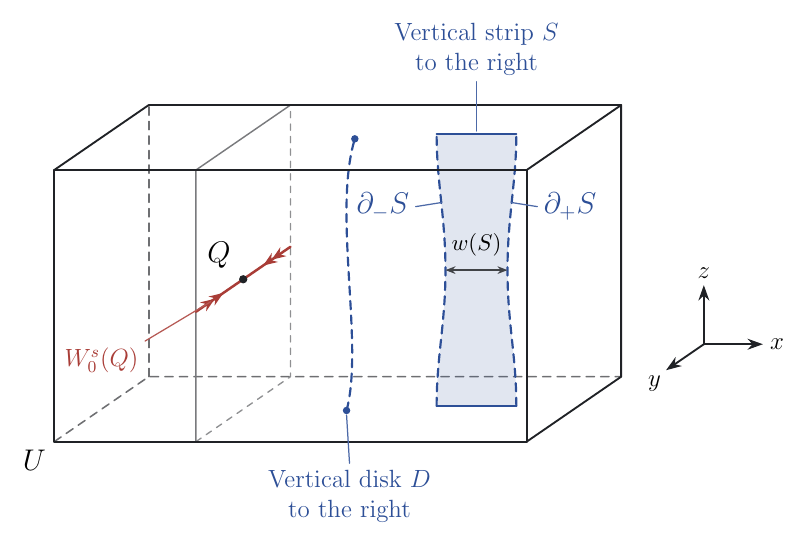}
        \caption[Vertical disks and strips]{\textbf{Vertical disks and strips.}
            The vertical disk $D$ and strip $S$ lie to the right of the reference stable disk $W_0^s(Q)$ and have their unstable boundaries on $\partial^uU$. The strip's vertical boundary components are $\partial_-S$ and $\partial_+S$, and $w(S)$ denotes its width.}
        \label{fig:superposition-region}
    \end{figure}

    A \emph{vertical strip through $U$} is an embedded image $S=\Psi([-1,1]\times D_u)\subset U$ tangent to $\mca{C}^{cu}$, such that each slice $\Psi(\{t\}\times D_u)$ is a vertical disk. Its two \emph{vertical boundary components} are $\partial_\pm S=\Psi(\{\pm1\}\times D_u)$, and the remaining boundary lies on $\partial^uU$. Write $\mca{S}$ for the family of these strips. A strip is \emph{to the right} if all its vertical slices are to the right. These definitions allow the disks and strips to curve within their cones, so they persist in the corresponding families after a small perturbation of the map $f$.

    The \emph{width} $w(S)$ is the infimum of the lengths of curves in $S$ joining $\partial_-S$ to $\partial_+S$. For a fixed box and sufficiently narrow coordinate cones, there is a uniform upper bound $\nu$ for the widths of all such strips, see \cite[Lemma~1.4]{BonattiDiaz1996}. Moreover, a vertical disk separated from the central faces can be thickened in the central direction to a strip. Taking thinner substrips gives strips converging to the disk in Hausdorff distance, as in \cite[Lemma~1.5]{BonattiDiaz1996}.

    \newpage
    \medskip\noindent{\bf Dynamical blenders}

    We now assemble these ingredients into a geometric criterion, following conditions (H1)--(H5) of \cite{BonattiDiaz1996} in the branchwise formulation above.
    \begin{definition}\label{def:dynamical blender}
        A \emph{center-unstable dynamical blender} is a tuple $\mca{B}=(f,Q,U,A,B,N,\varepsilon,\lambda)$ satisfying the following conditions.
        \begin{enumerate}[(H1)]
            \item The primary branch $f:A\to f(A)$ satisfies the $Q$-compatible crossing conditions above.
            \item The secondary branch $f^N:B\to f^N(B)$ satisfies the $W^s(Q)$-returning crossing conditions above.
            \item Both branches satisfy the center-unstable cone conditions, with the same cone size $\varepsilon$ and expansion constant $\lambda>1$.
            \item There is a neighborhood $U_-$ of $\partial^-U$ such that every vertical disk through $U$ to the right of $W_0^s(Q)$ is disjoint from $U_-$.
            \item There are neighborhoods $U_+$ of $\partial^+U$ and $U_0$ of $W_0^s(Q)$ such that every vertical disk $D\in\mca{R}'$ to the right satisfies at least one of the following alternatives:
            \begin{enumerate}[(i)]
                \item $f(D)\cap f(A)$ contains a vertical disk through $U$ to the right of $W_0^s(Q)$ and disjoint from $U_+$;
                \item $f^N(D)\cap f^N(B)$ contains a vertical disk through $U$ to the right of $W_0^s(Q)$ and disjoint from $U_0$.
            \end{enumerate}
        \end{enumerate}
        We call (H5) the \emph{component dichotomy}.
    \end{definition}

    To define a \emph{center-stable dynamical blender}, we apply this construction to $f^{-1}$ and exchange the roles of the strong stable and strong unstable directions. If $A,B$ are the domains of the corresponding forward branches $f|_A$ and $f^N|_B$, we call $(f,Q,U,A,B,N,\varepsilon,\lambda)$ a center-stable dynamical blender if the tuple for the inverse map $(f^{-1},Q,U,f(A),f^N(B),N,\varepsilon,\lambda)$ forms a center-unstable one.

    The remainder of this subsection concerns center-unstable blenders, whose center direction is expanding, as in Definition~\ref{def:topological blender}.

    The fixed neighborhoods in (H4)--(H5) provide the separation margins needed for persistence under perturbation. To see how these hypotheses act on a strip $S$, follow the images of its vertical slices. For each slice, the dichotomy provides a crossing component, and the boundary avoidance conditions allow the selected components to continue across the strip. If a secondary branch reaches the left central face, it must cross $W_0^s(Q)$, while the right central face cannot be reached. If no secondary right-hand component occurs, the primary branch provides the required family, whose uniqueness and continuous dependence are established in \cite[Lemma~1.6]{BonattiDiaz1996}.

    \begin{lemma}\label{lemma:superposition region}
        Every strip $S\in\mca{S}$ to the right of $W_0^s(Q)$ intersects $W^s(Q)$ transversely.
    \end{lemma}
    \begin{proof}
        By the component argument of \cite[Lemma~1.7]{BonattiDiaz1996}, one of the following two alternatives holds.
        \begin{enumerate}[(a)]
            \item either $f(S)\cup f^N(S)$ meets $W_0^s(Q)$, or
            \item one of these images contains a strip $S_1\in\mca{S}$ to the right whose two vertical boundary components are inherited from the images of the two vertical boundary components of $S$.
        \end{enumerate}
        In (b), every curve in $S_1$ joining its vertical boundary components pulls back under the selected branch to a curve in $S$ joining its vertical boundary components. Expansion on $\mca{C}^{cu}$ therefore gives
        \[
            w(S_1)\ge\lambda w(S).
        \]
        Suppose that $S$ is disjoint from $W^s(Q)$. Its forward images are then also disjoint from this invariant manifold, so repeated application of (b) produces strips $S_j$ with $w(S_j)\ge\lambda^j w(S)$. These widths contradict the uniform upper bound $\nu$. Thus some selected image meets $W_0^s(Q)$, and the complementary stable and center-unstable cones make the intersection transverse. Pulling this intersection back gives the required intersection with $W^s(Q)$, as in \cite[Lemma~1.8]{BonattiDiaz1996}.
    \end{proof}

    Fix now a vertical disk $D\in\mca{R}'$ to the right. By \cite[Lemma~1.5]{BonattiDiaz1996}, it admits arbitrarily thin strips $S_j\in\mca{S}$ around it, still to the right, with Hausdorff distance $d_H(S_j,D)\to0$. Choose $a_j\in S_j\cap W^s(Q)$ using Lemma~\ref{lemma:superposition region}. Compactness gives a subsequence converging to a point of $D\cap\overline{W^s(Q)}$. This is the passage from strips to disks in \cite[Lemma~1.9]{BonattiDiaz1996}.

    To obtain a superposition region open in the ambient space $\mca{D}^{d_u}(M)$, extend one such disk slightly beyond the unstable faces of $U$ while keeping positive margins from the other faces and from $W_0^s(Q)$. Every sufficiently $C^1$-close embedded disk then contains a vertical crossing subdisk of $U$ to the right. We may choose the extended disk with tangencies strictly inside the unstable cones and denote the resulting nonempty open family by $\mca{R}$. By \cite[Lemma~1.11]{BonattiDiaz1996}, the complete geometric conditions persist, including the fixed separation neighborhoods. After shrinking this family and the perturbation neighborhood if necessary, the same construction applies to every nearby map $G$, with reference disk $W_0^s(Q_G,G)$. Consequently, each $D\in\mca{R}$ meets $\overline{W^s(Q_G,G)}$, giving a topological blender in the sense of Definition~\ref{def:topological blender}.

    \begin{theorem}
        Every center-unstable dynamical blender is a center-unstable topological blender in the sense of Definition~\ref{def:topological blender}.
    \end{theorem}

    \begin{definition}
        A hyperbolic saddle $P$ is \emph{connected} to a blender $\mca{B}$ if $\ind_u(P)=d_u$ and $W^u(P)$ contains a disk $D$ in the superposition region $\mca{R}$ of $\mca{B}$.
    \end{definition}

    \begin{lemma}\label{lemma:blender saddle connection}
        If a saddle $P$ is connected to a blender associated to $Q$, then $W^s(P)\subset\overline{W^s(Q)}$.
    \end{lemma}
    \begin{proof}
        Let $p\in W^s(P)$ and let $V$ be an open neighborhood of $p$. Choose a small $d_u$-disk $D_1\subset V$ through $p$ and transverse to $W^s(P)$, and let $D_3\subset W^u(P)$ be a disk in $\mca{R}$. By the inclination lemma \cite{KatokHasselblatt1995}, for a suitable $m\ge1$, the iterate $f^m(D_1)$ contains a disk $D_2$ that is sufficiently $C^1$-close to $D_3$ to belong to the open family $\mca{R}$. Thus $D_2$ meets $\overline{W^s(Q)}$, and the open set $f^m(V)$ containing $D_2$ therefore meets $W^s(Q)$. Invariance gives $V\cap W^s(Q)\ne\varnothing$, proving the inclusion.
    \end{proof}

    The compact disk witnessing the connection continues under small perturbations of $f$, and the openness of $\mca{R}$ makes the connection robust. Lemma~\ref{lemma:blender saddle connection} therefore applies to the continuations. This is the mathematical essence of producing based maps in Section~\ref{section:based-maps}.
\subsection{Liouville fillings and the Floer conventions}\label{sec:floer-conventions}
    \begin{definition}[Liouville filling and completion]\label{def:liouville}
        Following \cite{McDuffSalamon1998}, a \emph{Liouville domain} is a compact manifold $W$ with boundary $Y=\partial W$, equipped with a one-form $\lambda$ such that $\omega=\D\lambda$ is symplectic and the Liouville vector field $Z$, defined by $\iota_Z\omega=\lambda$, points outward on $Y$. Its \emph{completion} $\widehat W$ is obtained by attaching the end $[1,\infty)\times Y$, on which $\lambda=r\alpha$ with $\alpha=\lambda|_Y$. We identify a two-sided collar of $Y$ with $(1-\delta,1+\delta)\times Y$.
    \end{definition}

    The \emph{symplectization} of $(Y,\alpha)$ is $((0,\infty)\times Y,\D(r\alpha))$ \cite{Geiges2008}. For any contactomorphism $\varphi$ preserving the coorientation of $\xi=\ker\alpha$, write $\varphi^*\alpha=\kappa_\varphi\alpha=e^{g_\varphi}\alpha$. When considering the fixed map $f$, we abbreviate $g=g_f$. We fix $\omega=\D\lambda$ and the symplectic Hamiltonian convention
    \begin{equation}\label{eq:hamconv}
        \iota_{X_H}\omega=-\D H.
    \end{equation}
    Denote by $R_s^\alpha$ the Reeb flow determined by $\alpha$. The identity \(\D(r-1)=\D r\) implies that \(r\) and \(r-1\) have the same Hamiltonian vector field. Our sign convention gives \(X_r=X_{r-1}=R_\alpha\) on the symplectization. We choose \(r-1\) as the defining Hamiltonian because \(\Sigma=\{r=1\}=(r-1)^{-1}(0)\).

    \begin{definition}[Contact conformal factor and translated point]\label{def:tp}
        Let $f\in\Cont(Y,\xi)$ satisfy $f^*\alpha=e^g\alpha$, and for $n\in\N$ let
        \begin{equation}\label{eq:cocycle}
            G_n(x)=\sum_{j=0}^{n-1}g(f^jx),\qquad (f^n)^*\alpha=e^{G_n}\alpha.
        \end{equation}
        A \emph{translated point pair} of $f^n$ is a pair $(x,\eta)\in Y\times\R$ for which
        \begin{equation}\label{eq:tp}
            f^n(x)=R_\eta^\alpha(x),\qquad G_n(x)=0.
        \end{equation}
        We denote this set by $\TP(f^n;\alpha)$ and use $\TP(\varphi;\alpha)$ analogously for any contactomorphism $\varphi$. When the contact form is fixed, we write $\TP(\varphi;I):=\{(x,\eta)\in\TP(\varphi;\alpha):\eta\in I\}$. A point $x\in Y$ is called a \emph{translated point} of $f^n$ with respect to $\alpha$ if some $\eta$ makes $(x,\eta)$ a translated point pair.
    \end{definition}

    \begin{definition}[Nondegenerate translated pairs]\label{def:regular}
        Let
        \[
        \begin{aligned}
            \mathcal T_n:Y\times\mathbb R &\longrightarrow Y\times Y\times\mathbb R,\\
            (x,\eta)&\longmapsto\bigl(f^n(x),R_\eta^\alpha(x),G_n(x)\bigr).
        \end{aligned}
        \]
        A translated point pair $(x,\eta)$ is \emph{nondegenerate} if $\mathcal T_n$ is transverse to the submanifold $\Delta_Y\times\{0\}$ at $(x,\eta)$, where \(\Delta_Y=\{(y,y):y\in Y\}\). Equivalently, the homogeneous linear system below has only the trivial solution $(\zeta,\sigma)\in T_xY\times\R$.
        \begin{equation}\label{eq:regular}
        \begin{cases}
            Df_x^n\zeta-DR_\eta^\alpha|_x\zeta-\sigma R_\alpha(f^n(x))=0,\\
            \D G_n(x)\zeta=0.
        \end{cases}
        \end{equation}
        A map is \emph{strongly translated-nondegenerate} if this condition holds for every pair of every positive iterate. This is the iteratewise form of the translated-pair nondegeneracy convention. Proposition~\ref{prop:normalized-hessian-translated-pair-kernel} proves the equivalence between translated-pair nondegeneracy and critical-point nondegeneracy.
    \end{definition}

    The defining Hamiltonian $\Fsig$ is a function on $\widehat{W}$, whereas the contact Hamiltonians $k_z$ are functions on $Y$ with time-one maps $\varphi_z$. We denote ordinary loops by $u$ and tomograph parameters by $z\in\Bpar\subset\R^N$. The Hamiltonians $K_{n,t}$ and $L_{z,t}$ generate the inverse lifts of $f^n$ and $\varphi_z$, respectively, with time supports specified in Subsection~\ref{section:analytic-admissibility-and-barcodes}. We write $X_A$ for a symplectic Hamiltonian vector field on $\widehat{W}$ and $X_k^\alpha$ for a contact Hamiltonian vector field on $Y$.

    \begin{lemma}\label{lem:radial-filling}
        If $\alpha=u\alpha_{\rm can}$ on $M=S^*S$, with $u>0$ smooth, then $(M,\alpha)$ is the contact boundary of a Liouville domain in $T^*S$.
    \end{lemma}
    \begin{proof}
        Let $\iota_u(q,p)=(q,u(q,p)p)$ for $|p|_h=1$. The tautological form satisfies $\iota_u^*\lambda_{\rm can}=u\alpha_{\rm can}$. The radial graph $\iota_u(S^*S)$ bounds the compact fiberwise star-shaped domain consisting of the zero section and all $rp$ with $0<r\le u(q,p)$, $|p|_h=1$. This domain has a smooth boundary along which the radial Liouville field is outward transverse. Equipped with the restricted tautological form, it therefore gives the required filling.
    \end{proof}
\section{Contact lemmas}\label{section:Contact-lemmas}
    We establish the contact perturbation and approximation lemmas needed to construct periodic points. These lemmas account for the derivative constraints arising from the contact condition, which obstruct the direct application of standard perturbation tools, such as Franks' lemma and the connecting lemma, formulated in $\Diff^1(M)$. We work in the standard contact model $(\rn{2d+1},\D{z}-x\cdot\D{y})$.
\subsection{A contact lift}
    In Section~\ref{subsection:flatten-center-germ} we apply the following lemma to the inverse return germ on a local invariant center curve. It realizes a one-dimensional germ $\psi$, with $\psi(0)=0$ and $\psi'(0)=1$, on the coordinate axis by a compactly supported, $C^1$-small contact map.

    Let $\alpha_0=\D{t}-x\cdot\D{y}$. For any orientation-preserving $C^2$ diffeomorphism $\psi:(\R,0)\to(\R,0)$, we define its standard \emph{contact lift} by
    \begin{equation}
        \mca{L}_{\psi}(x,y,t)=(\psi'(t)x,y,\psi(t)).
    \end{equation}
    A direct computation gives $\mca{L}_{\psi}^*\alpha_0=\psi'(t)\alpha_0$, so $\mca{L}_\psi$ is contact. 
    \begin{lemma}[Contact lifts]\label{lemma:contact lift}
        Consider the standard contact form $\alpha_0=\D{t}-x\cdot\D{y}$ on $(x,y,t)\in\rn{d}\times\rn{d}\times\R$. Let $\psi:(\R,0)\to(\R,0)$ be a local $C^2$ diffeomorphism such that $\psi'(0)=1$. Then for every $\varepsilon>0$ there are a $\delta>0$ and a lift $\mca{L}_{\psi,\delta}:\rn{2d+1}\to\rn{2d+1}$ such that
        \begin{itemize}
            \item $\mca{L}_{\psi,\delta}$ is contact and compactly supported in some small neighborhood $U_{cc}$ of the origin;
            \item $\tnorm{\mca{L}_{\psi,\delta}-\id}{C^1}\le\varepsilon$;
            \item $\mca{L}_{\psi,\delta}(0,0,t)=(0,0,\psi(t))$ for $\abs{t}\le\delta$.
        \end{itemize}
        More generally, for every integer $k\ge1$, if $\psi$ is $C^{k+1}$, the lift can be chosen $C^k$; it can be chosen $C^\infty$ when $\psi$ is $C^\infty$.
    \end{lemma}
    \begin{proof}
        Choose fixed smooth cutoffs $\zeta=1$ on $[-1,1]$, supported in $(-2,2)$, and $\chi=1$ on $B_1^{2d}(0)$, supported in $B_2^{2d}(0)$. Set $e_\delta(t)=\zeta(t/\delta)(\psi(t)-t)$ and $\varphi_s(t)=t+se_\delta(t)$, $0\le s\le1$. Taylor's theorem gives $\|e_\delta\|_{C^0}\le C\delta^2$, $\|e_\delta'\|_{C^0}\le C\delta$ and $\|e_\delta''\|_{C^0}\le C$. Thus $\varphi_s$ is increasing for small $\delta$, and its Eulerian velocity $v_s=e_\delta\circ\varphi_s^{-1}$ satisfies, uniformly in $s$,
        \begin{equation}\label{eq:velocity-estimates}
            \|v_s\|_{C^0}\le C\delta^2,\qquad
            \|v_s'\|_{C^0}\le C\delta,\qquad
            \|v_s''\|_{C^0}\le C.
        \end{equation}
        Put $\rho=\sqrt{\delta}$, $\chi_\rho(x,y)=\chi(x/\rho,y/\rho)$, and use the compactly supported contact Hamiltonian $K_s(x,y,t)=\chi_\rho(x,y)v_s(t)$. Its vector field is
        \begin{equation}\label{eq:X_{K_{delta,s}}}
            X_{K_s}=
                (v_s\partial_y\chi_\rho+x\chi_\rho v_s')\cdot\partial_x
                -v_s\partial_x\chi_\rho\cdot\partial_y
                +(\chi_\rho-x\cdot\partial_x\chi_\rho)v_s\partial_t.
        \end{equation}
        Since $|x|\le2\rho$ on the support and $\|D^j\chi_\rho\|_{C^0}\le C\rho^{-j}$ for $j=1,2$, \eqref{eq:velocity-estimates} and \eqref{eq:X_{K_{delta,s}}} give
        \begin{equation}\label{eq:localized-vector-field-small}
            \sup_s\|X_{K_s}\|_{C^0}\le C\delta^{3/2},\qquad
            \sup_s\|X_{K_s}\|_{C^1}\le C\sqrt{\delta}.
        \end{equation}
        The time-one map $h_{\delta,1}$ is therefore $C^1$-close to the identity by the flow equation and Gronwall's inequality. In particular, $\|h_{\delta,1}-\id\|_{C^1}\le e^{C\sqrt{\delta}}-1$. Near the center axis the vector field is $xv_s'(t)\cdot\partial_x+v_s(t)\partial_t$, whose flow is the standard lift of $\varphi_s$. Hence $h_{\delta,1}(0,0,t)=(0,0,\psi(t))$ for $|t|\le\delta$. The support of $h_{\delta,1}$ is contained in $B_{2\sqrt{\delta}}^{2d}(0)\times(-2\delta,2\delta)$, which can be contained in any prescribed neighborhood of the origin. Take $\mca{L}_{\psi,\delta}=h_{\delta,1}$ and shrink $\delta$ to achieve the desired bound.

        If $\psi$ is $C^{k+1}$, then $v_s$, $K_s$ are $C^{k+1}$ in space. The vector field $X_{K_s}$ and the flow $h_{\delta,s}$ are therefore $C^k$. This also applies to the smooth case, and the resulting map is also smooth.
    \end{proof}
\subsection{A contact connecting lemma}
    The following point-moving lemma realizes a small displacement by a contact perturbation supported in a fixed Darboux neighborhood. We use it in both Section~\ref{subsection:contact-blender-connection} and Section~\ref{subsection:Local contact perturbations and the center multiplier}.

    \begin{lemma}[Point-moving]\label{lemma:point-moving}
        Let $V$ be a Darboux neighborhood of a point $x_*$, and let $y_j\to x_*$. There are smooth contactomorphisms $h_j$, for all sufficiently large $j$, such that
        \[
            \supp(h_j)\subset\subset V,\qquad h_j(x_*)=y_j,\qquad
            h_j\longrightarrow\id\quad\text{in }C^1.
        \]
    \end{lemma}
    \begin{proof}
        Center Darboux coordinates at $x_*$, with $\alpha_0=\D{z}-x\cdot\D{y}$, and write $y_j=(a_j,b_j,c_j)$. Choose a fixed smooth cutoff $\chi$ supported in the chart and equal to one on a neighborhood of the origin. Set
        \[
            H_j(x,y,z)=\chi(x,y,z)
            \Bigl(-b_j\cdot x+a_j\cdot y+c_j-\frac12a_j\cdot b_j\Bigr).
        \]
        Since $\chi=1$ near $x_*$, the contact Hamiltonian trajectory starting at the origin is
        \[
            x(s)=s a_j,\qquad y(s)=s b_j,\qquad
            z(s)=s\Bigl(c_j-\frac12a_j\cdot b_j\Bigr)+\frac{s^2}{2}a_j\cdot b_j .
        \]
        For large $j$, this trajectory remains in the region $\chi=1$ for $0\le s\le1$, and its endpoint is $y_j$. Since the cutoff $\chi$ is fixed,
        \[
            \|H_j\|_{C^2}\le
            C_{\chi}\bigl(|a_j|+|b_j|+|c_j|+|a_j|\,|b_j|\bigr)\longrightarrow0.
        \]
        The contact Hamiltonian vector fields $X_{H_j}$ therefore tend to zero in $C^1$, and their time-one maps tend to the identity in $C^1$. Conjugating these maps by the Darboux chart and extending them by the identity proves the assertion.
    \end{proof}
\subsection{Contact approximation}\label{subsection:contact-approximation}
    We prove relative smoothing in the ordinary $C^1$ topology. A vector field or curve is \emph{horizontal} if it is tangent to $\xi=\ker\alpha$; see \cite[Section~2.1]{Geiges2008}. For a continuous function $\lambda$, $X\lambda$ denotes differentiation along the flow of a smooth horizontal field $X$. The fixed smooth contact form gives $TM=\xi\oplus\R R_\alpha$. In Darboux coordinates $(q,p,z)$, the Reeb line is $\R\partial_z$, whereas the \emph{vertical Lagrangian subspace} $\Span\{\partial_{p_i}\}$ lies in $\xi$.

    \begin{lemma}[H\"older continuity of the conformal factor]\label{lemma:conformal-factor-holder}
        Let $(M^{2d+1},\alpha)$ be a smooth contact manifold, with $d\ge1$, and let $f$ be a coorientation-preserving $C^1$ contactomorphism satisfying $f^*\alpha=\lambda\alpha$, where $\lambda>0$ is initially assumed only continuous. Then $\lambda\in C^{0,1/2}_{\loc}$.
    \end{lemma}
    \begin{proof}
        The identity $\D(f^*\alpha)=f^*(\D\alpha)$ holds distributionally. Consequently,
        \begin{equation}\label{eq:approximation-horizontal-factor}
            \D\lambda\wedge\alpha=f^*(\D\alpha)-\lambda\,\D\alpha.
        \end{equation}
        If $X$ is a smooth vector field tangent to $\xi=\ker\alpha$ and $R_\alpha$ is the Reeb field, evaluation on $(X,R_\alpha)$ gives
        \[
            X\lambda=f^*(\D\alpha)(X,R_\alpha).
        \]
        The right-hand side is continuous. At a point where $X\neq0$, a flow box for $X$ turns this distributional equation into the fundamental theorem of calculus along the flow, so the horizontal derivative is classical along the corresponding flow line. At a zero of $X$, its flow is stationary and both sides vanish.

        In a Darboux chart with $\alpha=\D{z}-p\cdot\D{q}$, use the horizontal fields $X_i=\partial_{p_i}$ and $Y_i=\partial_{q_i}+p_i\partial_z$, whose flows are
        \[
            \Phi_{X_i}^{s}(q,p,z)=(q,p+s e_i,z),\qquad
            \Phi_{Y_i}^{s}(q,p,z)=(q+s e_i,p,z+s p_i),
        \]
        where $e_i$ is the $i$th coordinate vector of $\R^d$. The derivatives $X_i\lambda$ and $Y_i\lambda$ are continuous and bounded on compact subcharts. For nearby points $x=(q,p,z)$ and $x'=(q',p',z')$, put $u=p'-p$ and $v=q'-q$.

        First change $p$ to $p'$ and then $q$ to $q'$ along the horizontal paths
        \[
            \gamma_p(t)=(q,p+t u,z),\qquad
            \gamma_q(t)=(q+t v,p',z+t p'\cdot v),\qquad 0\le t\le1.
        \]
        These are flow segments of $\sum_i u_iX_i$ and $\sum_i v_iY_i$, respectively, so the classical horizontal derivative identity applies along each. The second path ends at $(q',p',z_*)$, where $z_*=z+p'\cdot(q'-q)$, leaving the Reeb displacement
        \[
            \Delta:=z'-z_*=z'-z-p'\cdot(q'-q),\qquad
            |\Delta|\le C|x-x'|.
        \]

        Set $a=\sqrt{|\Delta|}$ and $b=\sgn(\Delta)\sqrt{|\Delta|}$, so that $ab=\Delta$. From $(q',p',z_*)$, follow successively the unit-time flows of $aX_1$, $bY_1$, $-aX_1$, and $-bY_1$. Their projection to the $(q_1,p_1)$-plane is a rectangle, while their net change in $z$ is $(p'_1+a)b-p'_1b=ab=\Delta$. Thus these four horizontal segments end at $x'$, forming the commutator rectangle associated with $[X_1,Y_1]=\partial_z$. When $\Delta=0$, take $a=b=0$ and omit the four constant segments.

        After shrinking the subchart, all six segments for sufficiently close $x,x'$ lie in a fixed larger compact subchart. Let $L$ bound $|X_i\lambda|$ and $|Y_i\lambda|$ there for $1\le i\le d$. Integrating the horizontal derivatives gives
        \begin{equation}\label{eq:approximation-factor-holder}
            \begin{aligned}
                |\lambda(x')-\lambda(x)|
                &\le L\biggl(\,\sum_{i=1}^d(|u_i|+|v_i|)+2|a|+2|b|\biggr)\\
                &\le C\bigl(|x-x'|+\sqrt{|\Delta|}\bigr)
                \le C|x-x'|^{1/2}.
            \end{aligned}
        \end{equation}
        Thus $\lambda\in C^{0,1/2}_{\loc}$.
    \end{proof}

    The next two lemmas collect the spatial and radius estimates needed. Ordinary mollification is standard; see \cite[Addendum on mollification, pp.~25--30]{Schmidt2026}. The weighted quotient below requires the additional estimate that we prove explicitly. 

    \medskip\noindent\textbf{Convolution notations}

    Let $\Omega\subset\R^m$ be open. Fix a nonnegative symmetric function $\rho\in C_c^\infty(\R^m)$ supported in the closed unit ball and satisfying $\int\rho=1$. Its rescaling is $\rho_t(h):=t^{-m}\rho(h/t)$, for $t>0$. Thus $\int\rho_t=1$ and $\supp\rho_t\subset\cball{t}{0}$. Fix compact sets $K_0\subset\interior K_1\subset K_1\subset\Omega$ and a radius $0<t_0\le1$ satisfying $\cball{t_0}{v}\subset K_1$ for every $v\in K_0$.

    For a continuous scalar- or vector-valued function $u$ on $\Omega$, convolution at these points is defined by
    \[
        (u*\rho_t)(v):=\int_\Omega u(y)\rho_t(v-y)\,\D y,
        \qquad v\in K_0,\quad 0<t\le t_0.
    \]
    The ball condition ensures that all function values and line segments used below stay in $K_1$. For a continuous function $g$ on $K_1$, define its \emph{modulus of continuity} by
    \[
        \omega_g(r):=\sup_{v,w\in K_1,\,|v-w|\le r}\|g(v)-g(w)\|.
    \]
    Uniform continuity of $g$ gives $\omega_g(r)\to0$ as $r\to0$.
    
    All estimates in the next two lemmas are uniform for $v\in K_0$ and $0<t\le t_0$. Their constants may depend on the fixed domains, kernel, and functions, but not on $v$ or $t$. Throughout the rest of this subsection, $D$ denotes the spatial derivative, and $\partial_v,\partial_t$ denote partial derivatives in the indicated variables. 

    \begin{lemma}\label{lemma:ordinary-convolution}
        In the preceding setting, let $Q\in C^1(\Omega;\R^k)$, with $k\ge1$, let $a\in C^{0,1/2}_{\loc}(\Omega)$, and let $W\in C^{1,1/2}_{\loc}(\Omega)$. Set $Q_t:=Q*\rho_t$, $a_t:=a*\rho_t$, $S(v,t):=(W*\rho_t)(v)$. Then the following estimates hold.
        \begin{enumerate}[(i)]
            \item For the $C^1$ function $Q$,
            \begin{equation}\label{eq:Q_t-est.}
                \|Q_t(v)-Q(v)\|\le Ct,\qquad
                \|DQ_t(v)-DQ(v)\|\le\omega_{DQ}(t),\qquad
                \|\partial_tQ_t(v)\|\le C\omega_{DQ}(t).
            \end{equation}
            In particular, $Q_t\to Q$ in $C^1$ and $\pps{t}Q_t\to0$ as $t\to0$.
            \item For the H\"older function $a$,
            \begin{equation}\label{eq:a_t-est.}
                |a_t(v)-a(v)|\le C\sqrt{t},\qquad \|Da_t\|+|\partial_ta_t|\le Ct^{-1/2}.
            \end{equation}
            If $a_*:=\min_{K_1}a>0$, then $a_t(v)\ge a_*$.
            \item For the function $W$,
            \begin{equation}\label{eq:S(v,t)-est.}
                \begin{aligned}
                    |S(v,t)-W(v)|&\le Ct^{3/2},\\
                    \|\partial_vS(v,t)-DW(v)\|+\abs{\pps{t}S(v,t)}&\le C\sqrt{t},\\
                    \|\partial_v\partial_tS(v,t)\|+|\partial_t^2S(v,t)|&\le Ct^{-1/2}.
                \end{aligned}
            \end{equation}
        \end{enumerate}
        Defining $S(v,0):=W(v)$ gives a $C^1$ extension to $t=0$, with $\partial_vS(v,0)=DW(v)$ and $\partial_tS(v,0)=0$.
    \end{lemma}
    \begin{proof}
        Write $h=t\zeta$. The rescaled kernels and their derivatives have support in $|h|\le t$ and satisfy $\|D\rho_t\|_{L^1}+\|\partial_t\rho_t\|_{L^1}\le C/t$ and $\|\partial_t^2\rho_t\|_{L^1}\le C/t^2$. Each of the derivative kernels $\pps{i}\rho_t$ and $\pps{t}\rho_t$ has integral zero, and symmetry of $\rho$ also gives $\int h\rho_t(h)\,\D h=0$. Differentiating this identity in $t$ provides the corresponding first-moment cancellations.

        For $Q$, a standard mean-value estimate gives $\|Q_t-Q\|\le Ct$, while $DQ_t=(DQ)*\rho_t$ gives the derivative bound. Moreover,
        \[
            \partial_tQ_t(v)=-\int
            [DQ(v-t\zeta)-DQ(v)]\zeta\rho(\zeta)\,\D\zeta,
        \]
        proving the radius-derivative estimate. For $a$, averaging gives $|a_t-a|\le C\sqrt{t}$ and preserves
        the positive lower bound. Subtracting $a(v)$ in each derivative integral and using the kernel bounds gives
        \[
            \|Da_t(v)\|+\|\partial_ta_t(v)\|
                \le C\sqrt{t}\bigl(\|D\rho_t\|_{L^1}
                    +\|\partial_t\rho_t\|_{L^1}\bigr)
                \le Ct^{-1/2}.
        \]

        For $W$, its first-order Taylor remainder $R_v(y):=W(y)-W(v)-DW(v)\cdot(y-v)$ at $v$ is bounded by
        $C|y-v|^{3/2}$. The first-moment cancellation $S(v,t)-W(v)=\int R_v(y)\rho_t(v-y)\,\D{y}$ therefore yields $|S(v,t)-W(v)|\le Ct^{3/2}$. Analogously, convolution of $DW$ gives $|\partial_vS-DW|\le C\sqrt{t}$, since $DW$ is $1/2$-H\"older. The displayed formula for $\partial_tQ_t$ applied to $S(v,t)$ gives $|\partial_tS|\le C\sqrt{t}$. Use the identities
        \[
            \partial_v\partial_tS(v,t)
            =\int[DW(y)-DW(v)]\,\partial_t\rho_t(v-y)\,\D{y},\qquad
            \partial_t^2S(v,t)
            =\int R_v(y)\,\partial_t^2\rho_t(v-y)\,\D{y}.
        \]
        to obtain the bounds $\abs{\pps{v}\pps{t}S}\le C\sqrt{t}\Lnorm{\pps{t}\rho_t}{1}\le Ct^{-1/2}$ and $|\partial_t^2S|\le Ct^{3/2}\Lnorm{\pps{t}^2\rho_t}{1}\le Ct^{-1/2}$, respectively. These estimates prove all assertions, including the $C^1$ extension at $t=0$ and its stated derivatives.
    \end{proof}

    \begin{lemma}[Weighted convolution]\label{lemma:weighted-convolution}
        In the same setting, let $a\in C^{0,1/2}_{\loc}(\Omega)$ be positive and let $b\in C^1(\Omega,\R^k)$, with $k\ge1$. With $a_t=a*\rho_t$ and $a_*=\min_{K_1}a>0$, define the weighted convolution by
        \[
            b_t(v):=\frac{((ab)*\rho_t)(v)}{a_t(v)}.
        \]
        There is a constant $C$ such that
        \begin{equation}\label{eq:approximation-weighted-convolution}
            |b_t(v)-b(v)|\le Ct,\qquad \|Db_t(v)-Db(v)\|\le C\bigl(\omega_{Db}(t)+\sqrt{t}\,\bigr),\qquad|\partial_t b_t(v)|\le C.
        \end{equation}
        In particular, $b_t\to b$ in $C^1$ on $K_0$ as $t\to0$.
    \end{lemma}
    \begin{proof}
        We use the kernel derivative formulas, $L^1$ bounds, and cancellations proved in Lemma~\ref{lemma:ordinary-convolution}, which also gives $a_t(v)\ge a_*>0$. Put $L=\sup_{K_1}\|Db\|$. When $|y-v|\le t$, we have $|b(y)-b(v)|\le L|y-v|\le Lt$. Since $a>0$ and $\rho_t\ge0$, the weighted-average 
        $$
        b_t(v)-b(v)=\frac{1}{a_t(v)}\int a(y)\rho_t(v-y)[b(y)-b(v)]\,\D{y}
        $$
        gives $|b_t(v)-b(v)|\le Lt$. In particular, $|b(y)-b_t(v)|\le\abs{b(y)-b(v)}+\abs{b(v)-b_t(v)}\le2Lt$ when $|y-v|\le t$. 
        
        Define $r_v(y)=b(y)-b(v)-Db(v)(y-v)$. The argument used to estimate the remainder $R_v(y)$ of $W$ in Lemma~\ref{lemma:ordinary-convolution} yields $|r_v(y)|\le t\,\omega_{Db}(t)$ when $\abs{y-v}\le t$.

        Differentiating $b_t$ gives
        \[
            \partial_i b_t(v)=\frac{1}{a_t(v)}\int
            \partial_i\rho_t(v-y)a(y)[b(y)-b_t(v)]\,\D y.
        \]
        Move $a_t(v)$ to the left and split $a(y)=a(v)+[a(y)-a(v)]$ in the integrand. Use the identities for the constant and first moments of the kernel $\partial_i\rho_t$ to obtain
        \[
            a_t(v)\partial_i b_t(v)
            =a(v)\pps{i}b(v)+a(v)\int\partial_i\rho_t(v-y)r_v(y)\,\D y+\int\partial_i\rho_t(v-y)[a(y)-a(v)][b(y)-b_t(v)]\,\D y.
        \]
        Since $a$ is bounded, we have $\abs{a_t\pps{i}b_t-a\pps{i}b}\le C(\omega_{Db}(t)+\sqrt{t})$. Dividing it by $a_t\ge a_*>0$ and using $|a_t-a|\le C\sqrt{t}$, 
        \[
            |\partial_i b_t(v)-\partial_i b(v)|
            \le\frac{1}{a_t}\abs{a_t\pps{i}b_t-a\pps{i}b}+\frac{\abs{a_t-a}}{a_t}\abs{\pps{i}b}
            \le C\bigl(\omega_{Db}(t)+\sqrt{t}\,\bigr).
        \]
        which proves the spatial derivative estimate. The same quotient differentiation in the radius gives
        \[
            \partial_t b_t(v)=\frac{1}{a_t(v)}\int \partial_t\rho_t(v-y)a(y)[b(y)-b_t(v)]\,\D y.
        \]
        Its absolute value is bounded by $C$, since $\|\partial_t\rho_t\|_{L^1}\le C/t$ and $\abs{b(y)-b_t(v)}\le 2Lt$ on the kernel support. This proves the lemma.
    \end{proof}

    \begin{proposition}[Relative contact approximation]\label{prop:relative-contact-approximation}
        Let $(M^{2d+1},\alpha)$ be a smooth closed contact manifold, with $d\ge1$, and let $f$ be a coorientation-preserving $C^1$ contactomorphism. For every compact set $K\subset O$, with $O\subset M$ open, and every $\varepsilon_0>0$, there is a coorientation-preserving $C^1$ contactomorphism $\tilde{f}$ such that
        \[
            d_{C^1}(\tilde{f},f)<\varepsilon_0,
            \qquad \tilde{f}=f\quad\text{on }M\setminus O,
        \]
        and $\tilde{f}$ is $C^\infty$ on a neighborhood of $K$. Moreover, $\tilde{f}$ is smooth at every point $q$ whenever $f$ is smooth on a neighborhood of $q$. In particular, taking $K=O=M$, smooth contactomorphisms are $C^1$-dense in the space of coorientation-preserving $C^1$ contactomorphisms.
    \end{proposition}
    \begin{proof}        
        We smooth a local generating function at a radius vanishing on the boundary, then patch finitely many such modifications. The smooth generating-function construction is standard \cite[Sections~3.1--3.2]{Sandon2011}; the estimates below justify the reconstruction at $C^1$ regularity.

        \emph{Step 1: a local generating function.}
        Write $f^*\alpha=\lambda\alpha$, with $\lambda>0$ continuous. Fix a point and choose source and target Darboux charts, both with contact form $\D{z}-p\cdot \D{q}$. Write
        \[
            f(q,p,z)=(Q(q,p,z),P(q,p,z),Z(q,p,z)).
        \]
        The charts may be chosen so that $\partial_pP$ is invertible at the base point. By Lemma~\ref{lemma:conformal-factor-holder}, $Df|_\xi$ is conformally symplectic, so the image under $Df$ of the source vertical Lagrangian subspace is Lagrangian in the target contact hyperplane. A symplectic linear change, together with its smooth contact lift, then makes this image the target vertical subspace. On a smaller neighborhood the mixed-coordinate projection
        \[
            (q,p,z)\longmapsto(q,P(q,p,z),z)
        \]
        is therefore a $C^1$ diffeomorphism. Use $(q,P,z)$ as independent variables on its image $\Omega$, and write $H(q,P,z)$ for the original source $p$-coordinate. The inverse mixed-coordinate map is $A(q,P,z):=(q,H(q,P,z),z)$. From now on, $Q(q,P,z)$ and $Z(q,P,z)$ denote the original target coordinate functions composed with $A$. Set $a(q,P,z)=\lambda\circ A(q,P,z)$ and $B(q,P,z)=f\circ A(q,P,z)$. With $U=A(\Omega)$ and $V=f(U)$, we have the following commutative diagram
        \[
            \xymatrix@C=4em@R=3em{
                \Omega: (q,P,z) \ar[r]^{A} \ar[dr]_{B}
                    & U: (q,p,z) \ar[d]^{f} \\
                    & V: (Q,P,Z)
            }
            \qquad B=f\circ A.
        \]
        Define the scalar \emph{local generating function} on $\Omega$ by
        \[
            W(q,P,z):=Z(q,P,z)-P\cdot Q(q,P,z).
        \]
        Pulling the contact equation back by $A$ gives $\D{Z}-P\cdot\D{Q}=a(\D{z}-H\cdot \D{q})$, where all functions and differentials are now expressed in the mixed variables $(q,P,z)$. Thus
        \begin{equation}\label{eq:approximation-generating-derivatives}
            \D{W}=\D{Z}-\D(P\cdot Q)=a\,\D{z}-aH\cdot \D{q}-Q\cdot \D{P}.
        \end{equation}
        Equivalently, $\partial_zW=a>0$, $\partial_qW=-aH$ and $\partial_PW=-Q$.

        Conversely, the identities
        \[
            H=-\frac{\partial_qW}{\partial_zW},\qquad
            Q=-\partial_PW,\qquad Z=W+P\cdot Q
        \]
        determine $A=(q,H,z)$ and $B=(Q,P,Z)$ from $W$, and hence determine $f=B\circ A^{-1}$ on the source neighborhood $U$.

        \emph{Step 2: H\"older regularity and fixed-radius convolution.}
        From this step to the end of the proof, abbreviate the mixed variables by $v=(q,P,z)$, and set $m=2d+1$. Lemma~\ref{lemma:conformal-factor-holder} gives $\lambda\in C^{0,1/2}_{\loc}$. Since the inverse mixed-coordinate map $A$ is $C^1$, $a=\lambda\circ A\in C^{0,1/2}_{\loc}$ is positive. The functions $H,Q$ are $C^1$, so \eqref{eq:approximation-generating-derivatives} gives $W\in C^{1,1/2}_{\loc}$. Use the mollifier and convolution conventions fixed before Lemma~\ref{lemma:ordinary-convolution}. For $t>0$ small enough that the convolution stays in $\Omega$, set
        \begin{equation}\label{eq:def. of S,a,H,Q}
            S(v,t)=(W*\rho_t)(v),\qquad a_t=a*\rho_t,\qquad H_t=\frac{(aH)*\rho_t}{a*\rho_t}, \qquad Q_t=Q*\rho_t.
        \end{equation}
        Applying Lemma~\ref{lemma:ordinary-convolution} to $Q\in C^1$, and Lemma~\ref{lemma:weighted-convolution} with weight $a\in C\ssp{\loc}{0,1/2}$, $b=H\in C^1$ gives $H_t\to H$ and $Q_t\to Q$ in $C^1$ as $t\to0$, respectively.

        \emph{Step 3: a smoothing radius that vanishes at the boundary.}
        Fix a smooth cutoff $0\le\beta\le1$ compactly supported in $\Omega$, and set
        \[
            \delta(v)=\eta\beta(v)^4,\qquad
            \widehat{W}(v)=
            \begin{cases}
                S(v,\delta(v)),&\beta(v)>0,\\
                W(v),&\beta(v)=0.
            \end{cases}
        \]
        The charts and $\beta$ remain fixed throughout this step, while the positive real parameter $\eta$ will tend to zero. Choose fixed compact subcharts $K_0\subset\interior K_1\subset K_1\subset\Omega$, with $\supp\beta\subset\interior K_0$, and take $0<\eta\le1$ sufficiently small that all convolution balls centered in $K_0$ stay in $K_1$. All norms and constants below refer to these fixed subcharts, with constants independent of $\eta$.
        
        Put $a_*=\min_{K_1}a>0$. Applying Lemma~\ref{lemma:ordinary-convolution} to $a$ and $W$ gives \eqref{eq:a_t-est.} and \eqref{eq:S(v,t)-est.}, as well as the $C^1$ extension $S(v,0)=W(v)$ with $\partial_tS(v,0)=0$. Composing with the smooth nonnegative function $\delta$ therefore gives $\widehat{W}\in C^1$.

        \emph{Reconstruction of $\hat{f}$ from $\widehat{W}$.}
        Define 
        \begin{equation}\label{eq:def-bar-a,H,Q}
            \bar{a}(v)=a_{\delta(v)}(v),\qquad
            \widebar{H}(v)=H_{\delta(v)}(v),\qquad
            \widebar{Q}(v)=Q_{\delta(v)}(v),
        \end{equation}
        At $\beta=0$ use the extensions $a_0=a$, $H_0=H$, $Q_0=Q$. Define the covector $E(v)=\partial_tS(v,\delta(v))D\delta(v)$, extended by zero at $\beta=0$. Then $D\widehat{W}(v)=\partial_vS(v,\delta(v))+E(v)$. Write $E(v)=(E^q(v),E^P(v),E^z(v))$, where $E^q,E^P,E^z$ are the components of the covector $E$ in the $\D{q},\D{P},\D{z}$ directions. At radius $t=\delta(v)$, \eqref{eq:approximation-generating-derivatives} and definitions \eqref{eq:def. of S,a,H,Q}, \eqref{eq:def-bar-a,H,Q}, together with $\widehat{W}(v)=S(v,\delta(v))$, give
        \begin{equation}\label{eq:partial-derivatives-hatW}
            \partial_q\widehat{W}=-\bar{a}\widebar{H}+E^q,\qquad
            \partial_P\widehat{W}=-\widebar{Q}+E^P,\qquad
            \partial_z\widehat{W}=\bar{a}+E^z.
        \end{equation}
        We can therefore reconstruct the coordinates by
        \begin{equation}\label{eq:approximation-reconstruction}
            \widehat{Q}:=-\partial_P\widehat{W}\qquad
            \widehat{H}:=-\frac{\partial_q\widehat{W}}{\partial_z\widehat{W}}\qquad
            \widehat{Z}:=\widehat{W}+P\cdot\widehat{Q},
        \end{equation} 
        provided that $\pps{z}\widehat{W}$ is positive, which will be verified later.

        \emph{$C^1$-convergence.}
        We prove that the modified coordinates $\widehat{Q}$, $\widehat{H}$ and $\widehat{Z}$ converge to $Q$, $H$ and $Z$, respectively, as $\eta\to0$. We first control the chain-rule correction $E$, and then verify the variable-radius convergence of $\widebar{Q}$ and $\widebar{H}$. Writing $E^i$ for the $i$th component of $E$ in the coordinates $v$, differentiation on $\{\beta>0\}$ gives
        \[
            \partial_jE^i(v)
                =\bigl[(\partial_j\partial_tS)(v,\delta(v))+\partial_t^2S(v,\delta(v))\partial_j\delta(v)\bigr]\partial_i\delta(v)\\
                +\partial_tS(v,\delta(v))\partial_j\partial_i\delta(v).
        \]
        Thus \eqref{eq:S(v,t)-est.} yields
        \begin{equation}\label{eq:approximation-cutoff-errors}
            \norm{E}{}\le C\sqrt{\delta}\,\norm{D\delta}{},\qquad
            \norm{DE}{}\le C\bigl(\delta^{-1/2}+\delta^{-1/2}\norm{D\delta}{}\bigr)\norm{D\delta}{}+C\sqrt{\delta}\tnorm{D^2\delta}{}.
        \end{equation}
        Since $\beta$ is fixed, $D\delta=4\eta\beta^3D\beta$ and $D^2\delta=12\eta\beta^2D\beta\otimes D\beta+4\eta\beta^3D^2\beta$. Allowing the constant $C$ to depend on $\beta$ yields $\norm{D\delta}{}\le C(\beta)\eta\beta^3$ and $\tnorm{D^2\delta}{}\le C(\beta)\eta\beta^2$. Substitution into \eqref{eq:approximation-cutoff-errors} gives the explicit bounds
        \begin{equation}\label{eq:C^1-est. for E}
            \norm{E}{}\le C\eta^{3/2}\beta^5,\qquad
            \norm{DE}{}\le C\sqrt{\eta}\,\beta.
        \end{equation}
        A $C^1$ correction on $\{\beta>0\}$ extends by zero to a $C^1$ function if its value is $O(\beta^{\nu})$ for some $\nu>1$ and its derivative tends to zero as $\beta\to0$. Indeed, at each zero $v_0$, the Lipschitz bound $\beta(v)\le C|v-v_0|$ makes the value $o(|v-v_0|)$, so the extension is differentiable with derivative zero, and the derivative bound gives continuity. Applying this criterion to $E$ yields
        \[
            \|E\|_{C^1(K_0)}\le C(\beta)(\eta^{3/2}\beta^5+\sqrt{\eta}\beta)\le C(\beta)\sqrt{\eta},
        \]
        which tends to zero as $\eta\to0$. In particular, for all sufficiently small $\eta$, $\bar{a}+E^z\ge a_*-\norm{E}{}>0$.

        On $\{\beta>0\}$, the chain rule gives $D\widebar{H}=(D_vH_t+\partial_tH_tD\delta)|_{t=\delta}$ and $D\widebar{Q}=(D_vQ_t+\partial_tQ_tD\delta)|_{t=\delta}$. Lemmas~\ref{lemma:ordinary-convolution} and \ref{lemma:weighted-convolution} therefore imply $\|\widebar{H}-H\|+\|\widebar{Q}-Q\|\le C\delta$ and
        \[
            \|D\widebar{H}-DH\|
                \le C\bigl(\omega_{DH}(\delta)+\sqrt{\delta}+\|D\delta\|\bigr),\quad
            \|D\widebar{Q}-DQ\|
                \le (1+C|D\delta|)\omega_{DQ}(\delta).
        \]
        The value errors are $O(\eta\beta^4)$ and the derivative errors tend to zero as $\beta\to0$, so the same criterion applies to $\widebar{H}-H$ and $\widebar{Q}-Q$. Thus both barred maps are $C^1$, with derivatives $DH,DQ$ on $\{\beta=0\}$, and
        \[
            \|\widebar{H}-H\|_{C^1}+\|\widebar{Q}-Q\|_{C^1}
            \le C\bigl(\omega_{DH}(\eta)+\omega_{DQ}(\eta)+\sqrt{\eta}\bigr)
            \longrightarrow0.
        \]
        Also $\|\widehat{W}-W\|_{C^1}\to0$ by the convolution estimates and the bound on $E$. It remains to control the quotient correction $\widehat{H}$. Separating $\widebar{H}$ in the expression for $\widehat{H}$ gives $\widehat{H}=\widebar{H}-\frac{E^q+\bar{H} E^z}{\bar{a}+E^z}$. On $\{\beta>0\}$, $D\bar{a}(v)=\bigl((\partial_va_t)(v)+(\partial_ta_t)(v)D\delta(v)\bigr)\big|_{t=\delta(v)}$ and hence by \eqref{eq:a_t-est.}, $\norm{D\bar{a}}{}\le C\delta^{-1/2}(1+\norm{D\delta}{})$. Combining this estimate with \eqref{eq:C^1-est. for E} gives
        \begin{equation}\label{eq:approximation-quotient-error}
            \begin{aligned}
                \norm{E}{}\|D\bar{a}\|
                &\le C\sqrt{\delta}\norm{D\delta}{}\cdot C\delta^{-1/2}(1+\norm{D\delta}{})\\
                &\le C\norm{D\delta}{}(1+\norm{D\delta}{})\le C\eta\beta^3(1+\eta\beta^3).
            \end{aligned}
        \end{equation}
        in which the negative power of $\delta$ is canceled. To differentiate the quotient, set on $\{\beta>0\}$
        \[
            N=E^q+\widebar{H} E^z,\qquad
            G=\bar{a}+E^z,\qquad F=N/G.
        \]
        We have $G\ge a_*/2$, $\|N\|\le C\norm{E}{}$, thus $\abs{F}\le(2C/a_*)\norm{E}{}\le C\eta^{3/2}\beta^5$. Furthermore, using the uniform $C^1$ bound for $\widebar{H}$ gives $\norm{DN}{}\le C(\norm{DE}{}+\norm{E}{})$. The quotient rule reads $DF=DN/G-N(D\bar{a}+DE^z)/{G^2}$. Therefore by \eqref{eq:C^1-est. for E} and \eqref{eq:approximation-quotient-error}, on $\set{\beta>0}$
        \[
            \norm{DF}{}\le\frac{2C}{a_*}(\norm{DE}{}+\norm{E}{})+\frac{4}{a_*^2}\cdot C\norm{E}{}(\norm{D\bar{a}}{}+\norm{DE}{})\le C\sqrt{\eta}\beta.
        \]
        The value and derivative bounds for $F$ satisfy the same criterion, so its zero extension is $C^1$, with derivative zero on $\{\beta=0\}$ and $\|F\|_{C^1(K_0)}\le C\sqrt{\eta}$. Since $\widehat{H}=\widebar{H}-F$ and $\widehat{Z}-Z=(\widehat{W}-W)+P\cdot(\widehat{Q}-Q)$, these estimates give
        \[
            \begin{aligned}
                \|(\widehat{H},\widehat{Q},\widehat{Z})-(H,Q,Z)\|_{C^1(K_0)}
                    &\le C\bigl(\omega_{DH}(\eta)+\omega_{DQ}(\eta)+\sqrt{\eta}\,\bigr)\longrightarrow0.
            \end{aligned}
        \]
        Thus, with the charts and $\beta$ fixed, decreasing $\eta$ makes the reconstructed coordinates arbitrarily $C^1$-close to the originals. All reconstructed quantities are smooth where $\beta>0$.

        \emph{Step 4: the local relative perturbation.}
        Put $\widehat{A}(v)=(q,\widehat{H}(v),z)$ and $\widehat{B}(v)=(\widehat{Q}(v),P,\widehat{Z}(v))$. The map $\widehat{A}\circ A^{-1}$ is $C^1$-close to the identity and equals the identity outside a compact subset of the source chart. Extend it by the identity in Euclidean coordinates. For sufficiently small $\eta$, this extension is a diffeomorphism, so $\widehat{A}$ maps $\Omega$ diffeomorphically onto $A(\Omega)$. Define $\hat{f}=\widehat{B}\circ\widehat{A}^{-1}$ on $A(\Omega)$, and use $f$ outside this chart. These definitions agree near its boundary. With $U=A(\Omega)=\widehat{A}(\Omega)$, the original and reconstructed maps form the two commuting triangles below, in which the target maps take values in $M$.
        \[
            \xymatrix@C=4em@R=3em{
                U \ar[d]_{f}
                    & \Omega \ar[l]_{A} \ar[r]^{\widehat{A}}
                        \ar[dl]^{B} \ar[dr]_{\widehat{B}}
                    & U \ar[d]^{\hat{f}} \\
                M & & M
            }
            \qquad
            \begin{gathered}
                B=f\circ A,\\
                \widehat{B}=\hat{f}\circ\widehat{A},\\
                d_{C^1}(f,\hat{f})\ll1.
            \end{gathered}
        \]
        The generating identities give
        \[
            \D\widehat{Z}-P\cdot\D\widehat{Q}
            =\D\widehat{W}+\widehat{Q}\cdot \D{P}
            =(\partial_z\widehat{W})(\D{z}-\widehat{H}\cdot \D{q}),
        \]
        so the resulting map is contact and coorientation-preserving. It tends to $f$ in $C^1$ and is a global diffeomorphism for sufficiently small $\eta$.

        If $C\subset U$ is a compact region inside the source Darboux chart, choose $\beta=1$ on a neighborhood of $A^{-1}(C)$. For small $\eta$, the positive margin and convergence of the inverse maps ensure $C\subset\widehat{A}(\{\beta>0\})$. Thus $\hat{f}$ is smooth near $C$.
        
        Existing smoothness of $f$ is also preserved. Indeed, for a source point $y$ where $f$ is smooth, write $v=\widehat{A}^{-1}(y)$. If $\beta(v)>0$, the new map is smooth near $y$. If $\beta(v)=0$, then $\widehat{A}(v)=A(v)=y$, and the original mixed-coordinate maps and $W$ are smooth near this point of $\Omega$. There the expression $S(v,t)=\int W(v-t\zeta)\rho(\zeta)\,d\zeta$ extends smoothly through $t=0$. Its composition with $t=\delta(v)$ is smooth near that point, so $\hat{f}$ is again smooth near $y$.

        \emph{Step 5. A finite covering and successive smoothing.}
        If $K=\varnothing$, take $\tilde{f}:=f$. Otherwise, Step 1 gives finitely many bounded open sets $U_i',U_i$, $1\le i\le N$, with $C_i:=\overline{U_i'}\subset U_i$ and $K\subset\bigcup_{i=1}^N U_i'$, such that each $\overline{U_i}$ is compactly contained in a source Darboux chart inside $O$, its image under $f$ is compactly contained in a target Darboux chart, and the mixed-coordinate projection is a diffeomorphism near $\overline{U_i}$. Choose $\varepsilon_{\mathrm{chart}}>0$ so that these chart conditions hold whenever $d_{C^1}(g,f)<\varepsilon_{\mathrm{chart}}$.

        Starting with $f_0:=f$, apply Step 4 successively in $U_i$. For the inverse mixed-coordinate map $A_i$ of $f_{i-1}$, choose a smooth cutoff $\beta_i$ compactly supported in $A_i^{-1}(U_i)$ and equal to one near $A_i^{-1}(C_i)$. A sufficiently small smoothing parameter $\eta_i>0$ gives a contactomorphism $f_i$ that is smooth near $C_i$, equals $f_{i-1}$ outside $U_i$, and satisfies $d_{C^1}(f_i,f_{i-1})<{\min\{\varepsilon_0,\varepsilon_{\mathrm{chart}}\}}/{2N}$. For every $1\le j\le N$, the triangle inequality gives
        \[
            d_{C^1}(f_j,f)
                \le\sum_{i=1}^j d_{C^1}(f_i,f_{i-1})
                <\frac{j}{2N}\min\{\varepsilon_0,\varepsilon_{\mathrm{chart}}\},
        \]
        so the chart conditions remain valid throughout the construction.

        Each modification preserves existing smoothness by Step 4, so induction gives smoothness near $C_1\cup\cdots\cup C_i$. Consequently, $\tilde{f}:=f_N$ is smooth near $K$, satisfies $d_{C^1}(\tilde{f},f)<\varepsilon_0$, and agrees with $f$ on $M\setminus O$. Taking $K=O=M$ proves the global approximation statement.
    \end{proof}

    \begin{remark}
        The generating function $W$ encodes a finite local contactomorphism through the mixed-coordinate reconstruction. A smooth time-dependent contact Hamiltonian $h_s$ instead specifies its infinitesimal motion through the vector field $X_{h_s}$ defined by $\alpha(X_{h_s})=h_s$ and $\iota_{X_{h_s}}\D\alpha=(R_\alpha h_s)\alpha-\D h_s$. Integrating this vector field gives a contact isotopy; see Geiges~\cite[Section~2.3]{Geiges2006}.
    \end{remark}
\section{A contact-open set and its based maps}\label{section:based-maps}
    This section constructs a nonempty contact-open set $\mca{U}$ with a persistent blender configuration and proves that its smooth based maps are $C^1$-dense. Relative contact approximation and the geometric strip argument of Bonatti--D\'iaz yield smooth maps with stable endpoints and localized forward orbits, after which a contact perturbation in a fixed Darboux box produces the strong connection.
\subsection{The initial model and the persistent configuration}\label{subsection:persistent-contact-configuration}
    \medskip\noindent{\bf Initial Reeb-flow setup}

    Let $(M^{2d+1},\xi)$, $d\ge1$, be a closed connected cooriented contact manifold, and fix a smooth contact form $\alpha_0$ with $\ker\alpha_0=\xi$. Assume that its Reeb flow $\varphi^t$ is Anosov and that its strong stable and strong unstable foliations are $C^\infty$. In particular, their tangent bundles are smooth. The invariant splitting is
    \[
        TM=E^{ss}_{\varphi}\oplus\R R_{\alpha_0}\oplus E^{uu}_{\varphi},
        \qquad \xi=E^{ss}_{\varphi}\oplus E^{uu}_{\varphi},
        \qquad \dim E^{ss}_{\varphi}=\dim E^{uu}_{\varphi}=d.
    \]
    The identity $\xi=E^{ss}_{\varphi}\oplus E^{uu}_{\varphi}$ follows from invariance of $\alpha_0$ and contraction in the corresponding time direction. Similarly, invariance of $\D\alpha_0$ makes both strong bundles isotropic in $\xi$, and each is Lagrangian because together they span $\xi$. Choose a closed Reeb orbit $\Gamma$ of period $T>0$ and set $\Phi_0:=\varphi^T$. Then $\Phi_0^*\alpha_0=\alpha_0$, and $\Gamma$ is a normally hyperbolic circle of fixed points of $\Phi_0$. Throughout this section, $\Cont^1(M,\xi)$ denotes coorientation-preserving $C^1$ contactomorphisms with the $C^1$ topology.

    \medskip\noindent\textbf{A concrete realization}

    Take a closed oriented surface $(S,h)$ of constant curvature $-1$ and let
    \[
       M=S^*S:=\{(q,p)\in T^*S:|p|_h=1\},\qquad
       \alpha_0=\lambda_{\mathrm{can}}|_{S^*S},
    \]
    where $\lambda_{\mathrm{can}}=p\,dq$ is the tautological (Liouville) one-form. The musical isomorphism $(q,v)\mapsto(q,h_q(v,\cdot))$ identifies $S^*S$ with the unit tangent bundle; under that identification $\alpha_0(\zeta)=h_q(v,D\pi(\zeta))$. Its Reeb flow is the unit-speed geodesic flow \cite[Section~1.5]{Geiges2008}, which is Anosov; for constant negative curvature its strong stable and unstable foliations are the smooth horocycle foliations (see \cite[Chapter~17]{KatokHasselblatt1995}). A closed geodesic provides $\Gamma$ and its period $T$. Thus the initial setup is realized with $\dim M=3$ and $d=1$.

    \medskip\noindent{\bf Preparation of the blender model}

    We use the preparation in the proof of \cite[Theorem 4.4]{Chaidez2025} to obtain the configuration of \cite[Setup 4.1]{Chaidez2025}. This preparation may replace $\alpha$ by a smooth positive multiple $\alpha=u\alpha_0$, equal to $\alpha_0$ near $\Gamma$, and use the time-$T$ map $\Phi$ of its Reeb flow. The modification, which adjusts the period of an auxiliary closed orbit, can be chosen arbitrarily small. The selected phases $P,Q\in\Gamma$ and their heteroclinic connection are provided by \cite[Lemma 4.5]{Chaidez2025}, which is essentially \cite[Lemma 4.3]{BonattiDiaz1996}. In particular, $\Phi$ is $C^1$-close to $\Phi_0$ and is strict with respect to $\alpha$. We fix the resulting configuration, which has the following properties.
    \begin{enumerate}[(a)]
        \item $\Phi$ is partially hyperbolic, with
        \[
            TM=E^{ss}_{\Phi}\oplus E^c_{\Phi}\oplus E^{uu}_{\Phi},\qquad
            \xi=E^{ss}_{\Phi}\oplus E^{uu}_{\Phi},\qquad
            E^c_{\Phi}=\R R_{\alpha}.
        \]
        \item $\Gamma$ is a closed Reeb orbit for $\alpha$ and a normally hyperbolic circle fixed pointwise by $\Phi$.
        \item On a neighborhood $V$ of $\Gamma$, there is a smooth integrable rank-$d$ subbundle $E^{ss}_{\rm sm}\subset\xi$, uniformly contracted by $\Phi$ where defined, which agrees with $E^{ss}_{\Phi}$ along the local stable manifold of $\Gamma$. The local smooth stable extension in \cite[Setup 4.1(c)]{Chaidez2025} suffices here without requiring global smoothness of the strong foliations of the modified flow.
        \item The distinct fixed points $P,Q\in\Gamma$ have a strong heteroclinic orbit from $P$ to $Q$. At a point $a$ of this orbit,
        \[
            a\in W^{uu}(P,\Phi)\cap W^{ss}(Q,\Phi),\qquad
            T_{a}M=T_{a}W^{uu}(P,\Phi)\oplus\R R_{\alpha}(a)\oplus T_{a}W^{ss}(Q,\Phi).
        \]
    \end{enumerate}
    The forms $\alpha$ and $\alpha_0$ define the same contact structure $\xi$. We keep $\alpha$ as the background form in the subsequent conformal-factor identities.

    For this configuration, \cite[Theorem 4.3]{Chaidez2025} provides a smooth family $\Psi_r$, with $\Psi_0=\Phi$, such that $\Psi_r$ has a contact blender for sufficiently small positive $r$. Fix one such parameter $r_*$ and put $f_*:=\Psi_{r_*}$. We write $p_*,q_*$ for the points denoted by $P,Q$ in that construction, with $p_*$ center-contracting and $q_*$ center-expanding. Our indices are unstable indices, so 
    \[
        \ind_u(p_*)=d,\qquad \ind_u(q_*)=d+1.
    \]
    This convention differs from the stable-index convention used in the corresponding statements of \cite{Chaidez2025}.

    More precisely, \cite[Theorem 4.3]{Chaidez2025} provides a blender box $B_*$ for an iterate $f_*^\ell$, with $\ell\ge1$, and a vertical disk of $W^u(p_*,f_*)$ in its superposition region. Here $W^u(p_*)=W^{uu}(p_*)$. The transverse connection from $q_*$ to $p_*$ and the transverse homoclinic points are provided by \cite[Lemma 5.14]{Chaidez2025}. Choose a point
    \[
        b_*\in W^u(q_*,f_*)\pitchfork W^s(p_*,f_*).
    \]

    The boundary condition for the long return deserves an explicit comparison with the cited blender construction. Let $C$ be the blender box, let $F$ be its defining map, and let $A'$ be the selected component of $C\cap F^m(C)$, with input domain $B=F^{-m}(A')$. These correspond to $U$, $f$, $f^N(B)$, and $B$ in Definition~\ref{def:dynamical blender}, with $m=N$. Injectivity gives $[A'\cap F^m(\partial^u C)=\varnothing]\Longleftrightarrow[B\cap\partial^u C=\varnothing]$, whereas $[A'\cap F(\partial^u C)=\varnothing]\Longleftrightarrow[F^{m-1}(B)\cap\partial^u C=\varnothing]$. The second disjointness condition is the one printed in \cite[p.~365, condition (H2)]{BonattiDiaz1996} and \cite[Definition~3.22(b) and Lemma~6.10]{Chaidez2025}, while our branchwise criterion uses the first. They concern different stages of the return and are not interchangeable.

    Separation of the input domain from $\partial^u C$ ensures that a continued crossing disk cannot terminate at the unstable boundary of its input disk. Together with the other boundary avoidance conditions, it allows the selected returned strip to contain both original end disks, so that center-unstable expansion yields the width inequality used in \cite[Lemma~1.7]{BonattiDiaz1996}. The subsequent orbit-selection argument uses this selected-strip conclusion and its two possible return times. The whole-component separation provides a sufficient geometric condition for that conclusion. It is also verified for the concrete construction in \cite[p.~377, proof of Lemma~2.11]{BonattiDiaz1996}.

    For the contact model, this separation follows from the estimates already established in \cite{Chaidez2025}. Write $F_r=\Psi_r^\ell$, where $\ell$ is denoted by $N$ in that paper, and write $m_0$ for the fixed excursion integer of Construction~5.5, reserving $m_r$ for the long-return integer of Notation~5.15. The unstable radius of $C_r=B_r(Q)$ is $\rho_r=l\mu^{-m_r}$. For $z\in A'_r$, put $y=F_r^{-m_r}(z)$ and $w=F_r^{-m_0}(z)$. By \cite[Lemma~6.14]{Chaidez2025}, $w$ lies in a fixed neighborhood $W$ and every intermediate iterate from $y$ to $w$ stays in the coordinate chart. Choose $K$ with $|u|\le K$ on $W$. Iterating the coordinate expansion in \cite[Lemma~5.10]{Chaidez2025} gives
    \[
        |u(F_r^{-m_r}(z))|
        \le K\mu^{-(m_r-m_0)}
        =\frac{K\mu^{m_0}}{l}\,\rho_r.
    \]
    The constants $\ell,\mu,m_0,W$ are fixed before $l$, and $m_r$ is independent of $l$. Construction~5.17 requires only a lower bound on $l$, so we choose $l>2K\mu^{m_0}$, together with that bound, and then take $r_*$ sufficiently small. Thus $F_{r_*}^{-m_{r_*}}(A'_{r_*})\subset\{|u|<\rho_{r_*}/2\}$, which establishes the long-return separation used in our criterion.

    After fixing $r_*$ and the return time, this strict separation persists for the component selected by the marked return point. The construction in \cite[proof of Lemma~6.5]{Chaidez2025}, together with the strict bound above, places this point and its preimage in the interiors of the output and input boxes. Keeping its preimage fixed therefore gives a continued marker for nearby maps. Any Hausdorff limit of the selected compact connected components is contained in the original intersection and contains the limiting marker, hence lies in the original selected component. Uniform convergence of the inverse return maps then preserves separation from the input unstable boundary. We keep these choices in the persistent configuration below.

    \begin{proposition}[Persistent geometric configuration]\label{prop:persistent-contact-data}
        There is a nonempty $C^1$-open neighborhood $\mca{U}$ of $f_*$ in $\Cont^1(M,\xi)$ such that every $f\in\mca{U}$ has the following properties.
        \begin{enumerate}[(i)]
            \item There is a partially hyperbolic splitting $TM=E_f^{ss}\oplus E_f^c\oplus E_f^{uu}$, with the one-dimensional bundle $E_f^c$ uniformly transverse to $\xi$.
            \item The fixed points $p_*,q_*$ have hyperbolic continuations $p_f,q_f$. For a constant $\eta>0$ independent of $f\in\mca{U}$,
            \[
                \zeta_c(p_f;f)<-\eta<0<\eta<\zeta_c(q_f;f),
            \]
            where $\zeta_c=\log\lambda_c$ is the logarithm of the positive center multiplier at a fixed point.
            \item There is a transverse connection $W^u(q_f,f)\pitchfork W^s(p_f,f)\neq\varnothing$.
            \item The blender configuration for $f_*^\ell$ persists for $f^\ell$, together with a disk of $W^{uu}(p_f,f)$ in the corresponding superposition region.
        \end{enumerate}
    \end{proposition}
    \begin{proof}
        The construction starts from a partially hyperbolic time map whose center is transverse to $\xi$, and the model $f_*$ can be chosen sufficiently close to this map. Persistence of partial hyperbolicity and continuity of its bundles then establish (i), including uniform transversality after shrinking the neighborhood $\mca{U}$. Hyperbolic continuation and continuity of the center multipliers imply (ii). These multipliers are positive because the coorientation is preserved and the invariant center line is transverse to $\xi$. A compact transverse intersection near $b_*$ persists, establishing (iii). Finally, the blender and the incoming disk provided by the construction persist under sufficiently small $C^1$ perturbations. Continuity of $f\mapsto f^\ell$ gives (iv) for the same integer $\ell$.
    \end{proof}

    We use the persistent partially hyperbolic splitting from Proposition~\ref{prop:persistent-contact-data} throughout. The cone invariance and expansion required for the vertical-strip argument are included in the geometric blender provided by \cite[Theorem 4.3]{Chaidez2025} and are used through \cite[Lemmas 1.4--1.8]{BonattiDiaz1996}.

    More precisely, we use the concrete return structure in Proposition \ref{prop:persistent-contact-data}(iv). For the fixed parameter $r_*$, Construction 5.17 of \cite{Chaidez2025} gives a compact box $C=B_{r_*}(Q)$ with $p_*\notin C$. Indeed, the center coordinate on $C$ satisfies $|t|\le r_*(1-e^{-8\ell r_*})<L$, whereas $t(p_*)=L$. There are two return branches, of lengths $1$ and $m$ for $F=f^\ell$, where $m$ is fixed after choosing $r_*$. After shrinking $\mca{U}$, the same box and return lengths satisfy the blender properties for every $f\in\mca{U}$, and $p_f\notin C$. This is the persistence of the geometric blender axioms; see \cite[Lemma 1.11]{BonattiDiaz1996} and \cite[Definition 3.22 and Remark 3.23]{Chaidez2025}. The primary branch defines the reference local stable disk $W_0^s(q_f)$, whose forward $F$-iterates stay in $C$ and converge to $q_f$; see \cite[Lemma 1.6]{BonattiDiaz1996}. The incoming vertical disk is a disk of $W^{uu}(p_f,f)$ to the right of $W_0^s(q_f)$.
\subsection{A contact strong connection}\label{subsection:contact-blender-connection}
    We use the orbit selection below to apply the point-moving Lemma \ref{lemma:point-moving}, with all perturbation supports contained in fixed neighborhoods chosen before the displacement is made small.

    For each $f\in\mca{U}\cap C^\infty$, choose one transverse heteroclinic point $b_f$ as in Proposition \ref{prop:persistent-contact-data}(iii), and put
    \[
        K_f^{\pitchfork}:=\{p_f,q_f\}\cup\Orb(b_f,f).
    \]
    This is a compact invariant set. The following lemma verifies the orbit-selection property for the chosen contact-open set.

    \begin{lemma}[Orbit selection through the blender]\label{lemma:blender-orbit-selection}
        Suppose $f\in\mca{U}\cap C^\infty$ satisfies $W^{uu}(p_f,f)\cap W^{ss}(q_f,f)=\varnothing$. Then there exist a nonperiodic point
        $x_*\in W^{uu}(p_f,f)$, a Darboux neighborhood $V_*$ of $x_*$, and points
        $y_j\in W^{ss}(q_f,f)$ converging to $x_*$, such that $\overline{V_*}\cap K_f^{\pitchfork}=\varnothing$, $f^{-i}(x_*)\notin\overline{V_*}\,(i\ge1)$, and, for all sufficiently large $j$, $f^i(y_j)\notin\overline{V_*}\,(i\ge1)$. The neighborhood $V_*$ is fixed as $j\to\infty$.
    \end{lemma}
    \begin{proof}
        Fix $f$ as in the statement, and suppress its subscript on $p,q$. Write $F=f^\ell$ and keep the box $C$, the return length $m$, and the local stable disk $W_{\loc}^{ss}(q,f)$ specified above. The preceding verification, including the long-return input-boundary separation, gives the branchwise geometric criterion of Definition~\ref{def:dynamical blender}. We use its selected-strip conclusion and the width and thickening estimates of \cite[Lemmas~1.4--1.8]{BonattiDiaz1996}. The cone comparison is \cite[Remark~3.24]{Chaidez2025}.

        \emph{Step 1. Stable endpoints and their return itineraries.}
        Let $D\subset W^{uu}(p,f)\cap C$ be an incoming vertical disk to the right of $W_{\loc}^s(q,f)$. After applying one of the blender branches and replacing $D$ by its selected vertical subdisk, we may assume that $D$ avoids both center boundary faces of $C$. The primary branch produces a disk separated from the right face, and the other branch lands in a component disjoint from that face; every right-hand vertical disk avoids a neighborhood of the left face. The resulting disk still belongs to $W^{uu}(p,f)$.

        By \cite[Lemma 1.5]{BonattiDiaz1996}, extend $D$ to a vertical strip and take successively thinner closed substrips $S_j$ about $D$. Thus $S_j$ are $(d+1)$-dimensional vertical strips through $C$, all to the right of $W_{\loc}^{ss}(q,f)$, and $d_H(S_j,D)\to0$ as $j\to\infty$. Here $d_H$ is Hausdorff distance.

        We use the localization argument in the proof of \cite[Lemma 1.8]{BonattiDiaz1996}. Starting with $S_j$, Lemma 1.7 of that paper either gives an intersection with $W_{\loc}^{ss}(q,f)$ after one of the two branches, or selects another right-hand strip in $C$ whose width is at least $\lambda>1$ times the previous width. In the latter case the selected strip lies in the image of the previous strip under $F$ or $F^m$, as in the component construction in that lemma. Since widths of vertical strips in $C$ have a common upper bound by Lemma 1.4 \cite{BonattiDiaz1996}, this iteration terminates after finitely many steps. Pulling the final intersection back through the selected branches produces $z_j\in S_j$ and iteration times
        \[
            0=t_{j,0}<t_{j,1}<\cdots<t_{j,k_j},
            \qquad t_{j,a+1}-t_{j,a}\in\{\ell,\ell m\},
        \]
        with $f^{t_{j,a}}(z_j)\in C$ for every $0\le a\le k_j$ and $f^{t_{j,k_j}}(z_j)\in W_{\loc}^{ss}(q,f)$. Consequently
        \begin{equation}\label{eq:blender-endpoint-tower}
            f^n(z_j)\in T_f:=\bigcup_{a=0}^{\ell m}f^a(C)
            \quad(n\ge0),\qquad z_j\in W^{ss}(q,f).
        \end{equation}
        Indeed, each intermediate iterate before the final intersection lies on a branch of length at most $\ell m$ starting in $C$. After that intersection, successive $\ell$-step iterates remain in $C$. The compact set $T_f$ avoids $p$, because if $p\in f^a(C)$, the identity $f^{-a}(p)=p$ would imply $p\in C$. We use $T_f$ to contain the selected forward orbits.

        \begin{figure}[H]
            \centering
            \includegraphics[scale=0.6]{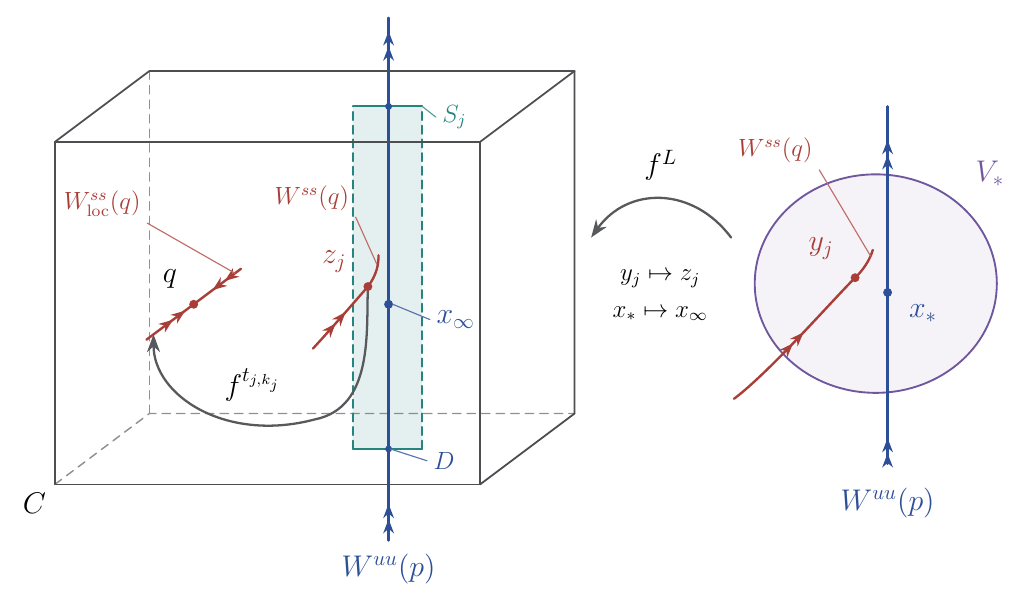}
            \caption[Orbit selection and a fixed perturbation domain]{\textbf{Orbit selection and a fixed perturbation domain.}
                Thin strips $S_j$ around $D\subset W^{uu}(p)$ yield points $z_j\to x_\infty\in D$ with $f^{t_{j,k_j}}(z_j)\in W_{\loc}^{ss}(q)$. For a fixed $L$ and sufficiently large $j$, the pullbacks $y_j=f^{-L}(z_j)\to x_*=f^{-L}(x_\infty)$ lie in a fixed Darboux neighborhood $V_*$ with the orbit-avoidance properties of Lemma~\ref{lemma:blender-orbit-selection}.}
            \label{fig:strong-connection-domain}
        \end{figure}

        Compactness and $d_H(S_j,D)\to0$ give, after passing to a subsequence, $z_j\to x_\infty\in D$. This point is distinct from $p$ and belongs to $W^{uu}(p,f)$, so its backward orbit converges to $p$ and it is nonperiodic.
        
        \emph{Step 2. A fixed perturbation domain.}
        Since $T_f$ is compact and avoids $p$, choose $L\ge1$ sufficiently large that
        \[
            x_*:=f^{-L}(x_\infty)\notin T_f,
            \qquad y_j:=f^{-L}(z_j)\longrightarrow x_*.
        \]
        The points $y_j$ belong to $W^{ss}(q,f)$. Consider the compact set $A_-=\{p\}\cup\{f^{-i}(x_*):i\ge1\}$ and the finite set $A_+=\{f^i(x_*):1\le i\le L\}$. By nonperiodicity, neither set contains $x_*$. Moreover, $x_*\notin K_f^{\pitchfork}$, since it differs from both fixed points and its backward limit is $p$, whereas every point of $\Orb(b_f,f)$ has backward limit $q$. Thus
        \[
            x_*\notin T_f\cup K_f^{\pitchfork}\cup A_-\cup A_+.
        \]
        Choose a Darboux neighborhood $V_*$ of $x_*$ whose closure avoids this compact union. This choice gives the required avoidance conditions for the backward orbit and the marked transverse connection. For sufficiently large $j$, continuity of the finitely many iterates $1\le i\le L$, together with the choice of $x_*$, ensures that $f^i(y_j)\notin\overline{V_*}$. For $i\ge L$, equation \eqref{eq:blender-endpoint-tower} gives $f^i(y_j)=f^{i-L}(z_j)\in T_f$, hence the remaining forward-orbit avoidance conditions. The choice of $V_*$ is independent of $j$, as required.
    \end{proof}

    \begin{proposition}[Contact realization of the strong connection]\label{prop:contact-blender-connection}
        For every $f\in\mca{U}\cap C^\infty$ and every $C^1$-neighborhood $\mca{V}\subset\mca{U}$ of $f$ in $\Cont^1(M,\xi)$, there exists $g\in\mca{V}\cap C^\infty$ with
        \[
            W^{uu}(p_g,g)\cap W^{ss}(q_g,g)\neq\varnothing.
        \]
        The transverse connection from $q_g$ to $p_g$ and the blender configuration remain persistent.
    \end{proposition}
    \begin{proof}
        If the required intersection already exists, take the smooth map $g=f$. Otherwise apply Lemmas \ref{lemma:blender-orbit-selection} and \ref{lemma:point-moving} to obtain smooth contactomorphisms $h_j$ supported in $V_*$ with $h_j(x_*)=y_j$ and $h_j\to\id$ in $C^1$. Put $g_j=h_j\circ f$. Each $g_j$ is smooth, and for large $j$, $g_j\in\mca{V}$. Furthermore, $g_j$ agrees with $f$ outside $f\inv(V_*)$. We omit the subscript $j$ for $g_j$.

        Let $z_*=f^{-1}(x_*)$. The orbit avoidance conditions established in Lemma~\ref{lemma:blender-orbit-selection} imply
        \[
            g^{-i}(z_*)=f^{-i}(z_*)\quad(i\ge0),\qquad
            g^i(z_*)=f^{i-1}(y_j)\quad(i\ge1).
        \]
        The first orbit converges to $p_f$ along its backward iterates and the second to $q_f$ along its forward iterates. Since the support avoids neighborhoods of $p_f,q_f$, these fixed points coincide with the original ones, namely $p_g=p_f$ and $q_g=q_f$. Thus $z_*$ is the desired strong heteroclinic point for $g$. Its intersection is quasi-transverse, since the strong tangent spaces belong to the transverse bundles $E_{g}^{uu}$ and $E_{g}^{ss}$. By $C^1$-smallness, both the transverse intersection and the blender configuration persist for large $j$.
    \end{proof}
\subsection{Based maps and their density}\label{subsection:Based-maps}
    \begin{definition}\label{def:based-maps}
        A \emph{based map} is a map in the set
        \[
            \mca{B}=\set{g\in\mca{U}\cap C^\infty:W^{uu}(p_g,g)\cap W^{ss}(q_g,g)\neq\varnothing
            }.
        \]    
    \end{definition}
    
    Every based map has a co-index-one cycle, because the transverse connection from $q_g$ to $p_g$ is already present throughout $\mca{U}$.

    \begin{proposition}[Density of the based maps]\label{prop:contact heterodim cycle}
        The smooth based maps are $C^1$-dense in the chosen contact-open set $\mca{U}$, namely, $\clos_{C^1}(\mca{B})\cap\mca{U}=\mca{U}$.
    \end{proposition}
    \begin{proof}
        Let $f_0\in\mca{U}$ and let $\mca{V}_0\subset\mca{U}$ be any $C^1$-open neighborhood of $f_0$. Proposition~\ref{prop:relative-contact-approximation}, with $K=O=M$, gives $f_\infty\in\mca{V}_0\cap C^\infty$. Choose a $C^1$-open neighborhood $\mca{V}_\infty$ of $f_\infty$ contained in $\mca{V}_0$. Apply Proposition~\ref{prop:contact-blender-connection} to $f_\infty$ and $\mca{V}_\infty$ to obtain a smooth map $g\in\mca{V}_\infty$ with the required strong connection. Then $g\in\mca{V}_0\cap\mca{B}$. Since $f_0$ and $\mca{V}_0$ were arbitrary, the density follows.
    \end{proof}

    Every based map is smooth, whereas the ambient open set $\mca{U}$ consists of $C^1$ contactomorphisms. Smoothness is used in Section~\ref{subsection:flatten-center-germ}, while $C^1$-regularity suffices for the other parts of the argument.

    Fix a based map $g$, and suppress the subscript $g$ on its saddles. Choose
    \[
        a\in W^{uu}(p,g)\cap W^{ss}(q,g),\qquad
        b\in W^u(q,g)\pitchfork W^s(p,g).
    \]
    The compact invariant cycle used for subsequent perturbations is $K_g:=\{p,q\}\cup\Orb(a,g)\cup\Orb(b,g)$. The choices of $a$ and $b$ need not depend continuously on $g$.

    Fix disjoint saddle neighborhoods $U_p$, $U_q$ and move the marked points along their orbits so that $a\in U_p$ and $b\in U_q$. Choose integers $T_a,T_b\ge1$ such that $a'=g^{T_a}(a)\in U_q$ and $b'=g^{T_b}(b)\in U_p$. On small neighborhoods of these points, the iterates define smooth contact transition maps
    \[
        A=g^{T_a}:U_a\longrightarrow U_{a'},\qquad
        B=g^{T_b}:U_b\longrightarrow U_{b'}.
    \]
    The transition neighborhoods can be chosen with closures compactly contained in coordinate neighborhoods. After shrinking, the neighborhoods of any fixed finite collection of distinct marked orbit points can be made pairwise disjoint.

    \begin{remark}[Choice of the transition times]\label{remark:transition-time-choice}
        The marked points and the integers $T_a,T_b$ are auxiliary choices. They may be adjusted when smaller saddle plaques or perturbation neighborhoods are selected. More explicitly, replacing $(a,a')$ by $(g^{-i}(a),g^j(a'))$, and $(b,b')$ by $(g^{-k}(b),g^\ell(b'))$ for nonnegative integers $i,j,k,\ell$ replaces the transition times by $T_a+i+j$ and $T_b+k+\ell$, respectively, and puts the endpoints closer to the corresponding saddles. We keep the same notation for the new endpoints and times. These changes preserve the heteroclinic orbits and hence $K_g$.

        This freedom is used during the preliminary choice of plaques, marked points, and boxes. Once these choices and the saddle family $F_t$ have been fixed for a construction, $T_a,T_b$ are fixed finite integers independent of $n,m,t$.
    \end{remark}

    Choose a Darboux chart at each marked point and consider the $Dg$-invariant partially hyperbolic splitting $TM=E_g^{ss}\oplus E_g^c\oplus E_g^{uu}$ from Proposition~\ref{prop:persistent-contact-data}. Its two strong bundles cannot both become constant coordinate distributions throughout an open neighborhood, because their sum is the nonintegrable contact distribution $\xi$. This obstruction concerns the bundles throughout a neighborhood and does not rule out aligning their fibers at a single point.

    Relative to the fixed global contact form $\alpha$, write $g^*\alpha=\kappa_g\alpha$. Then
    \[
        \kappa_{g^N}(x)=\prod_{j=0}^{N-1}\kappa_g(g^j(x)),
    \]
    for every $N\ge1$.
\section{A long periodic orbit near the cycle with flat center direction}\label{section:Long flat orbit near the cycle}
    The main result of this section is Theorem~\ref{thm:contact flat}, a contact analogue of the flat-germ statement in \cite[Theorem 2.3]{BonattiDiaz2008}. The proof follows the strategy of \cite{BonattiDiazFisher2008,BonattiDiaz2008}. In Sections~\ref{subsection:hyperbolic passages}--\ref{subsection:unit-center-multiplier}, we use $C^1$ estimates to construct long periodic orbits with center multiplier one. In Section~\ref{subsection:flatten-center-germ}, we then apply the contact lift to obtain an identity center germ. The construction starts from the smooth based class $\mca{B}$ defined in Section~\ref{subsection:Based-maps}, and all neighborhoods and perturbation sizes are measured in the $C^1$ topology.

    The contact condition prevents a direct application of the affine-model constructions and explicit center-eigenvalue calculations used in those works. In particular, transition maps for elements of $\mca{U}$ need not have a block form. To overcome these difficulties, we first construct long periodic orbits and estimate their center eigenvalues asymptotically. We then introduce a family of local perturbations supported near a saddle and apply the intermediate value theorem to make the center eigenvalue exactly $1$.

    More precisely, for each $g\in\mca{B}$, we construct perturbations $g_k\to g$ in $C^1$ that preserve the saddles $p,q$ and have controlled supports, together with periodic points $w_k$ of least periods $N_k\to\infty$, such that
    \[
        \mu_{w_k,g_k}\stackrel{*}{\longrightarrow}
        \theta_g\delta_p+(1-\theta_g)\delta_q,\quad\text{where}\quad
        \theta_g=\frac{\zeta_c(q;g)}{\zeta_c(q;g)-\zeta_c(p;g)}.
    \]
    Here, for a partially hyperbolic periodic point $w$ of $g$ with primitive period $N$, $\mu_{w,g}$ denotes the equidistributed measure on the periodic orbit $\Orb(w,g)$, and $\zeta_c(w;g)=\frac{1}{N}\log\lambda_c(w;g)$ is the normalized logarithm of its center eigenvalue. This notation agrees with that of Proposition~\ref{prop:persistent-contact-data}.

    Geometrically, the orbit $\Orb(w_k,g_k)$ spends approximately $\theta_g N_k$ iterates near $p$ and $(1-\theta_g)N_k$ iterates near $q$. Consequently, $\zeta_c(w_k;g_k)\to0$. Indeed, Lemma~\ref{lemma:center multiplier conformal factor} gives
    \[
        \zeta_c(w_k;g_k)=\int\log\kappa_{g_k}\,\D\mu_{w_k,g_k},
    \]
    and convergence $g_k\to g$ in $C^1$ implies uniform convergence of $\log\kappa_{g_k}$ to $\log\kappa_g$. The integrals therefore converge to $\theta_g\zeta_c(p;g)+(1-\theta_g)\zeta_c(q;g)=0$, since $\kappa_g\approx\kappa_g(\,\bigcdot\,)=\lambda_c(\,\bigcdot\,;g)$ near the saddles $\bigcdot=p,q$. Nevertheless, the center eigenvalue $\lambda_c(w_k;g_k)=\exp(N_ko_k(1))$ may be large, so a $C^1$-small modification at \emph{only} one point of the orbit might not suffice to achieve $\lambda_c=1$. We therefore modify the center eigenvalue throughout the neighborhood $U_p$, which contains approximately $\theta_gN_k$ points of the orbit.

    This correction over multiple points of the orbit is implemented through the continuous families of long periodic orbits and the uniform averaging estimates in Proposition~\ref{prop:almost flat periodic point}. The balanced selection in Proposition~\ref{prop:flat lambda_c(w)=1} then yields the required sequence with $\zeta_c(w_k;g_k)=0$ exactly. Lemma~\ref{lemma:uniform-parameter-passage} provides the passages near the saddles. Lemma~\ref{lemma:passage-transition-connection}, applied in Step~2 of the proof of Proposition~\ref{prop:almost flat periodic point}, connects these passages near $p,q$ to the transitions between them.

    \begin{figure}[H]
        \centering
        \includegraphics[scale=0.6]{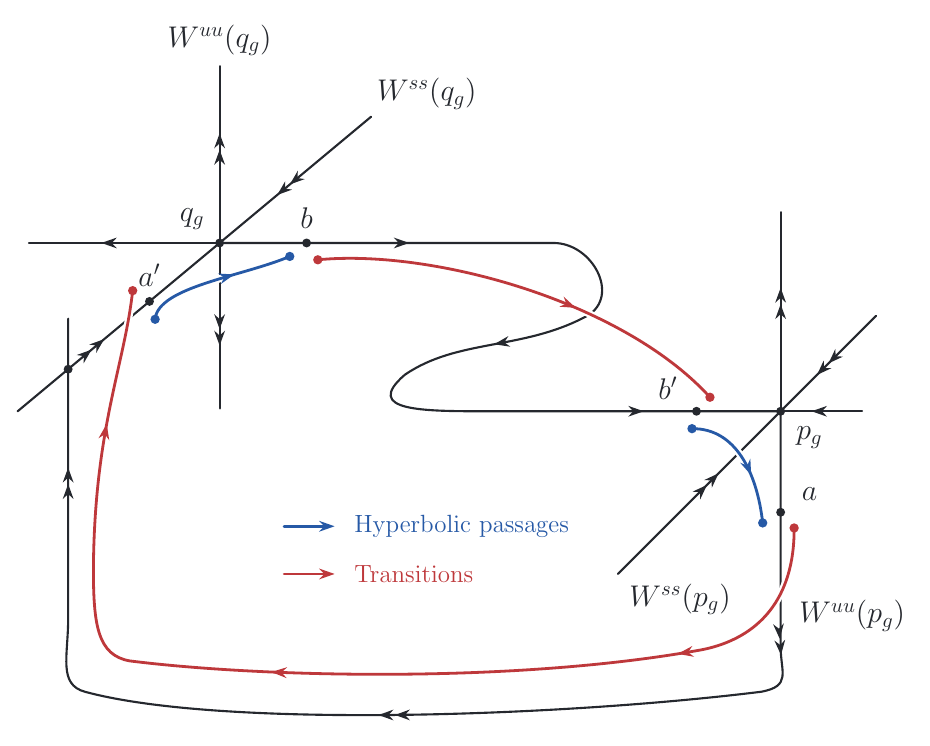}
        \caption[Two passages and two transitions]{\textbf{Two passages and two transitions.}
            Blue curves schematically represent local hyperbolic passages near $q_g$ and $p_g$, while red curves represent transitions near the heteroclinic itineraries $b\to b'$ and $a\to a'$. A $C^1$-small contact perturbation then connects the endpoints for the chosen passages and transitions.}
        \label{fig:transitions-passages}
    \end{figure}

    \begin{definition}
        Let $p$ be a periodic point for a partially hyperbolic diffeomorphism $f$, with period $N$. Then $T_pM$ splits into $E^{ss}(p)\oplus E^c(p)\oplus E^{uu}(p)$. We say $p$ is a \emph{saddle-node} if $Df^{N}|_{E^c(p)}=\id$.
    \end{definition} 
    
    \begin{theorem}\label{thm:contact flat}
        Let $g\in\mca{B}$. For every integer $r\ge1$, every $C^1$-neighborhood $\mca{V}$ of $g$, and every $N_*>0$, there exist
        \[
            G\in\mca{V}\cap\mca{U}\cap\Cont^r(M,\xi)
        \]
        and a periodic point $w$ for $G$, of least period $N>N_*$, such that there is a $C^r$ curve $C$ through $w$, with $T_wC=E^c(w,G)$, on which $G^N$ is the identity as a germ at $w$.
        In particular, $w$ is a saddle-node with an identity center germ.
    \end{theorem}
\subsection{Hyperbolic passages}\label{subsection:hyperbolic passages}
    A \emph{hyperbolic passage} is a finite orbit segment entering near a prescribed local stable point, remaining near the saddle, and leaving near a prescribed local unstable point. Let $W^s_{\loc}(p,g;V)$ denote the points whose forward orbits stay in $V$ and converge to $p$, while the unstable set $\lssp{W}{\loc}{u}{p,g;V}$ is defined using backward orbits. Let us consider the following toy model.

    \begin{example}
        Let $g=(S,U)$ be a hyperbolic map on $\rn{d_s}\times\rn{d_u}$ with coordinates $(s,u)$ and with fixed point $p=(0,0)$. Suppose $S=S(s)$ and $U=U(u)$ are uniformly contracting and expanding, respectively. Then $\lsp{W}{s}{p}=\rn{d_s}\times\set{0}$ and $\lsp{W}{u}{p}=\set{0}\times\rn{d_u}$. Use the maximum product norm and let $V=\ball{\rho}{p}$. For any $x=(x_s,0)\in\lssp{W}{\loc}{s}{p;V}$ and any $y=(0,y_u)\in\lssp{W}{\loc}{u}{p;V}$ there are points $(z_n)$ such that 
        \[
            z_n\to x,\qquad g^n(z_n)\to y,\qquad \tset{g^k(z_n):0\le k\le n}\subset V.
        \]
    \end{example}
    \begin{proof}
        It suffices to put $z_n=(x_s,U^{-n}(y_u))$. Since the pair $(S,U)$ is uniformly hyperbolic, $z_n\to(x_s,0)$ and $g^n(z_n)=(S^n(x_s),y_u)\to(0,y_u)$ as $n\to\infty$. Moreover, each $g^k(z_n)$ lies in $V$ by the construction.
    \end{proof}

    \begin{figure}[H]
        \centering
        \includegraphics[scale=0.6]{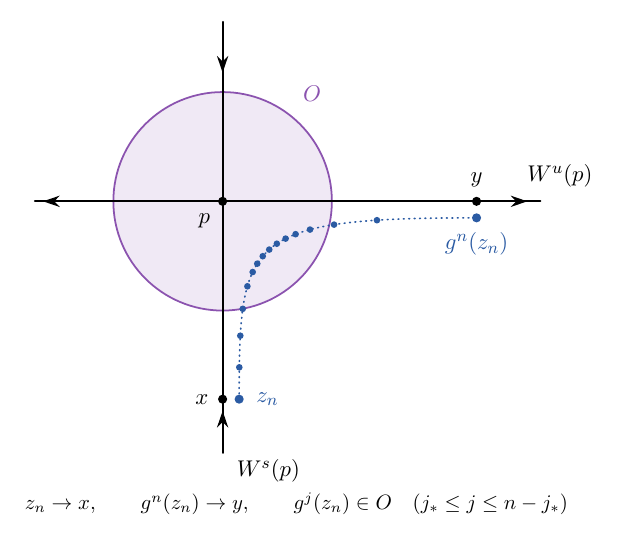}
        \caption[Passage through a hyperbolic saddle]{\textbf{Passage through a hyperbolic saddle.}
            The blue points depict a discrete orbit segment with $z_n\to x\in W^s_{\loc}(p)$ and $g^n(z_n)\to y\in W^u_{\loc}(p)$. For each neighborhood $O\subset U_p$ of $p$, all middle iterates $g^j(z_n)$, $j_*\le j\le n-j_*$, lie in $O$ for sufficiently large $n$, with $j_*$ independent of $n$.}
        \label{fig:hyperbolic-passage}
    \end{figure}

    For a continuous one-parameter family $(F_t)_{t\in[-\delta_0,\delta_0]}$ of hyperbolic maps with a common hyperbolic fixed point $p$, Lemma~\ref{lemma:uniform-parameter-passage} establishes an intrinsic version on the ambient manifold and constructs passages that depend continuously on $t$ and satisfy uniform estimates for the equidistributed measure. We then apply the passage construction at $p,q$ and identify the center multiplier with the conformal cocycle. These are the dynamical ingredients for the local perturbations and the periodic-orbit construction in the next two subsections.
    
    \begin{lemma}[Uniform saddle passage for a parameter family]\label{lemma:uniform-parameter-passage}
        Let $t\mapsto F_t$ be a continuous map into the space $\Diff^1(M)$, defined near $t=0$, with a common hyperbolic fixed point $p$ for each $F_t$. After restricting to a sufficiently small interval $I=[-\delta_0,\delta_0]$, there are local stable and unstable plaques of uniform size with the following property: If
        \[
            e_t\in W^s_{\loc}(p,F_t),\qquad d_t\in W^u_{\loc}(p,F_t)
        \]
        are continuous choices in sufficiently small compact subplaques, then for every sufficiently large $n$ there is a continuous choice $x_n(t)$ whose orbit segment $[F_t^j(x_n(t))]_{j=0}^n$ stays in a fixed neighborhood of $p$ and has the following properties.
        \begin{enumerate}[(i)]
            \item Uniformly for $t\in I$,
            \[
                x_n(t)\longrightarrow e_t,\qquad F_t^n(x_n(t))\longrightarrow d_t.
            \]
            In particular, for every fixed $i\ge0$, we have the uniform entrance and exit limits as $n\to\infty$
            \begin{equation}\label{eq:parameter-passage-ends}
                F_t^i(x_n(t))\longrightarrow F_t^i(e_t),\qquad 
                F_t^{n-i}(x_n(t))\longrightarrow F_t^{-i}(d_t).
            \end{equation}
            \item For every neighborhood $O$ of $p$, there is $j_*=j_*(O)$, independent of $n,t$, such that
            \begin{equation}\label{eq:parameter-passage-middle}
                F_t^j(x_n(t))\in O\qquad(j_*\le j\le n-j_*).
            \end{equation}
            In particular, for every continuous function $\varphi:I\times M\to\mathbb R$,
            \begin{equation}\label{eq:parameter-passage-average}
                \sup_{t\in I}\bigg|
                    \frac{1}{n}\sum_{j=0}^{n-1}\varphi(t,F_t^j(x_n(t)))-\varphi(t,p)
                \bigg|\longrightarrow0.
            \end{equation}
        \end{enumerate}
    \end{lemma}
    \begin{proof}
        \emph{Step 1. Uniform coordinates.}
        Let $d=d_s+d_u$. Fix a smooth ambient chart $\chi:(M,p)\to(\rn{d},0)$ with $\chi(p)=0$. Put $E^s=E^s(p,F_0)$ and $E^u=E^u(p,F_0)$, viewed as fixed vector spaces. The hyperbolic subspaces $E^s(p,F_t)$ and $E^u(p,F_t)$ have constant dimensions $d_s,d_u$, respectively, and vary continuously with $t$. Choose continuous linear isomorphisms $L_t:T_pM\to T_0\rn{d}$ sending the two factors $E^s$ and $E^u$ onto $D\chi(p)E^s(p,F_t)$ and $D\chi(p)E^u(p,F_t)$, respectively. The charts $\Theta_t=\chi^{-1}\circ L_t$ identify a neighborhood of $0$ in $T_pM\simeq\rn{d}$ with neighborhoods of $p$ in $M$. Write $\widetilde{F}_t=\Theta_t^{-1}\circ F_t\circ\Theta_t$. Thus $D\widetilde{F}_t(0)$ preserves the fixed coordinate factors.

        By the $C^1$ continuity, in the $\Theta_t$ coordinates the stable and unstable plaques are graphs of maps
        \[
            \sigma_t:B^s_{r_0}\longrightarrow E^u,\qquad
            \varrho_t:B^u_{r_0}\longrightarrow E^s,
        \]
        where $B^s_{r_0}\subset E^s$ and $B^u_{r_0}\subset E^u$ are fixed small balls. More precisely, the stable plaque consists of points $\Theta_t(s,\sigma_t(s))$, and the unstable plaque consists of points $\Theta_t(\varrho_t(u),u)$. By the factorization, $\sigma_t(0)=\varrho_t(0)=0$, $D\sigma_t(0)=D\varrho_t(0)=0$. Define the local coordinate change
        \[
            \varPsi_t(s,u)=(s+\varrho_t(u),\,u+\sigma_t(s)).
        \]
        Since $D\varPsi_t(0)=\id$, shrinking the domains uniformly in $t$ makes these maps $C^1$ diffeomorphisms onto their images. Moreover, the $s,u$-axes are sent to the graph plaques respectively.

        In the resulting coordinates put $\widehat{F}_t:=\varPsi_t^{-1}\circ\widetilde{F}_t\circ\varPsi_t$. Define linear maps $S_t:E^s\to E^s$ and $U_t:E^u\to E^u$ by $D\widehat{F}_t(0)(s,u)=(S_ts,U_tu)$. If $R_t=(R_t^s,R_t^u)$ is the nonlinear remainder of $\widehat{F}_t$, then 
        \begin{equation}\label{eq:hat-F_t-iterates}
            \widehat{F}_t(s,u)=(S_ts+R_t^s(s,u),\,U_tu+R_t^u(s,u)),
        \end{equation}
        in which $R_t(0)=0$ and $DR_t(0)=0$. Invariance of the plaques implies invariance of the coordinate axes wherever these local expressions are defined. Hence
        \begin{equation}\label{eq:parameter-axis-remainders}
            R_t^s(0,u)=0,\qquad R_t^u(s,0)=0.
        \end{equation}

        Choose fixed adapted norms on $E^s,E^u$ for the hyperbolic map at $t=0$. After shrinking $I$, continuity of $S_t$ and $U_t^{-1}$ ensures that there is $0<\rho<1$ such that $\|S_t\|,\|U_t^{-1}\|\le\rho$ for $t\in I$. Equip the product with $\|(s,u)\|_{\max}:=\max\{\|s\|,\|u\|\}$ and write $Q_r=\{(s,u):\|s\|\le r,\ \|u\|\le r\}$. For sufficiently small $r$, the map $\varPsi_t$ and hence $\widehat{F}_t$ are defined. 

        Since $t\mapsto R_t$ is continuous in $C^1$ and $DR_t(0)=0$, reducing $r$ gives $\sup_{t\in I}\sup_{z\in Q_r}\|DR_t(z)\|\le\varepsilon$. The convexity of $Q_r$ and the mean value estimate imply, for $z,z'\in Q_r$,
        \[
            \|R_t(z)-R_t(z')\|_{\max}\le\varepsilon\|z-z'\|_{\max},
            \qquad \|R_t(z)\|_{\max}\le\varepsilon\|z\|_{\max}.
        \]
        Choose $0<\varepsilon<1-\rho$ and reduce $r$ accordingly. In particular $\rho\varepsilon<1$, and
        \[
            \tau:=\max\{\rho+\varepsilon,\rho(1+\varepsilon)\}<1,\qquad
            \gamma:=\max\left\{\rho+\varepsilon,\frac{\rho}{1-\rho\varepsilon}\right\}<1.
        \]

        \emph{Step 2. A boundary-value contraction.}
        Shrinking $\lssp{W}{\loc}{s}{p,F_t}$ and $\lssp{W}{\loc}{u}{p,F_t}$ if necessary, we can write the prescribed endpoints as $e_t=\Theta_t\varPsi_t(e_t^s,0)$, $d_t=\Theta_t\varPsi_t(0,d_t^u)$ where $\|e_t^s\|,\|d_t^u\|\le r/2$. For a fixed integer $n$, define the phase space 
        \[
            \mca{X}_n:=(Q_r)^{n+1}=\{Z=[z_j]\ssp{j=0}{n}=[(s_j,u_j)]\ssp{j=0}{n}:(s_j,u_j)\in Q_r\text{ for }0\le j\le n\},
        \]
        with a norm
        \[
            \|Z-\widetilde{Z}\|_n:=\max_{0\le j\le n}\max\{\|s_j-\tilde{s}_j\|,\|u_j-\tilde{u}_j\|\}.
        \]
        For an arbitrary input sequence $Z$, define the output coordinates of $\mca{T}_{n,t}Z$ by
        \begin{equation}\label{eq:mcaT-def}
        \begin{aligned}
            (\mca{T}_{n,t}Z)^s_0&=e_t^s,&\quad
            (\mca{T}_{n,t}Z)^s_{j}&=S_ts_{j-1}+R_t^s(s_{j-1},u_{j-1}),&(1\le j\le n),\\
            (\mca{T}_{n,t}Z)^u_n&=d_t^u,&
            (\mca{T}_{n,t}Z)^u_j&=U_t^{-1}u_{j+1}-U_t^{-1}R_t^u(s_j,u_j),&\quad(0\le j<n).
        \end{aligned}
        \end{equation}
        Thus the stable coordinate is prescribed at the entrance and the unstable coordinate at the exit, while the remaining coordinates, namely $(\mca{T}_{n,t}Z)_0^u$, $(\mca{T}_{n,t}Z)_n^s$ and $(\mca{T}_{n,t}Z)_j$, $1\le j<n$, are determined by the input $Z$.
        
        By \eqref{eq:hat-F_t-iterates}, if $\mca{T}_{n,t}Z=Z$, then $Z=[(s_j,u_j)]\ssp{j=0}{n}$ is an orbit segment of $\widehat{F}_t$, with entrance stable coordinate $e_t^s$ and exit unstable coordinate $d_t^u$. We therefore apply the Banach fixed-point theorem to find such a fixed point in $\mca{X}_n$. 

        \emph{Step 3. Application of the Banach fixed-point theorem.} 
        To check that $\mca{T}_{n,t}$ is a self-map, take any $Z\in\mca{X}_n$. Each input coordinate has norm at most $r$, so $R_t(0)=0$ implies $\|R_t^s(s_j,u_j)\|,\|R_t^u(s_j,u_j)\|\le\varepsilon r$. By \eqref{eq:mcaT-def}, for $0\le j<n$ we obtain
        \[
            \|(\mca{T}_{n,t}Z)^s_{j+1}\|
                \le\rho\|s_j\|+\varepsilon r\le(\rho+\varepsilon)r\le\tau r,\qquad
            \|(\mca{T}_{n,t}Z)^u_j\|
                \le\rho\|u_{j+1}\|+\rho\varepsilon r\le\rho(1+\varepsilon)r\le\tau r.
        \]
        The two boundary output coordinates have norm at most $r/2$. Since $\tau<1$, every output coordinate has norm at most $r$, proving
        $\mca{T}_{n,t}(\mca{X}_n)\subset\mca{X}_n$.

        Next take $Z,\widetilde{Z}\in\mca{X}_n$. The boundary outputs are independent of the input, so their differences vanish. For the remaining coordinates and $0\le j<n$, we have
        \[
            \begin{aligned}
                \|(\mca{T}_{n,t}Z)^s_{j+1}-(\mca{T}_{n,t}\widetilde{Z})^s_{j+1}\|
                    &\le\rho\|s_j-\tilde{s}_j\|+\varepsilon\|z_j-\tilde{z}_j\|_{\max}\le(\rho+\varepsilon)\|Z-\widetilde{Z}\|_n,\\
                \|(\mca{T}_{n,t}Z)^u_j-(\mca{T}_{n,t}\widetilde{Z})^u_j\|
                    &\le\rho\|u_{j+1}-\tilde{u}_{j+1}\|+\rho\varepsilon\|z_j-\tilde{z}_j\|_{\max}\le\rho(1+\varepsilon)\|Z-\widetilde{Z}\|_n.
            \end{aligned}
        \]
        Taking the maximum over all coordinates proves
        \begin{equation}\label{eq:mcaT-Lipschitz}
            \|\mca{T}_{n,t}Z-\mca{T}_{n,t}\widetilde{Z}\|_n\le\tau\|Z-\widetilde{Z}\|_n,
        \end{equation}
        In particular, $\|\mca{T}_{n,t}\|_{\mathrm{Lip}}\le\tau<1$, where $\|\cdot\|_{\mathrm{Lip}}$ denotes the Lipschitz constant. The constant $\tau$ and the self-map estimate are independent of $n$ and $t$.

        The Banach fixed-point theorem yields a unique fixed sequence $Z_n(t)=[(s_j(t),u_j(t))]_{j=0}^n$, whose boundary coordinates satisfy $s_0(t)=e_t^s$ and $u_n(t)=d_t^u$. Thus $z_{j+1}(t)=\widehat{F}_t(z_j(t))$. In particular, defining $x_n(t)=\Theta_t\varPsi_t(z_0(t))$ gives an orbit of $F_t$ with $F_t^j(x_n(t))=\Theta_t\varPsi_t(z_j(t))$. 

        Furthermore, $t\mapsto Z_n(t)$ is continuous. To see this, note
        \[
        \begin{aligned}
            \norm{Z_n(t)-Z_n(t_0)}{}&\le\norm{\mca{T}_{n,t}Z_n(t)-\mca{T}_{n,t}Z_n(t_0)}{}+\norm{\mca{T}_{n,t}Z_{n}(t_0)-\mca{T}_{n,t_0}Z_n(t_0)}{}\\
            &\le\tau\norm{Z_n(t)-Z_n(t_0)}{}+\norm{\mca{T}_{n,t}-\mca{T}_{n,t_0}}{}\cdot r.
        \end{aligned}
        \]
        Because $t\mapsto\mca{T}_{n,t}$ is continuous, we conclude that $\|Z_n(t)-Z_n(t_0)\|\le\frac{r}{1-\tau}\norm{\mca{T}_{n,t}-\mca{T}_{n,t_0}}{}\to0$ as $t\to t_0$.

        \emph{Step 4. Uniform endpoint and middle portion estimates.}
        By \eqref{eq:parameter-axis-remainders}, we have $\|R_t^s(s,u)\|\le\varepsilon\|s\|$ and $\|R_t^u(s,u)\|\le\varepsilon\|u\|$. The stable recurrence therefore implies
        \[
            \|s_{j+1}\|\le\|S_t\|\,\|s_j\|+\|R_t^s(s_j,u_j)\|\le(\rho+\varepsilon)\|s_j\|.
        \]
        The unstable recurrence in its backward form gives
        \[
            \|u_j\|\le\rho\|u_{j+1}\|+\rho\varepsilon\|u_j\|\quad\text{so that}\quad\|u_j\|\le\frac{\rho}{1-\rho\varepsilon}\|u_{j+1}\|.
        \]
        Consequently, with the definition of $\gamma$, induction on $j$ yields
        \begin{equation}\label{eq:parameter-passage-estimates}
            \|s_j(t)\|\le r\gamma^j,\qquad
            \|u_j(t)\|\le r\gamma^{n-j}\qquad(0\le j\le n),
        \end{equation}
        uniformly in $t$.

        Returning to the ambient manifold, recall that in chart $(\Theta_t\varPsi_t)\inv$, $e_t$ and the entrance $x_n(t)$ have the same stable coordinate, while $d_t$ and the exit $F^n(x_n(t))$ have the same unstable coordinate. Therefore, with uniform coordinate bounds on the compact parameter interval, \eqref{eq:parameter-passage-estimates} gives
        \[
            d(x_n(t),e_t)+d(F_t^n(x_n(t)),d_t)\le C\gamma^n,\qquad
            d(F_t^j(x_n(t)),p)\le C(\gamma^j+\gamma^{n-j}).
        \]
        With uniform continuity of each fixed finite iterate $F_t^i$ or $F_t^{n-i}$, these prove \eqref{eq:parameter-passage-ends} and \eqref{eq:parameter-passage-middle}.
        
        For a jointly continuous test function $\varphi:I\times M\to\R$, omit the first and last $j_*$ terms of its average. The remaining terms are uniformly close to its value at $p$, while the omitted contribution is at most $4j_*\|\varphi\|_\infty/n$. This proves \eqref{eq:parameter-passage-average} and completes the proof.
    \end{proof}

    Applying Lemma~\ref{lemma:uniform-parameter-passage} to the constant family $g$, with $(x,y)=(b',a)$ at $p$ and $(x,y)=(a',b)$ at $q$, gives points $x_n,y_m$ for all independently large $n,m$ such that
    \begin{equation}\label{eq:two-saddle-passages}
        x_n\to b',\qquad g^n(x_n)\to a,\qquad
        y_m\to a',\qquad g^m(y_m)\to b.
    \end{equation}
    These passages stay in $U_p$ and $U_q$, respectively, with occupation measures tending to $\delta_p$ and $\delta_q$. For sufficiently large $n,m$, Lemma~\ref{lemma:passage-transition-connection} joins their transition endpoints by two $C^1$-small contact corrections that preserve the intermediate iterates. Proposition~\ref{prop:almost flat periodic point} carries out this matching continuously along a parameter family.
\subsection{Local contact perturbations and the center multiplier}\label{subsection:Local contact perturbations and the center multiplier}
    We prepare two local perturbation tools. For the passage points $x_n,y_m$ of \eqref{eq:two-saddle-passages}, the next lemma chooses two disjoint boxes at the transition endpoints and uses the point-moving Lemma~\ref{lemma:point-moving} to match each final arrival with the next passage entrance while preserving the intermediate iterates. The second tool produces a family that changes the center exponent at $p$. All support neighborhoods are fixed before the long passage lengths are chosen. We combine these tools and select the parameter in Section~\ref{subsection:unit-center-multiplier}.

    \begin{lemma}[Fixed boxes, avoidance, and contact corrections]\label{lemma:passage-transition-connection}
        There are disjoint Darboux boxes $V_{a'},V_{b'}$, with closures avoiding $p,q$, such that $\overline{V_{a'}}\cap K_g=\{a'\}$ and $\overline{V_{b'}}\cap K_g=\{b'\}$. Let $x_n$ and $y_m$ be points constructed after the proof of Lemma~\ref{lemma:uniform-parameter-passage}. For all sufficiently large $n,m$, among the following list of points, counted in the displayed order,
        \begin{equation}\label{eq:uncorrected-arrivals}
            g^j(x_n)\qquad(1\le j\le n+T_a),\qquad g^j(y_m)\qquad (1\le j\le m+T_b),
        \end{equation}
        the only point in $V_{a'}$ is $g^{n+T_a}(x_n)$, and the only point in $V_{b'}$ is $g^{m+T_b}(y_m)$. Also $y_m\in V_{a'}$ and $x_n\in V_{b'}$.

        Moreover, for all sufficiently large $n,m$, there are smooth contactomorphisms $h^{a'}_{n,m}$ and $h^{b'}_{n,m}$ supported compactly in the respective boxes, such that
        \[
            h^{a'}_{n,m}(g^{n+T_a}(x_n))=y_m,\qquad
            h^{b'}_{n,m}(g^{m+T_b}(y_m))=x_n.
        \]
        Both corrections, and their composition $h^{b'}_{n,m}\circ h^{a'}_{n,m}$, tend to $\id$ in $C^1$ as $\min\{n,m\}\to\infty$.
    \end{lemma}
    \begin{proof}
        The only accumulation points of the two marked heteroclinic orbits are $p,q$. Thus $a',b'$ are isolated in $K_g$. First choose small disjoint boxes isolating these two points, and then take small closed saddle neighborhoods disjoint from the closures of the boxes. By \eqref{eq:parameter-passage-middle}, the middle portions of the passages lie in these saddle neighborhoods, uniformly for all sufficiently large passage times. Each of the finitely many remaining entrance, exit, and transition points converges to its corresponding marked orbit point. In the list \eqref{eq:uncorrected-arrivals}, the only limiting occurrences of $a',b'$ are the two final arrivals. Shrinking the boxes slightly, if necessary, proves the avoidance assertion. The endpoint inclusions follow from \eqref{eq:two-saddle-passages}.

        To construct the contact corrections, we join two varying endpoints within a fixed box. If $r_k,s_k\to c$ in a fixed Darboux box $V$, apply Lemma~\ref{lemma:point-moving} twice to obtain $R_k(c)=r_k$ and $S_k(c)=s_k$, with both maps supported compactly in $V$ and tending to $\id$ in $C^1$. Then
        \[
            h_k:=S_k\circ R_k^{-1},\qquad h_k(r_k)=s_k,
        \]
        is supported compactly in $V$ and tends to $\id$ in $C^1$, by continuity of inversion and composition. The cutoffs can be fixed on a compact subset of $V$.

        Apply this observation in $V_{a'}$ to $g^{n+T_a}(x_n),y_m\to a'$ and in $V_{b'}$ to $g^{m+T_b}(y_m),x_n\to b'$. These limits follow from \eqref{eq:two-saddle-passages} and continuity of the fixed transitions. The construction therefore provides the required corrections indexed by $(n,m)$. Their composition also tends to $\id$ in $C^1$. Since the boxes are disjoint and each correction maps its own box to itself, the composition realizes both prescribed endpoint moves without interference.
    \end{proof}

    The closing corrections will later be combined with a parameter family near $p$. We construct that family in a fixed neighborhood of the saddle, before selecting the long periodic orbit.

    \begin{lemma}[A local family changing the center exponent at $p$]\label{lemma:center-exponent-family}
        There are $\delta_0>0$, a fixed Darboux neighborhood $B_p$ of $p$, and a family $F_t=h_t\circ g$, $|t|\le\delta_0$, continuous in the $C^1$ topology, such that $h_0=\id$ and $\supp h_t\subset\subset B_p$. Consequently,
        \[
            F_t=g\quad\text{on }M\setminus g^{-1}(B_p).
        \]
        Both $p$ and $q$ remain fixed hyperbolic saddles with the same indices, and
        \begin{equation}\label{eq:saddle-exponent-parameter}
            \zeta_c(p;F_t)=\zeta_c(p;g)+t,\qquad \zeta_c(q;F_t)=\zeta_c(q;g).
        \end{equation}
        The set $B_p\cup g^{-1}(B_p)$ may be chosen inside any prescribed neighborhood of $p$. Moreover,
        \[
            \sup_{|t|\le\delta}\|F_t-g\|_{C^1}\longrightarrow0 \quad\text{as }\delta\to0.
        \]
    \end{lemma}
    \begin{proof}
        Choose $B_p$ so small that both it and its inverse image under $g$ lie in the prescribed neighborhood of $p$. In particular, they can avoid the closing boxes, a closed neighborhood of $q$, and all the marked points of the finite transitions. In Darboux coordinates centered at $p$, with $\alpha_{\rm std}=\D{z}-x\cdot\D{y}$, choose a fixed smooth cutoff $\eta$ supported in $B_p$ and equal to one near the origin. Let $h_t$ be the time-one map of the contact Hamiltonian
        \[
            H_t(x,y,z)=t\eta(x,y,z)\Bigl(z-\frac12x\cdot y\Bigr).
        \]
        On a sufficiently small neighborhood of the origin, uniformly for small $|t|$, its flow gives
        \[
            h_t(x,y,z)=\bigl(e^{t/2}x,e^{t/2}y,e^t z\bigr),\qquad h_t^*\alpha_{\rm std}=e^t\alpha_{\rm std}.
        \]
        In particular, $h_t(p)=p$. When we change from Darboux coordinates back to the manifold, the conformal factors of the coordinate change and its inverse cancel at this fixed point, so $\kappa_{h_t}(p)=e^t$ also holds relative to the global contact form. Thus $\kappa_{F_t}(p)=e^t\kappa_g(p)$, and Lemma~\ref{lemma:center multiplier conformal factor} proves \eqref{eq:saddle-exponent-parameter}. The perturbed map agrees with $g$ near $q$. The fixed cutoff makes $h_t$ smooth in $t$ and $C^1$-close to $\id$ as $t\to0$. Hyperbolicity, the indices, and membership in $\mca{U}$ persist after reducing $\delta_0$.
    \end{proof}
    
    The remaining auxiliary fact identifies the center multiplier as the quantity to be controlled by the later parameter choice. When the center direction along a periodic orbit is transverse to the contact distribution, this multiplier equals the product of the conformal factors along the orbit.

    \begin{lemma}\label{lemma:center multiplier conformal factor}
        Let $G$ be a contactomorphism with $G^*\alpha=\kappa_G\alpha$, and let $w$ be a periodic point of primitive period $N$. Suppose $G$ admits along $\Orb(w)$ a partially hyperbolic splitting $TM=E^{ss}\oplus E^c\oplus E^{uu}$ with $\dim E^c=1$ and $E^c\pitchfork\xi$, where $\xi=\ker\alpha$. Then
        \begin{equation}\label{eq:center eigenvalue conformal factor}
            \lambda_c(w;G)=\kappa_{G^N}(w)=\prod_{j=0}^{N-1}\kappa_G(G^j(w)).
        \end{equation}
    \end{lemma}
    \begin{proof}
        If $v_c\neq0$ belongs to $E^c(w)$, then applying $\alpha$ to $DG^N(w)v_c=\lambda_c(w)v_c$ gives
        \[
            \kappa_{G^N}(w)\alpha(v_c)=\alpha(DG^N(w)v_c)=\lambda_c(w)\alpha(v_c).
        \]
        Since $E^c\pitchfork\xi$, $\alpha(v_c)\neq0$ so that $\lambda_c(w)=\kappa_{G^N}(w)$.
    \end{proof}

    Recall that Proposition~\ref{prop:persistent-contact-data} provides the partially hyperbolic splitting and uniform transversality of $E^c$ to $\xi$ for every map in $\mca{U}$. Thus the lemma applies to the periodic orbits constructed below, whose center multipliers are positive.
\subsection{Long periodic orbits with center multiplier one}\label{subsection:unit-center-multiplier}
    Fix a based map $g\in\mca{B}$. We combine the uniform saddle passages of Section~\ref{subsection:hyperbolic passages} with the local contact perturbations of Section~\ref{subsection:Local contact perturbations and the center multiplier}. In Proposition~\ref{prop:almost flat periodic point}, we choose the saddle family and the two closing boxes, then construct periodic orbits for every sufficiently large pair $(n,m)$. These orbits depend continuously on $t$ and satisfy uniform exponent and orbit-average estimates. In Proposition~\ref{prop:flat lambda_c(w)=1}, we then choose a balanced pair of passage lengths $(n,m)$ and keep it fixed while varying $t$ to make the center multiplier equal to one.

    \begin{proposition}[A continuous family of long periodic orbits]\label{prop:almost flat periodic point}
        Fix $g\in\mca{B}$, and write
        \[
            \zeta_p=\zeta_c(p;g)<0,\qquad \zeta_q=\zeta_c(q;g)>0.
        \]
        There are $\delta_0>0$ and a family $F_t$ as in Lemma~\ref{lemma:center-exponent-family} such that, for every sufficiently large $n,m$, there are contactomorphisms $G_{n,m,t}\in\mca{U}$ and points $w_{n,m}(t)$, continuous in $t\in[-\delta_0,\delta_0]$ in the $C^1$ and manifold topologies, respectively, with the following properties.
        \begin{enumerate}[(i)]
            \item The saddles $p,q$ remain fixed, and $w_{n,m}(t)$ has least period $N_{n,m}=n+m+T_a+T_b$. There are $\omega(\delta)\to0$ as $\delta\to0$ and $\varepsilon_{n,m}\to0$ as $\min\{n,m\}\to\infty$ such that
            \begin{equation}\label{eq:parameter-closing-size}
                \sup_{|t|\le\delta}\|G_{n,m,t}-g\|_{C^1}\le\omega(\delta)+\varepsilon_{n,m}
                \qquad(0<\delta\le\delta_0).
            \end{equation}
            Given any neighborhood of $K_g$, all these periodic orbits lie in it after first restricting $\delta$ and then taking $n,m$ sufficiently large.
            \item The normalized center multiplier
            \[
                E_{n,m}(t):=\zeta_c(w_{n,m}(t);G_{n,m,t})=\frac{1}{N_{n,m}}\log\lambda_c(w_{n,m}(t);G_{n,m,t})
            \]
            is continuous, and
            \begin{equation}\label{eq:uniform-periodic-exponent}
                \sup_{|t|\le\delta_0}
                \left|E_{n,m}(t)-\frac{n(\zeta_p+t)+m\zeta_q}{N_{n,m}}\right|
                \longrightarrow0
                \quad\text{as }\min\{n,m\}\to\infty.
            \end{equation}
        \end{enumerate}
        Moreover, let $\mu_{n,m,t}$ denote the equidistributed measure $\mu_{w_{n,m}(t),G_{n,m,t}}$. For every continuous $\psi:M\to\R$,
        \begin{equation}\label{eq:balanced-periodic-measures}
            \sup_{|t|\le\delta_0}
            \left|\int\psi\,\D\mu_{n,m,t}-\frac{n\psi(p)+m\psi(q)}{N_{n,m}}\right|
            \longrightarrow0
            \quad\text{as }\min\{n,m\}\to\infty.
        \end{equation}
    \end{proposition}
    \begin{proof}
        \emph{Step 1: choose the family and the passages.}
        Choose the small plaques for $g$ in Lemma~\ref{lemma:uniform-parameter-passage} and use the freedom in Remark~\ref{remark:transition-time-choice} to place the marked points $a,a',b,b'$ in smaller compact subplaques. Then fix the transition times and the boxes $V_{a'},V_{b'}$ from Lemma~\ref{lemma:passage-transition-connection}. Choose the support neighborhood in Lemma~\ref{lemma:center-exponent-family} so that $B_p\cup g^{-1}(B_p)$ avoids the closures of the closing boxes, a closed neighborhood containing the $q$ passage, and the finite transition itineraries. These compact sets avoid $p$, so the required choice is possible. The same avoidance holds on small neighborhoods of the finite transition points.

        After constructing $F_t$, use $C^1$ continuity at $t=0$ to restrict the parameter interval so that the chosen smaller plaques lie in common coordinate neighborhoods satisfying the estimates of Lemma~\ref{lemma:uniform-parameter-passage}. At this stage, no further change of the marked points or transition times is required. The stable and unstable plaques of $p$ vary continuously in $C^1$ with $F_t$. Choose continuous marked points $b'_t,a_t$ on them with $b'_0=b'$ and $a_0=a$, for example by fixing their plaque coordinates. Apply Lemma~\ref{lemma:uniform-parameter-passage} with $e_t=b'_t$ and $d_t=a_t$ to obtain $x_n(t)$.

        At $q$ use the points $y_m$ from Lemma~\ref{lemma:uniform-parameter-passage} for the constant family. Their passages and subsequent length-$T_b$ transitions agree for $F_t$ and $g$, for all sufficiently large $m$, since the passages and all the transition arrivals avoid $B_p$. Define
        \[
            c_t=F_t^{T_a}(a_t),\qquad
            r_n(t)=F_t^{n+T_a}(x_n(t)),\qquad
            s_m=F_t^{m+T_b}(y_m).
        \]
        The $q$ passage and its following finite transition avoid $B_p$, so $s_m=g^{m+T_b}(y_m)$. Compactness of the parameter interval and continuity of the fixed transition iterates give
        \begin{equation}\label{eq:parameter-transition-limits}
            r_n(t)\longrightarrow c_t,\qquad s_m\longrightarrow b',\qquad c_0=a',
        \end{equation}
        uniformly in $t$. In particular $c_t\to a'$ as $t\to0$. The starting points $a_t$ vary with $t$, although the finite transition maps coincide with the original maps near their original itineraries.

        \emph{Step 2: uniform avoidance and continuous corrections.}
        After reducing $\delta_0$, the only arrivals in the closing boxes $V_{a'},V_{b'}$ given by Lemma~\ref{lemma:passage-transition-connection} in the lists
        \[
            F_t^j(x_n(t))\quad(1\le j\le n+T_a),\qquad
            g^j(y_m)\quad(1\le j\le m+T_b)
        \]
        are $r_n(t)=F_t^{n+T_a}(x_n(t))\in V_{a'}$ and $s_m=g^{m+T_b}(y_m)\in V_{b'}$, respectively; also $y_m\in V_{a'}$ and $x_n(t)\in V_{b'}$. For uniformity in the unbounded passage lengths, first choose closed saddle neighborhoods disjoint from both closing boxes. Lemma~\ref{lemma:uniform-parameter-passage} places every middle iterate there, with a single entrance/exit bound $j_*$ for all parameters. The remaining entrance, exit and transition lists are finite. At $t=0$ their limiting marked points have positive margins from the closed boxes, except at the designated final arrivals. Shrink the parameter interval to preserve these margins, then increase $n,m$ using the uniform finite-iterate limits. This proves the stated avoidance uniformly in $t,n,m$.

        The construction in Lemma~\ref{lemma:point-moving} gives, for each $c\in\{a',b'\}$, a smooth family $P_{c;z}$ of smooth contactomorphisms, defined for $z$ near $c$, such that
        \[
            P_{c;z}(c)=z,\qquad P_{c;c}=\id,\qquad\supp P_{c;z}\subset L_c\subset\subset V_c.
        \]
        The fixed cutoff and the smooth dependence of the Hamiltonian coefficients on the target coordinates provide the claimed smooth family. Its support is contained in the fixed set $L_c$. Set
        \begin{equation}\label{eq:parameter-endpoint-corrections}
            H^{a'}_{n,m,t}=P_{a';y_m}\circ P_{a';r_n(t)}^{-1},\qquad
            H^{b'}_{n,m,t}=P_{b';x_n(t)}\circ P_{b';s_m}^{-1}.
        \end{equation}
        These maps send $r_n(t)$ to $y_m$ and $s_m$ to $x_n(t)$, respectively. They are continuous in $t$ in $C^1$, by continuity of composition and inversion, and have fixed compact supports in the two disjoint boxes. Put
        \begin{equation}\label{eq:Gnmt-and-wnmt}
            G_{n,m,t}=H^{b'}_{n,m,t}\circ H^{a'}_{n,m,t}\circ F_t,
            \qquad w_{n,m}(t)=x_n(t).
        \end{equation}
        As $n\to\infty$ and $m\to\infty$, we have $x_n(t)\to b_t'$, $y_m\to a'$, $r_n(t)\to c_t$, and $s_m\to b'$. Consequently, the limiting family of maps $G_{n,m,t}$ as $n,m\to\infty$ is
        \[
            G_{\infty,t}=P_{b';b'_t}\circ P_{a';c_t}^{-1}\circ F_t.
        \]
        Uniform endpoint convergence and uniform continuity of the point-moving families on compact parameter sets imply
        \begin{equation}\label{eq:error-varepsilon_{n,m}}
            \varepsilon_{n,m}:=\sup_{|t|\le\delta_0}\|G_{n,m,t}-G_{\infty,t}\|_{C^1}\longrightarrow0.
        \end{equation}
        Since $G_{\infty,0}=g$ and $t\mapsto G_{\infty,t}$ is $C^1$-continuous, setting
        $\omega(\delta)=\sup_{|t|\le\delta}\|G_{\infty,t}-g\|_{C^1}$ proves \eqref{eq:parameter-closing-size}. In particular, for fixed nonzero $t$ the corrections need not tend to the identity. Shrinking $\delta_0$ and increasing the common threshold for $n,m$ ensure $G_{n,m,t}\in\mca{U}$.

        \emph{Step 3: exact itinerary and least period.}\phantomsection\label{claim:periodic time control}
        The avoidance just proved and the postcomposition order give the exact identities
        \begin{equation}\label{eq:two-box-periodic-itinerary}
            \begin{aligned}
                G_{n,m,t}^j(x_n(t))&=F_t^j(x_n(t))
                    &&(0\le j<n+T_a),\\
                G_{n,m,t}^{n+T_a+j}(x_n(t))&=g^j(y_m)
                    &&(0\le j<m+T_b),\\
                G_{n,m,t}^{n+m+T_a+T_b}(x_n(t))&=x_n(t).
            \end{aligned}
        \end{equation}
        Among the $n+m+T_a+T_b$ points indexed from $0$ through $n+m+T_a+T_b-1$, the only point in $V_{b'}$ is $x_n(t)$ and the only point in $V_{a'}$ is $y_m$. If the least period were a proper divisor of $n+m+T_a+T_b$, the initial point would recur in this list. Thus its least period is exactly $N_{n,m}=n+m+T_a+T_b$. The corrections fix neighborhoods of $p,q$, which consequently remain fixed for $G_{n,m,t}$.

        To check orbit containment in a given neighborhood of $K_g$, choose small saddle neighborhoods inside it. The entrance and exit portions outside these saddle neighborhoods have uniformly bounded lengths, and their limits at $t=0$, together with the finite transition points, belong to $K_g$. First restrict the parameter interval so that these finitely many limits remain in the prescribed neighborhood, and then increase $n,m$. Formula \eqref{eq:two-box-periodic-itinerary} proves (i).

        \emph{Step 4: uniform averages and the conformal cocycle.}\phantomsection\label{claim:spending time control}
        Put $\varphi(t,z)=\log\kappa_{F_t}(z)$, which is jointly continuous on $I\times M$. Relative to $F_t$, the two closing corrections contribute additional conformal factors only at $r_n(t)$ and $s_m$. Lemma~\ref{lemma:center multiplier conformal factor} therefore yields the exact equality
        \begin{equation}\label{eq:center eigenvalue long orbit}
            \begin{aligned}
                \log\lambda_c(w_{n,m}(t);G_{n,m,t})
                 ={}&\underbrace{\sum_{j=0}^{n-1}
                        \varphi(t,F_t^j(x_n(t)))}_{\text{passage near }p}
                    +\underbrace{\sum_{j=n}^{n+T_a-1}
                        \varphi(t,F_t^j(x_n(t)))}_{\text{transition }a\to a'}\\
                    &+\underbrace{\sum_{j=0}^{m-1}
                        \varphi(t,g^j(y_m))}_{\text{passage near }q}
                    +\underbrace{\sum_{j=m}^{m+T_b-1}
                        \varphi(t,g^j(y_m))}_{\text{transition }b\to b'}\\
                    &+\underbrace{\log\kappa_{H^{a'}_{n,m,t}}(r_n(t))}
                        _{\text{point-moving correction in }V_{a'}}
                    +\underbrace{\log\kappa_{H^{b'}_{n,m,t}}(s_m)}
                        _{\text{point-moving correction in }V_{b'}}.
            \end{aligned}
        \end{equation}
        Every term depends continuously on $t$ for fixed $n,m$, proving continuity of $E_{n,m}$.

        For the passage near $p$, Lemma~\ref{lemma:uniform-parameter-passage} applied to the jointly continuous $\varphi$ shows that the first sum divided by $n$ equals $\varphi(t,p)+o_n(1)=\zeta_p+t+o_n(1)$, uniformly in $t$. For the passage near $q$, the third sum divided by $m$ equals $\zeta_q+o_m(1)$ uniformly in $t$, since $F_t=g$ on the $q$ passage neighborhood. The second and fourth sums correspond to the transitions and contain only $T_a$ and $T_b$ terms, respectively, so they are uniformly bounded. The two closing logarithms are uniformly bounded as well, because their endpoints range in fixed compact neighborhoods and the point-moving families and their inverses are continuous in $C^1$ with positive conformal factors. Consequently there are $o_n(1),o_m(1)\to0$ and a constant $C$, independent of $n,m,t$, such that
        \[
            \left|E_{n,m}(t)-\frac{n(\zeta_p+t)+m\zeta_q}{N_{n,m}}\right|
            \le\frac{no_n(1)+mo_m(1)+C}{N_{n,m}}.
        \]
        This proves (ii) for arbitrary sufficiently large $n,m$.

        To prove \eqref{eq:balanced-periodic-measures}, fix a continuous $\psi:M\to\R$. The same passage estimates give
        \[
            \eta_n(\psi):=\sup_{|t|\le\delta_0}\bigg|\frac{1}{n}\sum_{j=0}^{n-1}\psi(F_t^j(x_n(t)))-\psi(p)\bigg|\longrightarrow0,\quad
            \zeta_m(\psi):=\bigg|\frac{1}{m}\sum_{j=0}^{m-1}\psi(g^j(y_m))-\psi(q)\bigg|\longrightarrow0.
        \]
        By the exact itinerary \eqref{eq:two-box-periodic-itinerary}, the remaining $T_a+T_b$ orbit points belong to the two transitions. Each contributes at most $\|\psi\|_\infty$ in absolute value. Therefore
        \[
            \sup_{|t|\le\delta_0}
                \left|\int\psi\,d\mu_{n,m,t}-\frac{n\psi(p)+m\psi(q)}{N_{n,m}}\right|\\
            \le\frac{n\eta_n(\psi)+m\zeta_m(\psi)+(T_a+T_b)\|\psi\|_\infty}{N_{n,m}}
                \longrightarrow0.
        \]
        Thus \eqref{eq:balanced-periodic-measures} also holds true.
    \end{proof}

    We now use the uniform estimates to choose a balanced sequence of passage lengths. For each sufficiently long pair, we keep the period fixed as the parameter varies and use the intermediate value theorem to obtain center multiplier one.

    \begin{proposition}[Long periodic orbits with center multiplier one]\label{prop:flat lambda_c(w)=1}
        For every $g\in\mca{B}$, every $C^1$ neighborhood $\mca{V}$ of $g$, every neighborhood $\mca{O}$ of $K_g$, and every $N_*>0$, there are a contactomorphism $G\in\mca{V}\cap\mca{U}$ and a periodic point $w$ of least period $N>N_*$ such that
        \[
            \Orb(w,G)\subset\mca{O},\qquad \lambda_c(w;G)=1.
        \]
        The points $p,q$ remain fixed points for $G$. In particular, one can choose a sequence $G_k\to g$ in $C^1$ with least periods $N_k\to\infty$, multipliers $\lambda_c(w_k;G_k)=1$, and
        \[
            \mu_{w_k,G_k}\stackrel{*}{\longrightarrow}\theta_g\delta_p+(1-\theta_g)\delta_q.
        \]
        Here $\theta_g=\zeta_q/(\zeta_q-\zeta_p)\in(0,1)$.
    \end{proposition}
    \begin{proof}
        Use the family of Proposition~\ref{prop:almost flat periodic point}. Choose $n_k=\lfloor\theta_g k\rfloor$, $m_k=k-n_k$ for large $k$. Then $N_k=k+T_a+T_b\to\infty$, and the uniform exponent estimate gives
        \[
            u_k:=\sup_{|t|\le\delta_0}|E_{n_k,m_k}(t)-\theta_g t|\longrightarrow0.
        \]
        For large $k$, the endpoint exponents at $-\delta_0$ and $\delta_0$ have opposite signs. Choose a zero $t_k$ by the intermediate value theorem. Crucially, $\theta_g|t_k|=|\theta_g t_k-E_{n_k,m_k}(t_k)|\le u_k$, so $t_k\to0$. Set $G_k=G_{n_k,m_k,t_k}$ and $w_k=w_{n_k,m_k}(t_k)$. The closing estimate, with errors \eqref{eq:error-varepsilon_{n,m}} and $\omega(0)=0$, gives
        \[
            \|G_k-g\|_{C^1}\le\omega(|t_k|)+\varepsilon_{n_k,m_k}\longrightarrow0,
        \]
        where the case $t_k=0$ follows directly from $G_{\infty,0}=g$. The multiplier is exactly one, the least period is $N_k$, and the perturbation preserves the saddles. The orbit-containment conclusion of Proposition~\ref{prop:almost flat periodic point}, together with $t_k\to0$, places the orbits in any prescribed neighborhood of $K_g$ for all large $k$. Its uniform averaging estimate gives
        \[
            \int\psi\,\D\mu_{w_k,G_k}
            =\frac{n_k\psi(p)+m_k\psi(q)}{N_k}+o(1)
            \longrightarrow\theta_g\psi(p)+(1-\theta_g)\psi(q).
        \]
        Taking $k$ sufficiently large proves the individual existence assertion as well.
    \end{proof}
\subsection{Flattening the center germ}\label{subsection:flatten-center-germ}
    By the construction in Section~\ref{section:based-maps}, every based map in $\mca{B}$ is smooth, as are all contact perturbations used there and in Sections~\ref{subsection:Local contact perturbations and the center multiplier}--\ref{subsection:unit-center-multiplier}. Hence each map $G_k$ provided by Proposition~\ref{prop:flat lambda_c(w)=1} is smooth.

    Fix one such map and orbit, and write
    \[
        G=G_k,\qquad w=w_k,\qquad N=N_k,\qquad R=G_k^{N_k}.
    \]
    Thus $R(w)=w$ and $DR(w)|_{E^c(w)}=\id$. We reserve $R$ for the return map. For any prescribed finite regularity, the smoothness of $R$ provides a local invariant center curve of sufficiently high finite regularity. We straighten this curve by contact coordinates and apply Lemma~\ref{lemma:contact lift} to the inverse of the return germ. The resulting map has the prescribed regularity and is arbitrarily $C^1$-close to $G$. Only finite regularity is needed for the subsequent orbit-splitting construction.

    \begin{proposition}[Flattening the center germ]\label{prop:flatten-center-germ}
        Let $G\in\Cont^\infty(M,\xi)$, and let $w$ be a periodic point of least period $N$. Suppose that the return map $R=G^N$ has the invariant partially hyperbolic splitting $T_wM=E^{ss}(w)\oplus E^c(w)\oplus E^{uu}(w)$, with one-dimensional center $E^c(w)$ and
        \[
            DR(w)|_{E^c(w)}=\id,\qquad E^c(w)\pitchfork\xi_w.
        \]
        For every integer $r\ge1$, every $C^1$-neighborhood $\mca{W}$ of $G$, and every neighborhood $V$ of $w$ whose closure avoids the other points of $\Orb(w,G)$, there exist $\widetilde{G}\in\mca{W}\cap\Cont^r(M,\xi)$ and a $C^r$ curve $C$ through $w$, tangent to $E^c(w)$, such that
        \[
            \widetilde{G}=G\text{ on }M\setminus V,\qquad
            \Orb(w,\widetilde{G})=\Orb(w,G),\qquad
            D\widetilde{G}(x)=DG(x)\quad(x\in\Orb(w,G)),
        \]
        and the perturbed return $\widetilde{R}=\widetilde{G}^N$ satisfies $\widetilde{R}|_C=\id$ as germs at $w$.
    \end{proposition}
    \begin{proof}
        Fix the integer $r\ge1$ and choose a smaller neighborhood $V_0$ of $w$, with $\overline{V_0}\subset V$, such that the sets $G^j(\overline{V_0})$, $0\le j<N$, are pairwise disjoint. We first construct an invariant center curve $C$ and a contact chart adapted to $C$. Since the return map $R=G^N$ is smooth, the finite-regularity center-manifold theorem \cite{HirschPughShub1977} gives a local $C^{r+3}$ invariant center curve $C$ satisfying
        \[
        w\in C,\qquad T_{w}C=E^c(w),\qquad R(C)\subset C
        \]
        as germs at \(w\). Note that the curve is not claimed to be unique. Since \(E^c(w)\pitchfork\xi_{w}\), after shrinking \(C\) we have $TC\pitchfork\xi$ at every point of \(C\). Choose a $C^{r+3}$ parametrization $\varGamma:(-\varepsilon,\varepsilon)\to C$ with $\varGamma(0)=w$, $\alpha(\dot\varGamma(t))=1$ for every sufficiently small \(t\).
        
        The invariance of \(C\) under $R$ defines a one-dimensional germ \(R_c\) by
        \begin{equation}\label{eq:R.Gamma=Gamma.Rc}
            R(\varGamma(t))=\varGamma(R_c(t)).
        \end{equation}
        In particular, we have $R_c(0)=0$. Differentiating this identity at \(t=0\) gives $DR(w)\dot\varGamma(0) = R_c'(0)\dot\varGamma(0)$. Because \(\dot\varGamma(0)\in E^c(w)\) and \(DR(w)|_{E^c(w)}=\id\), it follows that $R_c'(0)=1$.
        
        There is an explicit contact chart adapted to this curve. In a smooth Darboux chart for $\alpha$, write $\varGamma(t)=(a(t),b(t),c(t))$. The normalization gives $c'(t)-a(t)\cdot b'(t)=1$. Define
        \begin{equation}\label{eq:explicit-contact-chart}
            \vartheta(x,y,t)=\frac{1+a'(t)\cdot y}{1+x\cdot b'(t)},\qquad
            \Phi(x,y,t)=
            \bigl(a(t)+\vartheta(x,y,t)x,\ b(t)+y,
                    \ c(t)+a(t)\cdot y\bigr).
        \end{equation}
        A direct calculation gives $\Phi^*(\D z-x\cdot\D y)=\vartheta(\D t-x\cdot\D y)$. The factor $\vartheta$ is positive near the axis, and the differential at $(0,0,0)$ is invertible because $c'-a\cdot b'=1$. After composing with the ambient Darboux chart and shrinking, this gives the required contact chart $\varPhi$ into $V_0$, with $\varPhi(0,0,t)=\varGamma(t)$. Formula \eqref{eq:explicit-contact-chart} loses at most one derivative, so it is at least $C^{r+1}$. The return germ $R_c$ and its inverse have the regularity required by Lemma~\ref{lemma:contact lift}.

        Put $\widehat{R}=\varPhi\inv\circ R\circ\varPhi$. The invariance of \(C\) gives $\widehat{R}(0,0,t)=(0,0,R_c(t))$. The identity $\ker\varPhi^*\alpha=\xi_{\rm std}$ implies the inclusion $\varPhi\circ\Cont_c^r(U_{\varPhi},\xi_{\rm std})\circ\varPhi\inv\subset\Cont_c^r(V_0,\xi)$.
        
        Applying Lemma~\ref{lemma:contact lift} to $\psi:=R_c\inv\in C^{r+1}$ produces a $\delta>0$ and a $C^1$-small contact perturbation $\mca{L}_{\psi,\delta}\in\Cont_c^r(U_\varPhi,\xi_{\rm std})$ satisfying 
        $$
        \mca{L}_{\psi,\delta}(0,0,t)=(0,0,R_c\inv(t))
        $$
        for $\abs{t}\le\delta$. This lift fixes the origin with derivative $\id$. Choose $\delta$ small enough that the support of the lift is compactly contained in $U_\varPhi$.

        \begin{figure}[H]
            \centering
            \includegraphics[scale=0.7]{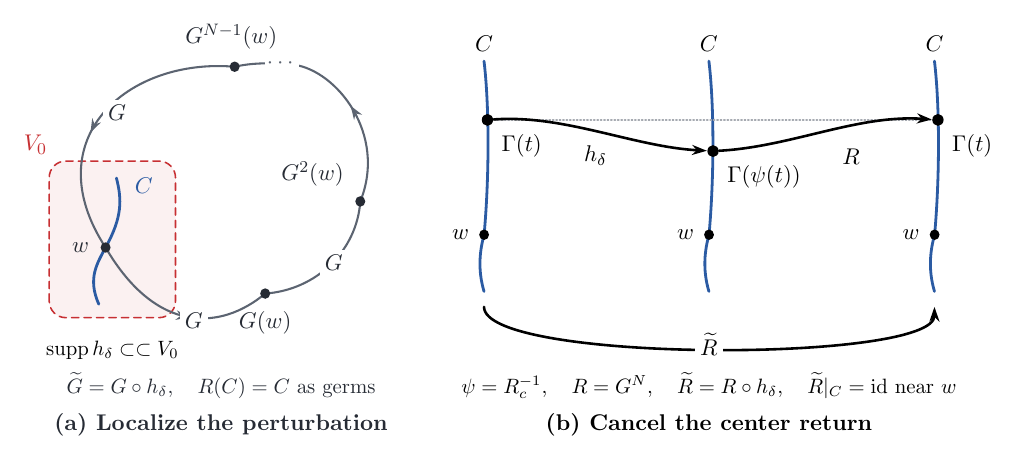}
            \caption[Flattening the center germ]{\textbf{Flattening the center germ.}
                (a) The correction $h_\delta$ is supported in $V_0$, which contains only $w$ from its period-$N$ orbit, and $\widetilde G=G\circ h_\delta$ preserves this orbit and its derivative. (b) Three copies of the local center curve $C$ illustrate the correction $\psi=R_c^{-1}$ followed by the return $R=G^N$, giving $\widetilde R|_C=\id$ as germs at $w$.}
            \label{fig:center-germ-flattening}
        \end{figure}
        
        Conjugate the localized perturbation back to \(M\) by defining $h_\delta:=\varPhi\circ\mca{L}_{\psi,\delta}\circ\varPhi^{-1}$ on $\varPhi(U_\varPhi)$, and extend this map by the identity outside the chart. Thus $h_\delta\in\Cont_c^r(V_0,\xi)$. The chart is fixed before choosing $\delta$, so $h_\delta$ is also arbitrarily $C^1$-close to the identity and 
        \begin{equation}\label{eq:h.Gamma=Gamma.psi}
            h_\delta(\varGamma(t))=\varGamma(\psi(t))
        \end{equation}
        for $t$ sufficiently small. The choice of $V_0$ ensures that $G^j(\supp h_\delta)\cap\supp h_\delta=\varnothing$ for every $1\le j\le N-1$. 
        
        Define $\widetilde{G}=G\circ h_\delta$. Then $\widetilde{G}=G$ outside $V$, and $\widetilde{G}$ is a $C^r$ contactomorphism arbitrarily $C^1$-close to $G$. Thus we can ensure $\widetilde{G}\in\mca{W}$. Since $h_\delta(w)=w$, $Dh_\delta(w)=\id$, and the support contains no other point of $\Orb(w,G)$, we have $\Orb(w,\widetilde{G})=\Orb(w,G)$, and $D\widetilde{G}(x)=DG(x)\,(x\in\Orb(w,G))$. In particular, the least period remains $N$. Moreover, for every point sufficiently close to \(w\), the disjointness gives
        \[
        \widetilde{R}:=\widetilde{G}^N=G^N\circ h_\delta=R\circ h_\delta.
        \]
        Combining this with \eqref{eq:h.Gamma=Gamma.psi} and \eqref{eq:R.Gamma=Gamma.Rc} yields, for all sufficiently small \(t\),
        \[
        \widetilde{R}(\varGamma(t))=R(h_\delta(\varGamma(t)))=R(\varGamma(\psi(t)))=\varGamma(R_c(\psi(t))).
        \]
        Since $\psi=R_c\inv$, we obtain $\widetilde{R}(\varGamma(t))=\varGamma(t)$ for all sufficiently small \(t\). Equivalently,
        \[
        \widetilde{R}|_C=\id
        \]
        as germs at \(w\). This completes the proof.
    \end{proof}

    We are now able to prove the main result of this section.

    \begin{proof}[Proof of Theorem~\ref{thm:contact flat}]
        Let $g\in\mca{B}$, fix an integer $r\ge1$, a $C^1$-neighborhood $\mca{V}$ of $g$, and $N_*>0$. Choose a smaller neighborhood $\mca{V}_0\subset\mca{V}\cap\mca{U}$. Proposition~\ref{prop:flat lambda_c(w)=1} gives a map $G_0\in\mca{V}_0$ with a periodic point $w$ of least period $N>N_*$ and center multiplier one. The map $G_0$ is smooth because $g$ and all corrections used in its construction are smooth.

        Apply Proposition~\ref{prop:flatten-center-germ} with this regularity $r$ to $G_0$ in a neighborhood of $w$ avoiding the remaining orbit points, choosing the perturbation small enough to remain in $\mca{V}\cap\mca{U}$. Denote the resulting $C^r$ map by $G$. The perturbation preserves the orbit and its derivative, so its least period is still $N$, and the proposition provides a $C^r$ center curve $C$ with $G^N|_C=\id$ as germs. This proves the theorem.
    \end{proof}
\section{Growth of periodic orbits}\label{section:superexponential}
    We continue to use the contact-open set $\mca{U}$, the dense set of smooth based maps $\mca{B}$, and the notation $O_n(f)$ from Section~\ref{section:based-maps}. Thus $O_n(f)$ counts distinct periodic orbits of least period $n$, and $\Cont^1(M,\xi)$ denotes coorientation-preserving contactomorphisms with the $C^1$ topology. We first perturb an identity center germ to create many hyperbolic orbits of one large least period. Hyperbolic continuation and the Baire category theorem then give arbitrarily fast growth along an unbounded sequence of periods. This follows the strategy of \cite[Sections~3--4]{BonattiDiazFisher2008}. The contact version of the construction in \cite{BonattiDiazFisher2008} yields the two adjacent unstable indices $d$ and $d+1$.

    \begin{proposition}[Creation of many periodic orbits]\label{prop:many-contact-periodic-orbits}
        Let $g\in\mca{B}$, let $r\ge1$ be an integer, and let $(a_n)_{n\ge1}$ be a sequence of positive integers. For every $C^1$-open neighborhood $\mca{V}\subset\mca{U}$ of $g$ and every integer $N_*\ge1$, there exist an integer $j\ge N_*$ and a contactomorphism $\widetilde{G}\in\mca{V}\cap C^r$ with hyperbolic periodic points $r_i^s,r_i^u$, $1\le i\le ja_j$, satisfying the following properties.
        \begin{enumerate}[(i)]
            \item Every $r_i^s$ and $r_i^u$ has least period $j$, and $\ind_u(r_i^s)=d$, $\ind_u(r_i^u)=d+1$;
            \item all $2ja_j$ periodic orbits represented by these points are distinct.
        \end{enumerate}
    \end{proposition}
    \begin{proof}
        Apply Theorem~\ref{thm:contact flat} at the finite order $r_0=r+2$. There are $G\in\mathcal V\cap\Cont^{r_0}(M,\xi)$ and a periodic point $w$ of least period $j>N_*$ such that $R=G^j$ fixes a $C^{r_0}$ transverse center curve $C$ pointwise as a germ at $w$. Choose a neighborhood $V$ of $w$ such that $G^\ell(\overline V)$, $0\le\ell<j$, are pairwise disjoint. Normalize the parameter on $C$ by $\alpha(\dot\varGamma)=1$.
        Formula~\eqref{eq:explicit-contact-chart} gives a contact chart
        of class $C^{r_0-1}=C^{r+1}$
        \[
            \varPhi:U_\varPhi\longrightarrow V,\qquad
            \varPhi(0,0,t)=\varGamma(t)\in C,\qquad
            \varPhi^*\alpha=\vartheta\alpha_{\rm std},\qquad
            \alpha_{\rm std}=\D t-x\cdot\D{y},
        \]
        with $\vartheta>0$. Write $\widehat{R}=\varPhi^{-1}\circ R\circ\varPhi$. Choose $\delta>0$ so that the compact arc $\varGamma([-\delta,\delta])$ lies in the interior of a larger arc on which $\widehat{R}(0,0,t)=(0,0,t)$. Shrinking $\delta$ if necessary, the transverse return derivatives along this arc are uniformly hyperbolic, with $d$ contracting and $d$ expanding directions. Its tangent is the one-dimensional center direction.

        Choose a smooth cutoff $\chi\in C_c^\infty(U_\varPhi;[0,1])$ equal to $1$ on a neighborhood of $\{(0,0,t):|t|\le\delta\}$. Fix an even integer $k>ja_j$ and set
        \[
            \widehat{H}(x,y,t)=\chi(x,y,t)v(t),\qquad v(t)=\sin(k\pi t/\delta).
        \]
        Let $\hat{h}^s$ be the contact Hamiltonian flow generated by $\widehat{H}$ with respect to $\alpha_{\rm std}$. Where $\chi=1$, equation~\eqref{eq:Hamiltonian eq. in xyz} gives $X_{\widehat{H}}=v'(t)x\cdot\partial_{x}+v(t)\partial_t$. Consequently, the points
        \[
            \begin{aligned}
                \hat{r}_\ell^u&=(0,0,2\ell\delta/k),&\qquad \ell&\in\Z\cap(-k/2,k/2),\\
                \hat{r}_\ell^s&=(0,0,(2\ell+1)\delta/k),&\qquad \ell&\in\Z\cap[-k/2,k/2),
            \end{aligned}
        \]
        are equilibria of this flow and fixed points of $\widehat{R}$. There are $k-1$ points of the first type and $k$ of the second. Since $k-1\ge ja_j$, choose $ja_j$ points of each type and relabel them $\hat{r}_i^u,\hat{r}_i^s$.

        At any of these points $\hat{r}=(0,0,t_i)$,
        \[
            D\hat{h}^s(\hat{r})=\diag\bigl(e^{s v'(t_i)}I_d,I_d,e^{s v'(t_i)}\bigr).
        \]
        The center tangent line is invariant under both derivatives, and $D\widehat{R}$ is the identity on that line. Thus for $s>0$, the center multipliers of $\hat{h}^s\circ\widehat{R}$ are
        \[
            \lambda_c(\hat{r}_i^u)=e^{sk\pi/\delta}>1,\qquad
            \lambda_c(\hat{r}_i^s)=e^{-sk\pi/\delta}<1.
        \]
        The induced derivatives on the quotient by the center tangent line remain hyperbolic for sufficiently small $s$, with $d$ contracting and $d$ expanding eigenvalues. Hence the selected fixed points are hyperbolic with unstable indices $d+1$ and $d$, respectively.

        Choose $k$ first. For this fixed smooth compactly supported Hamiltonian, its flow $\hat h^s$ tends to the identity in every finite $C^q$ topology as $s\to0$. Conjugation by the fixed $C^{r+1}$ chart therefore gives convergence in $C^r$. Put
        \[
            h^s=\varPhi\circ\hat{h}^s\circ\varPhi^{-1},
            \qquad \widetilde{G}=h^s\circ G,
        \]
        extending $h^s$ by the identity outside $V$. Conjugation by the contact chart preserves the contact distribution even though $\vartheta$ need not be constant. Thus $\widetilde G$ is of class at least $C^r$ and belongs to $\mathcal V$ for sufficiently small $s>0$. Figure~\ref{fig:many-periodic-saddles} shows how the alternating equilibria of $h^s|_C$ become distinct periodic saddles.

        \begin{figure}[H]
            \centering
            \includegraphics[scale=0.7]{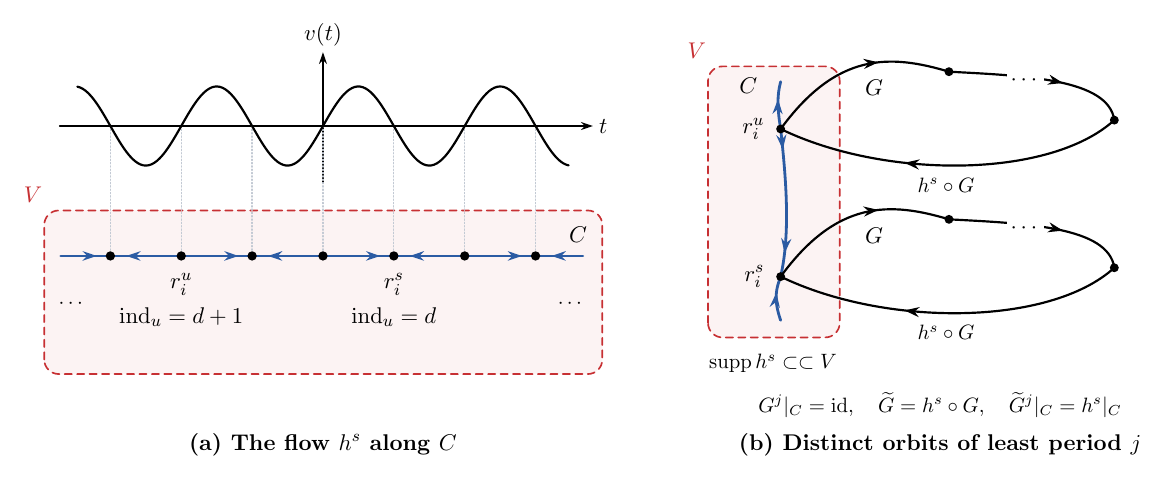}
            \caption[Creating many periodic saddles]{\textbf{Creating many periodic saddles.}
                (a) The flow $h^s|_C$ has alternating center-attracting and center-repelling equilibria, yielding saddles $r_i^s$ and $r_i^u$ of unstable indices $d$ and $d+1$ for sufficiently small $s>0$. (b) The return of $\widetilde G=h^s\circ G$ on the chosen center arc is $\widetilde G^j|_C=h^s|_C$. The chosen $2ja_j$ points have distinct orbits of least period $j$.}
            \label{fig:many-periodic-saddles}
        \end{figure}

        On $V$, disjointness of the intermediate images gives
        \[
            \widetilde{G}^\ell=G^\ell\quad(1\le\ell<j),
            \qquad \widetilde{G}^j=h^s\circ G^j.
        \]
        For $\bigcdot=s,u$, each $r_i^{\bigcdot}:=\varPhi(\hat{r}_i^{\bigcdot})$ is therefore a periodic point of $\widetilde{G}$. Its first $j-1$ iterates lie outside $V$, so its least period is $j$. Because each such orbit visits $V$ exactly once during a period, distinct selected points represent distinct periodic orbits. This completes the proof.
    \end{proof}

    Write $O_n^{\rm hyp}(f)$ for the number of distinct hyperbolic periodic orbits of least period $n$.

    \begin{proposition}[Baire property and open density]\label{prop:open-dense-periodic-growth}
        For a closed cooriented contact manifold $(M,\xi)$, the space $\Cont^1(M,\xi)$ with its $C^1$ topology is completely metrizable and hence Baire. In particular, the open set $\mca{U}$ is Baire. For every sequence $a=(a_n)_{n\ge1}$ of positive integers and every integer $N\ge1$, the set
        \[
            \mca{D}_N(a):=
            \bigcup_{n\ge N}
            \left\{f\in\mca{U}:O_n^{\rm hyp}(f)\ge na_n\right\}
        \]
        is open and dense in $\mca{U}$.
    \end{proposition}
    \begin{proof}
        For the topological assertion, take a complete metric $d$ for $C^1(M;M)$, obtained from a smooth Euclidean embedding of $M$. The compatible metric
        $d_{\rm Diff}(f,h)=d(f,h)+d(f^{-1},h^{-1})$ on $\Diff^1(M)$ is complete. Indeed, limits of a Cauchy sequence and its inverses satisfy both inverse identities by continuity of composition. The coorientation-preserving contactomorphisms form a closed subset for this metric. To see this, suppose $f_m\to f$ in $\Diff^1(M)$ and $f_m^*\alpha=\kappa_m\alpha$, $\kappa_m>0$. Evaluating on the Reeb vector field $R_\alpha$ shows that
        \[
            \kappa_m(x)=\alpha_{f_m(x)}(Df_m(x)R_\alpha(x))
                \longrightarrow
            \kappa(x)=\alpha_{f(x)}(Df(x)R_\alpha(x))
        \]
        uniformly. Thus $f^*\alpha=\kappa\alpha$ and $\kappa\ge0$. Since $Df$ is invertible and $\alpha$ is nowhere zero, $f^*\alpha$ is nowhere zero, which forces $\kappa>0$. This proves closedness and complete metrizability. The Baire category theorem applies to $\Cont^1(M,\xi)$ and to its open subset $\mca{U}$.

        Fix $n$ and a map $f$ with at least $na_n$ distinct hyperbolic periodic orbits of least period $n$. Choose these orbits. Hyperbolic continuation preserves this finite collection throughout a sufficiently small $C^1$-neighborhood of $f$. Pairwise disjoint neighborhoods of the finitely many orbit points ensure that the continued orbits have the same least periods and remain distinct. Therefore $\{f\in\mca{U}:O_n^{\rm hyp}(f)\ge na_n\}$ is open, and so is $\mca{D}_N(a)$.

        To prove density, let $\mca{V}\subset\mca{U}$ be nonempty and open. Proposition~\ref{prop:contact heterodim cycle} provides $g\in\mca{B}\cap\mca{V}$. Applying Proposition~\ref{prop:many-contact-periodic-orbits} in $\mca{V}$ with threshold $N$ gives $\widetilde{G}\in\mca{V}$ and $j\ge N$ with at least $ja_j$ distinct hyperbolic periodic orbits of least period $j$. Hence $\widetilde{G}\in\mca{V}\cap\mca{D}_N(a)$, proving density.
    \end{proof}

    \begin{theorem}[Arbitrarily fast growth of periodic orbits]\label{thm:arbitrarily-fast-contact-growth}
        Let $(M^{2d+1},\xi)$ be a closed contact manifold admitting an Anosov Reeb flow with $C^\infty$ stable and unstable foliations, as in \cite[Theorem~4.4]{Chaidez2025}. There exists a nonempty $C^1$-open set $\mca{U}\subset\Cont^1(M,\xi)$ such that, for every sequence $(a_n)_{n\ge1}$ of positive real numbers, there is a residual subset $\mca{R}_a\subset\mca{U}$ satisfying
        \[
            \limsup_{n\to\infty}\frac{O_n(f)}{a_n}=\infty
            \qquad\text{for every }f\in\mca{R}_a.
        \]
        Here $O_n(f)$ counts distinct periodic orbits of least period $n$. In particular, this conclusion holds on the unit tangent bundle of a closed hyperbolic surface with its canonical contact structure.
    \end{theorem}
    \begin{proof}
        Under the present hypothesis, the construction in \cite[Theorem~4.3 and the proof of Theorem~4.4]{Chaidez2025} provides the initial map $f_*$, the saddles $p_*,q_*$, the blender with its incoming strong unstable disk, and the transverse connection from $q_*$ to $p_*$ used in Proposition~\ref{prop:persistent-contact-data}. This configuration persists throughout the open set $\mca{U}$ in $\Cont^1(M,\xi)$. The remaining constructions of Section~\ref{section:based-maps} use closedness of $M$, the geometric blender and its attachments, and local contact perturbations to establish density of the smooth based maps $\mca{B}$ in $\mca{U}$. Section~\ref{section:Long flat orbit near the cycle} uses this persistent geometry and the local contact lemmas in dimension $2d+1$, so Theorem~\ref{thm:contact flat} and the preceding two propositions apply to this setting as well.

        Given $(a_n)$, define the positive integers $b_n=\max\{n,\lceil a_n\rceil\}$, which satisfy $b_n\ge a_n$. By Proposition~\ref{prop:open-dense-periodic-growth}, each $\mca{D}_N(b)$ is open and dense in the Baire space $\mca{U}$. Thus
        \[
            \mca{R}_a:=\bigcap_{N\ge1}\mca{D}_N(b)
        \]
        is a dense residual subset of $\mca{U}$. If $f\in\mca{R}_a$, then for every $N$ there is $n\ge N$ with $O_n^{\rm hyp}(f)\ge nb_n$. Choosing these periods to form a strictly increasing sequence $n_\ell\to\infty$, we obtain
        \[
            \frac{O_{n_\ell}(f)}{a_{n_\ell}}
                \ge\frac{O_{n_\ell}^{\rm hyp}(f)}{a_{n_\ell}}
                \ge n_\ell\frac{b_{n_\ell}}{a_{n_\ell}}
                \ge n_\ell\longrightarrow\infty.
        \]
        This proves the theorem.
    \end{proof}
\section{Growth of geometrically distinct translated points}\label{section:translated-point-proliferation}
    Fix a closed cooriented contact manifold $(M^{2d+1},\xi)$ and a smooth contact form $\alpha$ with $\ker\alpha=\xi$. We use the translated-point-shift convention of \eqref{eq:TP} and the counting functions $\TP_n(f;s_*)$ and $N_n(f;s_*)$ of \eqref{eq:intro-TP-count} and \eqref{eq:intro-N-count}. The positive conformal factor $\kappa_f$ and the Reeb flow $R_s^\alpha$ are also as in the Introduction. Distinct shifts at the same translated point contribute separately to the pair count, while each underlying point $x$ contributes only once to the geometric count. Our construction produces new pairs with shift $0$, although their shifts may become nonzero after a $C^1$ perturbation.

    \begin{theorem}[Translated-point proliferation from a contact blender]\label{thm:translated-point-proliferation}
        Let $(M,\xi)$ be a closed connected cooriented contact manifold admitting an Anosov Reeb flow whose strong stable and unstable foliations, including their tangent bundles, are $C^\infty$. Choose the contact form $\alpha$ and nonempty $C^1$-open set $\mca{U}\subset\Cont^1(M,\xi)$ of the blender construction in Proposition~\ref{prop:persistent-contact-data}. For each $s_*>0$ and every sequence $(a_n)_{n\ge1}$ of positive integers, there is a residual set $\mca{R}_a\subset\mca{U}$ for which
        \begin{equation}\label{eq:translated-point-arbitrarily-fast-growth}
            \limsup_{n\to\infty}\frac{N_n(f;s_*)}{a_n}
            =\limsup_{n\to\infty}\frac{\TP_n(f;s_*)}{a_n}
            =+\infty
            \qquad (f\in\mca{R}_a).
        \end{equation}
        In particular, there is a residual set $\mca{R}_\infty\subset\mca{U}$ whose elements satisfy
        \begin{equation}\label{eq:translated-point-superexponential-growth}
            \limsup_{n\to\infty}\frac{1}{n}\log\max\{1,N_n(f;s_*)\}
            =\limsup_{n\to\infty}\frac{1}{n}\log\max\{1,\TP_n(f;s_*)\}
            =+\infty.
        \end{equation}
    \end{theorem}

    Proposition~\ref{prop:contact heterodim cycle} gives a $C^1$-dense set of smooth based maps, while Theorem~\ref{thm:contact flat} produces arbitrarily long flat-center germs with any prescribed finite regularity. For the local construction below, fix an integer $r\ge2$ and apply Theorem~\ref{thm:contact flat} at order $r_0=r+2$. Thus, in any prescribed $C^1$ neighborhood of a based map, we obtain $G\in\Cont^{r+2}(M,\xi)$ and a $C^{r+2}$ curve $C$ fixed pointwise by $G^j$ near a point $w$ of arbitrarily large least period $j$, with $T_wC=E_G^c(w)\pitchfork\xi_w$. The contact coordinates adapted to $C$ lose at most one derivative, so the perturbations below have at least $C^r$ regularity and their translated-point equations are $C^1$. Nonzero local degree then yields the $C^1$-open conditions used in the Baire argument.

\subsection{Bounded shift multiplicity}
    The following elementary compactness observation converts the pair-count conclusion into a statement about geometrically distinct points. It applies independently of any blender or hyperbolicity assumption.

    \begin{lemma}\label{lemma:uniform-shift-multiplicity}
        Let $M$ be a closed contact manifold, let $\alpha$ be a smooth contact form, and fix $s_*>0$. There exists $\tau_\alpha>0$, depending only on $\alpha$, such that no Reeb trajectory has a nonzero return time of absolute value at most $\tau_\alpha$.  Consequently, with $K_{\alpha,s_*}:=1+\lfloor{2s_*}/{\tau_\alpha}\rfloor$, for every $f\in\Cont^1(M,\ker\alpha)$ and $n\ge1$ one has
        \begin{equation}\label{eq:uniform-shift-multiplicity}
            N_n(f;s_*)\le \TP_n(f;s_*)
            \le K_{\alpha,s_*}\,N_n(f;s_*).
        \end{equation}
        These inequalities also hold for infinite counts, interpreted as extended nonnegative integers. In particular, for fixed $s_*$, arbitrarily fast and super-exponential limsup growth for $\TP_n$ are equivalent to the respective growth statements for $N_n$.
    \end{lemma}
    \begin{proof}
        Since the Reeb field is everywhere nonzero, for each $x\in M$ the flowbox theorem gives a flowbox $U_x$ and a smaller neighborhood $V_x\subset\subset U_x$ of $x$ such that, for some $\tau_x>0$, the Reeb flow of any $y\in V_x$ stays inside the same flowbox for $|t|\le\tau_x$. In flowbox coordinates this flow is translation in one coordinate, so $R_t^\alpha(y)\neq y$ whenever $0<|t|\le\tau_x$. Choose a finite cover by the $V_x$ using compactness, and let $\tau_\alpha$ be the minimum of the corresponding positive $\tau_x$. Then $R_t^\alpha(y)\neq y$ for all $y\in M$ and $0<|t|\le\tau_\alpha$.

        Fix $f$, $n$, and a translated point $x$ counted by $N_n(f;s_*)$. If $s$ and $s'$ are two distinct admissible shifts for this point, then $R_{s-s'}^\alpha(x)=x$, which implies $|s-s'|>\tau_\alpha$. The interval $(-s_*,s_*)$ has length $2s_*$ and therefore contains at most $K_{\alpha,s_*}$ distinct admissible shifts for this point. Summing this uniform bound over the distinct underlying points gives \eqref{eq:uniform-shift-multiplicity}. The same finite-fiber bound shows that if either count is infinite, then both are infinite.

        For the last assertion, division by any positive sequence and comparison through the fixed factor $K_{\alpha,s_*}$ preserve divergence of the limsups to $+\infty$. Likewise,
        \[
            \max\{1,N_n(f;s_*)\}\le\max\{1,\TP_n(f;s_*)\}
            \le K_{\alpha,s_*}\max\{1,N_n(f;s_*)\}.
        \]
        Taking logarithms and dividing by $n$ changes the two growth rates by at most $(\log K_{\alpha,s_*})/n$ when the counts are finite; when they are infinite the corresponding logarithms are both $+\infty$. This proves the equivalence.
    \end{proof}

    \begin{proof}[Proof of Theorem~\ref{thm:intro-translated}]
        For the original geodesic Reeb flow on $S^*S$, both strong stable and strong unstable foliations are smooth, and the Anosov splitting satisfies the hypotheses of Theorem~\ref{thm:translated-point-proliferation}. Section~\ref{subsection:persistent-contact-configuration} constructs the model and its persistent contact-open set; see also Theorem~\ref{thm:arbitrarily-fast-contact-growth}. Applying Theorem~\ref{thm:translated-point-proliferation} with the contact form $\alpha$ chosen in Section~\ref{subsection:persistent-contact-configuration} yields both growth statements. This construction uses the same open set $\mca{U}$ as the periodic-orbit argument, with the residual subsets selected afresh.
    \end{proof}
\subsection{The flat-center contact perturbation}
    Fix $r\ge2$ and a flat center germ $F=G^j$ with $G$ and $C$ of class $C^{r+2}$, as above. Since $F$ fixes $C$ and $T_wC\pitchfork\xi_w$, its induced multiplier on $T_wM/\xi_w$ is $1$. Its strongly contracting and expanding generalized eigenspaces therefore lie in $\xi_w$ and together span it, so $DF_w|_{\xi_w}$ is hyperbolic. After shrinking $C$, continuity ensures that $TC\pitchfork\xi$ and $1\notin\operatorname{Spec}(DF_x|_{\xi_x})$ at every $x\in C$.

    Normalize a $C^{r+2}$ parametrization $\Gamma$ by $\alpha(\Gamma')=1$. Formula~\eqref{eq:explicit-contact-chart} gives a $C^{r+1}$ contact chart $\Phi$ about a compact subarc of $C$ such that
    \begin{equation}\label{eq:flat-center-darboux-chart}
        \Phi(0,0,t)=\Gamma(t)\in C,\qquad
        \Phi^*\alpha=\vartheta(x,y,t)\alpha_{\rm std},\qquad\vartheta>0,
    \end{equation}
    where $(x,y,t)\in\R^d\times\R^d\times\R$, and shrink the chart so that, for $|t|\le\delta$,
    \begin{equation}\label{eq:flat-return-axis-and-transverse-spectrum}
        \widehat{F}(0,0,t)=(0,0,t),\qquad
        1\notin\operatorname{Spec}\Bigl(D\widehat{F}_{(0,0,t)}\big|_{\xi_{\rm std}}\Bigr).
    \end{equation}
    The contact-plane restriction is well-defined because $\widehat{F}$ preserves the contact distribution.  Since $F$ fixes $C$ pointwise and $TC$ is transverse to $\xi$,
    \begin{equation}\label{eq:flat-return-conformal-factor}
        \kappa_F\big(\Gamma(t)\big)=1\qquad (|t|\le\delta).
    \end{equation}
    Choose $\chi\in C_c^\infty$ that is identically $1$ in a neighborhood of $\{(0,0,t):|t|\le\delta\}$.

    \begin{lemma}\label{lemma:regular-translated-point-creation}
        For any integer $m\ge1$ and any prescribed $C^1$ neighborhood of $G$, there is a contactomorphism $\widetilde{G}\in\Cont^r(M,\xi)$ in that neighborhood, obtained by a perturbation in the return chart, such that $\widetilde{G}^j$ has at least $m$ distinct translated-point pairs $(p_i,0)$. Each pair is an isolated regular zero of the full $(2d+2)$-dimensional translated-point equation.  These points have pairwise distinct least-period-$j$ orbits for $\widetilde{G}$.
    \end{lemma}
    \begin{proof}
        Fix a small $\eta>0$, to be chosen according to the desired $C^1$ tolerance, and choose an even integer $k>m+1$. Set
        \begin{equation}\label{eq:translated-point-oscillatory-hamiltonian}
            v_k(t)=\frac{\eta\delta^2}{\pi^2 k^2}
            \bigl(1-\cos(k\pi t/\delta)\bigr),\qquad
            H_k(x,y,t)=\chi(x,y,t)v_k(t).
        \end{equation}
        Let $h_k$ be the time-one map of the compactly supported contact Hamiltonian $H_k$ with respect to $\alpha_0=\D t-x\cdot \D y$.  Where $\chi=1$,
        \begin{equation}\label{eq:translated-point-contact-vector-field}
            X_{H_k}=v_k'(t)\sum_{r=1}^{d}x_r\partial_{x_r}+v_k(t)\partial_t.
        \end{equation}
        This follows directly from $\alpha_0(X_H)=H$ and $\iota_{X_H}\D\alpha_0=(\partial_tH)\alpha_0-\D H$. For every integer $\ell$ with $|\ell|<k/2$, put
        \begin{equation}\label{eq:translated-point-selected-axis-points}
            t_\ell=\frac{2\ell\delta}{k},\qquad
            z_\ell=(0,0,t_\ell),\qquad p_\ell=\Phi(z_\ell).
        \end{equation}
        There are $k-1\ge m$ such integers, and
        \begin{equation}\label{eq:oscillatory-profile-double-zeros}
            v_k(t_\ell)=v_k'(t_\ell)=0,\qquad
            v_k''(t_\ell)=\eta>0.
        \end{equation}
        By \eqref{eq:translated-point-contact-vector-field}, the local contact flow fixes $z_\ell$ at every time. Its time-one differential at $z_\ell$ is the identity because the linearization $DX_{H_k}(z_\ell)$ vanishes, as the potentially nonzero $v_k''$ term is multiplied by $x=0$. The flow also has conformal factor $1$ at $z_\ell$.

        The Hamiltonian $H_k$ and its local flow $h_k$ are smooth. Conjugating by the fixed $C^{r+1}$ contact chart $\Phi$ and extending by the identity gives a compactly supported $C^{r+1}$ contactomorphism $h_k^M$ of $M$. Choose the return neighborhood $V$ before localizing the perturbation so that $G^\ell(V)$, $0\le\ell<j$, are pairwise disjoint. Put
        \[
            \widetilde{G}=h_k^M\circ G.
        \]
        For points sufficiently close to the chosen arc, the disjointness implies
        \begin{equation}\label{eq:perturbed-flat-center-return}
            \widetilde{G}^j=h_k^M\circ G^j=h_k^M\circ F.
        \end{equation}
        Since $F$ fixes $p_\ell$ and $h_k^M$ fixes $p_\ell$, we have $\widetilde{G}^j(p_\ell)=p_\ell$. Equation~\eqref{eq:flat-return-conformal-factor}, the contact conformal chain rule, and $\kappa_{h_k^M}(p_\ell)=1$ then give $\kappa_{\widetilde{G}^j}(p_\ell)=1$. Thus $(p_\ell,0)$ is a translated-point pair.

        We next verify regularity in the \emph{full} translated-point equation, whose nondegeneracy condition differs from hyperbolicity of a fixed point. Choose a smooth coordinate chart $\varphi:U\to\mathbb{R}^{2d+1}$ around $p_\ell$. For $(z,s)$ sufficiently close to $(p_\ell,0)$, both $\widetilde{G}^j(z)$ and $R_s^\alpha(z)$ lie in $U$, so we may define
        \begin{equation}\label{eq:local-translated-point-map}
            \mca{T}_{\widetilde{G}^j}(z,s)=
            \bigl(\varphi(\widetilde{G}^j(z))
            -\varphi(R_s^\alpha(z)),\;
            \log\kappa_{\widetilde{G}^j}(z)\bigr).
        \end{equation}
        The map $\widetilde{G}$ has at least $C^r$ regularity, so its return has a $C^{r-1}$ conformal factor and $\mca{T}_{\widetilde{G}^j}$ is $C^{r-1}$, hence $C^1$. Its target and source both have dimension $2d+2$. In the coordinates $\varphi$, at $(p_\ell,0)$,
        \[
            D_z\bigl(\widetilde{G}^j-R_s^\alpha\bigr)\big|_{s=0}
            =D(\widetilde{G}^j)_{p_\ell}-I,
            \qquad \partial_s\bigl(\widetilde{G}^j-R_s^\alpha\bigr)
            =-R_\alpha(p_\ell).
        \]
        Use the decomposition $T_{p_\ell}M=\xi_{p_\ell}\oplus T_{p_\ell}C$. The derivative of $\widetilde{G}^j$ preserves both summands at $p_\ell$, with center derivative $1$ and contact-plane restriction $B_\ell$ satisfying $1\notin\operatorname{Spec}(B_\ell)$. The Reeb vector has a nonzero center component because $\alpha(R_\alpha)=1$ and $\alpha|_{TC}\ne0$. Finally, \eqref{eq:flat-return-conformal-factor} makes the derivative of $\log\kappa_F$ vanish along $TC$, while the perturbation gives
        \begin{equation}\label{eq:perturbed-conformal-center-derivative}
            \D(\log\kappa_{\widetilde{G}^j})_{p_\ell}
            \bigl(\Gamma'(t_\ell)\bigr)=v_k''(t_\ell)=\eta\ne0.
        \end{equation}
        To verify this identity, observe that on the center axis the time-one map of \eqref{eq:translated-point-contact-vector-field} is the flow $\psi_1$ of $\dot{t}=v_k(t)$. At an equilibrium satisfying $v_k=v_k'=0$, we have $\psi_1'(t_\ell)=1$ and $\psi_1''(t_\ell)=v_k''(t_\ell)$. For the conformal chart $\Phi^*\alpha=\vartheta\alpha_0$, the extra contribution to $\log\kappa_{h_k^M}$ on the axis is $\log\vartheta(0,0,\psi_1(t))-\log\vartheta(0,0,t)$. Its derivative vanishes at $t_\ell$ because $\psi_1(t_\ell)=t_\ell$ and $\psi_1'(t_\ell)=1$.

        Suppose $(u+av,\sigma)$ lies in the kernel of $D\mca{T}_{\widetilde{G}^j}(p_\ell,0)$, where $u\in\xi_{p_\ell}$ and $v\in T_{p_\ell}C$ is a chosen nonzero center vector. Apply $\alpha$ to the first component of the linearized equation. Because $D\widetilde{G}^j-I$ maps $\xi$ into $\xi$ and annihilates $v$, it follows that $\sigma=0$.  Then $(B_\ell-I)u=0$, so $u=0$. The last component and \eqref{eq:perturbed-conformal-center-derivative} imply $a=0$. Thus the Jacobian is invertible, proving regularity and local degree $\pm1$.

        Finally, the supports lie in a return box whose first $j-1$ forward images are disjoint. Each $p_\ell$ therefore has least period $j$. These points belong to distinct periodic orbits because a period-$j$ orbit visits the return box exactly once.

        Since the cutoff and $\delta$ were fixed before $k$, the perturbation satisfies
        \begin{equation}\label{eq:translated-point-perturbation-size}
            \|H_k\|_{C^2}\le C_{\chi,\delta}\eta,
            \quad \|X_{H_k}\|_{C^1}\le C'_{\chi,\delta}\eta,
            \quad \|h_k^M-\id\|_{C^1}\le C''\eta,
        \end{equation}
        with constants independent of $k$ for $\eta$ sufficiently small. Choosing $\eta$ small therefore places the perturbation in any prescribed $C^1$ neighborhood, \emph{uniformly in the number $k-1$ of translated pairs}.
    \end{proof}

    \begin{remark}
        Each selected $p_\ell$ satisfies $D(\widetilde{G}^j)_{p_\ell}v=v$ for the center tangent $v$, so the fixed point is necessarily nonhyperbolic. Indeed, if a hyperbolic fixed point $p$ satisfied $\kappa_{f^j}(p)=1$, then $\alpha_p\circ Df^j_p=\alpha_p$ would force $1$ to be an eigenvalue of $Df^j_p$. Regularity of the \emph{translated pair} follows from the extra conformal-factor equation, which is transverse to the center direction.
    \end{remark}
\subsection{Robustness through local degree}
    The translated-point map \eqref{eq:local-translated-point-map} involves the first derivative of the contactomorphism and therefore depends continuously on $f$ in the $C^1$ topology. Since its derivative at a regular zero involves second derivatives, we formulate robustness through local degree. The \emph{Brouwer (topological) degree} of a continuous map on a bounded domain, provided the target value is avoided on the boundary, is an integer invariant under homotopies that avoid that value on the boundary; see \cite[Chapter~1]{Deimling1985}. Its nonvanishing provides the $C^1$-open condition that derivative nondegeneracy alone does not provide.

    \begin{definition}\label{def:degree-essential-translated-pair}
        Let $f\in\Cont^1(M,\xi)$ and $(x_0,s_0)$ satisfy \eqref{eq:TP} for $f^n$. It is \emph{degree-essential in a neighborhood $W$} if $W$ is a relatively compact coordinate ball in $M\times(-s_*,s_*)$, if the local map $\mca{T}_{f^n}$ of \eqref{eq:local-translated-point-map} is nonzero on $\partial W$, and if
        \[
            \deg(\mca{T}_{f^n},W,0)\ne0.
        \]    
    \end{definition}

    The definition allows several roots in $W$ and applies to $C^1$ contactomorphisms because $\mca{T}_{f^n}$ is continuous.

    \begin{lemma}\label{lemma:degree-essential-pair-openness}
        Fix $n,m\ge1$ and $s_*>0$.  The set $\mca{E}_{n,m}(s_*)$ of contactomorphisms having $m$ disjoint neighborhoods as in Definition~\ref{def:degree-essential-translated-pair} is $C^1$-open.  Every contactomorphism provided by Lemma~\ref{lemma:regular-translated-point-creation} belongs to $\mca{E}_{j,m}(s_*)$.
    \end{lemma}
    \begin{proof}
        Choose $m$ pairwise disjoint closed isolating neighborhoods $\widebar{W}_i$ for the given degree-essential zeros.  For fixed $n$, the assignment
        \[
            f\mapsto \mca{T}_{f^n}|_{\widebar{W}_i}
        \]
        is continuous from the $C^1$ topology on contactomorphisms to the uniform topology on continuous maps, since compositions, first derivatives, the positive conformal factors, and the fixed Reeb flow vary continuously. Thus for every $f'$ sufficiently $C^1$ close to $f$, the maps $\mca{T}_{f'^n}$ and $\mca{T}_{f^n}$ are joined by a straight-line homotopy avoiding zero on every $\partial W_i$. Homotopy invariance of Brouwer degree implies
        \[
            \deg(\mca{T}_{f'^n},W_i,0)
            =\deg(\mca{T}_{f^n},W_i,0)\neq0.
        \]
        Hence each $W_i$ still contains a translated-point pair, yielding at least $m$ distinct pairs and proving openness. For the maps in Lemma~\ref{lemma:regular-translated-point-creation}, regularity of each constructed zero gives degree $\pm1$ in sufficiently small disjoint balls, which we choose inside $|s|<s_*$.
    \end{proof}
\subsection{Baire argument for exponential and arbitrarily fast growth}
    For $n,m\ge1$ denote by $\mca{E}_{n,m}(s_*)$ the $C^1$-open set of maps admitting at least $m$ pairwise disjoint isolating neighborhoods with nonzero local translated-point degree, as in Definition~\ref{def:degree-essential-translated-pair}. The open dense sets used in the proof are
    \begin{equation}\label{eq:translated-point-open-dense-sets}
        \mca{D}_N^{\mathrm{tr}}(a;s_*)
        :=\bigcup_{n\ge N}\mca{E}_{n,n a_n}(s_*).
    \end{equation}

    \begin{proof}[Proof of Theorem~\ref{thm:translated-point-proliferation}]
        For fixed $N$, equation \eqref{eq:translated-point-open-dense-sets} defines an open set by Lemma~\ref{lemma:degree-essential-pair-openness}. To prove density, let $\mca{V}\subset\mca{U}$ be nonempty and $C^1$ open. Proposition~\ref{prop:contact heterodim cycle} gives a smooth based map $g\in\mca{V}\cap\mca{B}$. Fix any integer $r\ge2$ and apply Theorem~\ref{thm:contact flat} at order $r+2$, in a sufficiently small neighborhood of $g$ contained in $\mca{V}$, to obtain $G\in\mca{V}\cap\Cont^{r+2}(M,\xi)$ with a $C^{r+2}$ flat-center curve and a point of least period $j\ge N$. Lemma~\ref{lemma:regular-translated-point-creation}, applied with $m=j a_j$, gives a $C^r$ contactomorphism $\widetilde{G}\in\mca{V}$ with the required regular translated pairs. Hence $\widetilde{G}\in\mca{E}_{j,ja_j}(s_*)\subset\mca{D}_N^{\mathrm{tr}}(a;s_*)$, proving density.

        By Proposition~\ref{prop:open-dense-periodic-growth}, the set $\mca U$ is Baire. Thus
        \[
            \mca{R}_a=\bigcap_{N\ge1}\mca{D}_N^{\mathrm{tr}}(a;s_*)
        \]
        is residual in $\mca{U}$.  For $f\in\mca{R}_a$, one can choose a strictly increasing sequence $n_i\to\infty$ with $\TP_{n_i}(f;s_*)\ge n_i a_{n_i}$.  By Lemma~\ref{lemma:uniform-shift-multiplicity}, with the constant $K=K_{\alpha,s_*}$ independent of $n_i$ and $f$,
        \[
            \frac{N_{n_i}(f;s_*)}{a_{n_i}}
            \ge \frac{\TP_{n_i}(f;s_*)}{K a_{n_i}}
            \ge \frac{n_i}{K}\longrightarrow+\infty.
        \]
        This proves both assertions of \eqref{eq:translated-point-arbitrarily-fast-growth}.

        For the last assertion, take $a_n^{(r)}=\lceil e^{rn}\rceil$ for each integer $r\ge1$ and let $\mca{R}_\infty=\bigcap_{r\ge1}\mca{R}_{a^{(r)}}$, which is residual by the Baire theorem. For each fixed $r$ and every $f\in\mca{R}_\infty$, an unbounded subsequence of $n$ satisfies
        \[
            \frac{1}{n}\log\max\{1,\TP_n(f;s_*)\}\ge r.
        \]
        Lemma~\ref{lemma:uniform-shift-multiplicity} further gives
        \[
            \frac{1}{n}\log\max\{1,N_n(f;s_*)\}
            \ge \frac{1}{n}\log\max\{1,\TP_n(f;s_*)\}
            -\frac{\log K_{\alpha,s_*}}{n}.
        \]
        Since $r$ is arbitrary and the correction tends to zero, both limsups in \eqref{eq:translated-point-superexponential-growth} equal $+\infty$.
    \end{proof}
\section{Contact graphs and fixed-window barcode entropy}\label{section:contact-graphs-and-barcode-entropy}
    We prove the entropy bound by combining a smooth compactification of the cocycle graph, a finite-dimensional family of contact perturbations, and stability of the completed action filtration. Analytic admissibility specifies the Floer-theoretic hypotheses under which we establish local finiteness of the associated barcodes.

\subsection{A smooth compactification of the contact graph}\label{section:cocycle-graph-compactification}
    Throughout this subsection, $Y$ is any closed smooth manifold. We write $\log^+x=\log\max\{1,x\}$ for $x\ge0$. Put $m=\dim Y$ and fix a Riemannian metric $q$ on $Y$. For $g\in C^\infty(Y,\R)$, define its \emph{additive cocycle} $G_n=\sum_{j=0}^{n-1}g\circ f^j$, so that $G_{n+k}=G_n+G_k\circ f^n$. Its \emph{cocycle graph} is the embedded $m$-manifold
    \begin{equation}\label{eq:cocycle-graph-definition}
        \Gamma_n=\{(x,f^n(x),G_n(x)):x\in Y\}\subset Y\times Y\times\R.
    \end{equation}
    We measure its $m$-volume using the product metric $q\oplus q\oplus \D\theta^2$. For a measurable subset of $\Gamma_n$, volume denotes $m$-dimensional Hausdorff measure. In the contact case, $g=\log\kappa_f$ and $(f^n)^*\alpha=e^{G_n}\alpha$.

    \begin{lemma}\label{lemma:circle-extension-entropy}
        Let $Y$ be a compact metric space, $f:Y\to Y$ a homeomorphism and $g:Y\to\R$ continuous. For the circle $\mathbb{T}=\R/\Z$ define
        \[
            f_{\mathbb{T}}:Y\times\mathbb{T}\to Y\times\mathbb{T},\qquad
            (x,\vartheta)\mapsto(f(x),\vartheta+g(x)\bmod\Z).
        \]
        Then $h_{\rm top}(f_{\mathbb{T}})=h_{\rm top}(f)$.
    \end{lemma}
    \begin{proof}
        Projection onto $Y$ makes $f_{\mathbb{T}}$ a $(\mathbb{T},\id_{\mathbb{T}})$-extension in Bowen's terminology: it intertwines the dynamics, and fiber translations commute with $f_{\mathbb{T}}$. Bowen's extension formula \cite[Theorem~19]{Bowen1971} therefore gives $h_{\rm top}(f_{\mathbb{T}})=h_{\rm top}(f)+h_{\rm top}(\id_{\mathbb{T}})=h_{\rm top}(f)$.
    \end{proof}

    \begin{theorem}[Cocycle-graph Yomdin estimate]\label{thm:cocycle-graph-volume-growth}
        Let $Y$ be a closed smooth manifold, $f\in\operatorname{Diff}^\infty(Y)$ and $g\in C^\infty(Y,\R)$. For the graphs in \eqref{eq:cocycle-graph-definition},
        \begin{equation}\label{eq:cocycle-graph-volume-entropy-bound}
            \limsup_{n\to\infty}\frac{1}{n}\log^+\operatorname{Vol}_m(\Gamma_n)
            \le h_{\rm top}(f).
        \end{equation}
        For every fixed $C>0$, the same bound holds with $\Gamma_n$ replaced by $\Gamma_n\cap\{|\theta|\le C\}$.
    \end{theorem}
    \begin{proof}
        On the compact manifold $X=Y_x\times Y_y\times\mathbb{T}_\vartheta$ define
        \[
            T(x,y,\vartheta)=(x,f(y),\vartheta+g(y)\bmod\Z),\qquad
            D=\{(x,x,0):x\in Y\}.
        \]
        Then $T=\id_Y\times f_{\mathbb{T}}$, so the product formula for entropy and Lemma~\ref{lemma:circle-extension-entropy} give $h_{\rm top}(T)=h_{\rm top}(f)$. Moreover
        \[
            T^nD=\{(x,f^n(x),[G_n(x)]):x\in Y\}.
        \]
        The parametrizations of $\Gamma_n$ and $T^nD$ induce the same pullback metric, namely,
        \[
            q_x+(Df_x^n)^*q_{f^nx}+\D G_n(x)\otimes \D G_n(x).
        \]
        Indeed, $\R\to\mathbb{T}, \theta\mapsto [\theta]$ is a local isometry, and the first factor $x$ ensures that both parametrizations are injective. Passing modulo one therefore preserves multiplicities, giving $\operatorname{Vol}(\Gamma_n)=\operatorname{Vol}(T^nD)$.

        Yomdin's smooth volume-growth inequality, applied to a fixed submanifold, states that for a smooth map $T$ of a compact manifold and a compact smooth submanifold $D$,
        \[
            \limsup_{n\to\infty}n^{-1}\log^+\operatorname{Vol}(T^nD)
            \le h_{\rm top}(T).
        \]
        This inequality follows from \cite{Yomdin1987}, while a graph-volume version is stated in \cite[Theorem~4.4.1]{Pelias2024}. This proves \eqref{eq:cocycle-graph-volume-entropy-bound}, and the same bound holds for the truncated subset because its measure is at most that of the whole graph.
    \end{proof}
\subsection{The Reeb diagonal and a contact tomograph}\label{section:reeb-diagonal-contact-tomograph}
    Return to a closed contact manifold $(Y,\alpha)$. Its \emph{contact product} is
    $Q=Y_x\times Y_y\times\R_\theta$ with $\beta=\alpha_y-e^\theta\alpha_x$; see \cite{Geiges2008}. Indeed,
    \[
        \ker\beta=\xi_x\oplus\xi_y\oplus
        \operatorname{span}\{\partial_\theta,R_\alpha(x)+e^\theta R_\alpha(y)\},
    \]
    and $\D\beta$ is nondegenerate on this space because the first two summands carry $-e^\theta\D\alpha_x\oplus\D\alpha_y$, while the last pair has pairing $-e^\theta$ and is orthogonal to them. The \emph{contact graph} of $\varphi$ is $\Gamma_\varphi=\{(x,\varphi(x),g_\varphi(x))\}$. In particular,
    \begin{equation}\label{eq:legendrian-contact-graph}
        \Gamma_n=\{(x,f^n(x),G_n(x)):x\in Y\}\subset Q
    \end{equation}
    is Legendrian because its dimension is $(\dim Q-1)/2=m$ and the pullback of $\beta$ is $(f^n)^*\alpha-e^{G_n}\alpha=0$.

    For a compact interval $I\subset\R$, consider the \emph{parameterized Reeb diagonal}
    \[
        \mathscr{D}_I:Y\times I\to Q,\qquad
        (x,s)\mapsto(x,R_s^\alpha (x),0).
    \]
    Closed Reeb trajectories may cause its image to self-intersect, so we always count \emph{parameterized} intersections. For each fixed $s$, the image is Legendrian. The definitions of the contact graph and the parameterized Reeb diagonal naturally identify their fiber product with
    \begin{equation}\label{eq:translated-pair-graph-incidence}
        \Gamma_n\cap_{\rm par}\mathscr{D}_I
        \longleftrightarrow
        \left\{
            (x,s): f^n(x)=R_s^\alpha(x),\ G_n(x)=0,\ s\in I
        \right\}.
    \end{equation}
    Thus these incidences are precisely the translated-point-shift pairs of Definition~\ref{def:tp}, where the equation $G_n(x)=0$ expresses preservation of $\alpha$ at $x$. Their expected dimension is zero, since $m+(m+1)=\dim Q$.

    \begin{lemma}\label{lemma:contact-hamiltonian-spanning-jets}
        There are smooth contact Hamiltonians $k_1,\ldots,k_N$ on $Y$ such that at every $y\in Y$ the vectors
        \[
            \bigl(\pi_\xi X^\alpha_{k_j}(y),R_\alpha k_j(y)\bigr),
            \qquad 1\le j\le N,
        \]
        span $\xi_y\oplus\R$. Here $X^\alpha_k$ is determined by $\alpha(X^\alpha_k)=k$ and $\iota_{X^\alpha_k}\D\alpha=(R_\alpha k)\alpha-\D k$, and $\pi_\xi(v)=v-\alpha(v)R_\alpha$.
    \end{lemma}
    \begin{proof}
        At a fixed $y$, a local smooth function can have arbitrarily prescribed first jet $k(y)=0$, $\D k_y|_{\xi_y}=\mu$ and $R_\alpha k(y)=c$. Since $\D\alpha|_{\xi_y}$ is nondegenerate, $(\pi_\xi X_k^\alpha(y),R_\alpha k(y))$ can therefore be prescribed arbitrarily. Cutoffs extend finitely many such functions to $Y$ without changing their jets at $y$. Spanning persists on a neighborhood, and a finite cover of the compact manifold $Y$ gives the desired global list.
    \end{proof}

    Set $k_z=\sum_{j=1}^Nz_jk_j$ and $\varphi_z=\phi_{k_z}^{1,\alpha}$, the time-one contact flow of $X_{k_z}^\alpha$. Smooth dependence of flows on parameters and compactness of $Y$ give $\varphi_0=\id$, and $\|\varphi_z-\id\|_{C^\ell}\le C_\ell|z|$ for each fixed $\ell$ and all sufficiently small $z$. We call this finite-dimensional probing family a \emph{contact tomograph}, following the integral-geometric terminology of \cite[Section~5]{CineliGinzburgGurel2024Entropy}.

    \begin{lemma}[Contact tomograph]\label{lemma:contact-tomograph-submersion}
        For a sufficiently small closed parameter ball $\Bpar\subset\mathbb{R}^N$ centered at zero, the map
        \begin{equation}\label{eq:contact-tomograph-map}
            \Psi:\Bpar\times Y\times I\to Q,
            \qquad\Psi(z,x,s)=\bigl(x,\varphi_z^{-1}(R_s^\alpha (x)),
            g_{\varphi_z^{-1}}(R_s^\alpha (x))\bigr)
        \end{equation}
        is a submersion, and the ball can be chosen independently of the compact shift interval $I$. For almost every $z$, simultaneously for all $n\ge1$ and all real shifts, $\Psi(z,\bigcdot,\bigcdot)$ intersects $\Gamma_n$ transversely. These parameterized intersections correspond to translated-point-shift pairs of $\varphi_z\circ f^n$.
    \end{lemma}
    \begin{proof}
        The contact lift
        \[
            \mathscr C_\varphi(x,y,\theta)=(x,\varphi(y),\theta+g_\varphi(y))
        \]
        satisfies $\mathscr C_\varphi^*\beta=e^{g_\varphi(y)}\beta$ and sends $\Gamma_{f^n}$ to $\Gamma_{\varphi\circ f^n}$. Since $g_{\varphi^{-1}}=-g_\varphi\circ\varphi^{-1}$, the map $\Psi_z$ is $\mathscr C_{\varphi_z}^{-1}\circ\mathscr D_I$. This proves the asserted correspondence of parameterized incidences.

        The first component of $D_x\Psi$ is the identity. After eliminating it, surjectivity is equivalent to spanning $T_yY\oplus\R$ by the following $z$ and $s$ columns, where $y=R_s^\alpha(x)$:
        \[
            \bigl(\partial_{z_j}\varphi_z^{-1}(y),\partial_{z_j}g_{\varphi_z^{-1}}(y)\bigr),\quad
            \bigl(D\varphi_z^{-1}(y)R_\alpha(y),\D g_{\varphi_z^{-1}}(y)R_\alpha(y)\bigr).
        \]
        At $z=0$ these are $(-X_{k_j}^\alpha(y),-R_\alpha k_j(y))$ and $(R_\alpha(y),0)$; they span by Lemma~\ref{lemma:contact-hamiltonian-spanning-jets}. The columns depend only on $(z,y)$, so compactness of $Y$ and openness of surjectivity give a parameter ball, contained in a slightly larger open ball, valid for every shift $s$.

        Parametric transversality \cite[Chapter~3]{Hirsch1976}, applied to the incidence manifold for each $n$ and each open interval $(-j,j)$, gives a null exceptional parameter set. The countable union over $n,j\ge1$ is still null, proving the simultaneous assertion.
    \end{proof}

    \begin{lemma}[Compact-family Crofton bound]\label{lemma:compact-family-crofton-bound}
        Let $\mca{Z}$ and $\mca{P}$ be compact smooth manifolds with boundary, viewed inside open extensions, and $Q$ a smooth Riemannian manifold. Suppose $\Psi:\mca{Z}\times \mca{P}\to Q$ extends to a smooth submersion on a neighborhood of its compact domain and $\dim \mca{P}+\dim V=\dim Q$, where $V\subset Q$ is an embedded smooth submanifold. Then there is a constant $C_\Psi$ independent of $V$ (depending on the fixed map $\Psi$, the metrics and its compact domain) such that, whenever the indicated intersections are discrete for almost every $z$,
        \begin{equation}\label{eq:compact-family-crofton-bound}
            \int_{\mca{Z}}\#\{p\in \mca{P}:\Psi(z,p)\in V\}\,\D z
            \le C_\Psi\operatorname{Vol}_{\dim V}
            (V\cap\Psi(\mca{Z}\times \mca{P})).
        \end{equation}
        One may allow an immersed $\mca{P}$ and count parameterized incidences.
    \end{lemma}
    \begin{proof}
        This is the integral-geometric argument of \cite[Lemma~5.3]{CineliGinzburgGurel2024Entropy}, with uniformity made explicit. We cover the compact domain by finitely many precompact submersion charts with inverse coordinate maps $\Xi(y,w)$ satisfying $\Psi(\Xi(y,w))=y$. In each chart the incidence manifold $\Psi^{-1}(V)$ is $(V\cap U)\times U_{\rm fib}$, of dimension $d=\dim\mca Z$. The incidence projection to $\mca Z$ is $p=\pi_{\mca Z}\circ\Xi$. Smoothness on the fixed precompact charts gives $\|Dp\|\le M$, independently of $V$, and therefore $J_dp\le M^d$ on every incidence tangent plane.

        Restrict to the measurable set $E$ of coordinate points lying over the original domain. The area formula \cite[Section~3.2.3]{Federer1969} gives
        \[
            \int_{\mca Z}\#(p^{-1}(z)\cap E)\,\D z
            =\int_E J_dp\,\D\operatorname{vol}_d
            \le C\operatorname{Vol}_{\dim V}(V\cap\Psi(\mca Z\times\mca P)).
        \]
        Here $C$ incorporates the bounded coordinate distortion and the volume of $U_{\rm fib}$, all fixed independently of $V$. Summing over the finite cover proves the bound; overlaps only overcount the left side. Working in open extensions and restricting to $E$ also handles boundaries. All counts are on the parameter domain, so the same argument applies to immersed families.
    \end{proof}

    \begin{theorem}\label{thm:translated-point-crofton-bound}
        For every fixed compact Reeb-shift interval $I$ there exist a ball $\Bpar\subset\R^N$, the family $\varphi_z$ above, and constants $C_I,c_I>0$, independent of $n$, such that
        \begin{equation}\label{eq:translated-point-crofton-bound}
            \int_{\Bpar}\#\TP(\varphi_z\circ f^n;I)\,\D z
            \le C_I\operatorname{Vol}_m(\Gamma_n\cap\{|\theta|\le c_I\})
            \le C_I\operatorname{Vol}_m(\Gamma_n).
        \end{equation}
        For almost every $z$, the integrand counts finitely many nondegenerate translated-point-shift pairs.
    \end{theorem}
    \begin{proof}
        Apply Lemma~\ref{lemma:compact-family-crofton-bound} with $\mca{Z}=\Bpar$, $\mca{P}=Y\times I$, $V=\Gamma_n$, and the submersion of Lemma~\ref{lemma:contact-tomograph-submersion}. Its image has uniformly bounded $\theta$ coordinate because $|g_{\varphi_z^{-1}}(y)|\le c_I$ on the compact parameter set. Differentiating the two defining equations in \eqref{eq:translated-pair-graph-incidence} shows that the translated-point equation is transverse exactly when the corresponding graph incidence is transverse. The transverse incidence set is zero-dimensional and compact on the compact parameter domain, hence finite. To include the endpoints of $I$, we apply transversality first on a slightly larger open shift interval. The isolated zeros in the resulting compact subinterval still form a finite set, and the volume bound follows.
    \end{proof}
\subsection{Analytic admissibility and the associated barcodes}\label{section:analytic-admissibility-and-barcodes}

    We next specify a \emph{normalized filtered package} consisting of Hamiltonian representatives, regular Floer choices, continuation maps (chain maps associated to homotopies of these choices), and action filtrations within perturbed Rabinowitz Floer homology (RFH) \cite[Theorems~2.9 and 2.16]{AlbersFrauenfelder2010Leafwise} and \cite[Proposition~2.8]{AlbersMerry2013}. Analytic admissibility concerns these choices, while the barcode construction uses persistence measures \cite[Section~2]{ChazalDeSilvaGlisseOudot2016}.

    Fix a nonnegative $\chi\in C_c^\infty((0,1/3))$ of integral one and the radial defining function $\Fsig$ of Lemma~\ref{lemma:radial-defining-hamiltonian}. For each $n$, Lemma~\ref{lemma:compactly-supported-hamiltonian-lift} gives a compactly supported Hamiltonian $K_{n,t}$, supported in $(2/3,1)$, generating the inverse homogeneous lift $\hat{f}^{-n}$ on the relevant trajectories. Its radial cutoff is chosen after $n$. Similarly choose a smooth family $L_{z,t}$, supported in $(1/3,2/3)$ and a fixed radial annulus, generating $\hat{\varphi}_z^{-1}$ along trajectories from $\Sigma$, with $L_{0,t}=0$. A fixed cutoff and bounded parameter derivatives give
    \begin{equation}\label{eq:tomograph-hamiltonian-size}
        \int_0^1\|L_{z,t}\|_\infty\,\D t\le C_L|z|,
    \end{equation}
    and, after shrinking the closed parameter ball,
    \begin{equation}\label{eq:tomograph-hamiltonian-parameter-lipschitz}
        \int_0^1\|L_{z,t}-L_{w,t}\|_\infty\,\D t\le C_L|z-w|.
    \end{equation}
    The constant is independent of $n$. Choose the annulus for $K_{n,\bigcdot}$ to contain all inverse-$f^n$ homogeneous isotopy trajectories from the compact set $\bigcup_{z\in\Bpar}\hat{\varphi}_z^{-1}(\Sigma)$. Thus $K_{n,t}+L_{z,t}$ generates $\hat{f}^{-n}\circ\hat{\varphi}_z^{-1}=\widehat{\varphi_z\circ f^n}^{-1}$ along every relevant critical trajectory.

    For fixed $n$ and $z\in\Bpar$, define on $C^\infty(S^1,\widehat{W})\times\mathbb{R}$ the functional
    \begin{equation}\label{eq:normalized-rabinowitz-functional}
        \mathscr{A}_{n,z}(u,\eta)=
        \int_{S^1}u^*\lambda-\int_0^1(K_{n,t}+L_{z,t})(u(t))\,\D t
        -\eta\int_0^1\chi(t)\Fsig(u(t))\,\D t.
    \end{equation}
    Here $\Fsig$ is the fixed real-valued defining Hamiltonian, while $\varphi_z$ is the contactomorphism used to construct $L_{z,t}$. We write $K_{n,\bigcdot}$ for the family $t\mapsto K_{n,t}$. A parameter is called \emph{regular} when $\mathscr{A}_{n,z}$ has nondegenerate critical points and the auxiliary Floer choices are regular.

    The \emph{action spectrum} is the set of critical values of $\mathscr{A}_{n,z}$. For regular endpoints outside this spectrum, $\CF_{n,z}^{(a,b)}$ denotes the \emph{finite-window Floer complex}, generated over $\F=\Z/2\Z$ by critical points of action in $(a,b)$. Its differential counts the regular index-one Floer trajectories, modulo the lower-action subcomplex \cite[Section~2.3, Eq.~(2.43)]{AlbersFrauenfelder2010Leafwise}.

    \begin{definition}[Analytically admissible family]\label{def:analytically-admissible-family}
        The family $\{\mathscr{A}_{n,z}\}_{n\ge1,\;z\in\Bpar}$, with its chosen Hamiltonian representatives and auxiliary Floer choices, is \emph{admissible} if the following hold for every fixed $n$.
        \begin{enumerate}[label=(\roman*),leftmargin=*]
            \item The regular parameters have full Lebesgue measure in $\Bpar$ (and hence are dense). Regular auxiliary choices can be made at the nondegenerate parameters used below. For every regular member, each compact action window has finitely many generators, and the finite-window Floer complexes over $\F=\Z/2\Z$ exist with strictly action-decreasing differential and compatible window inclusions and quotient maps.
            \item The Floer trajectories, continuation trajectories between regular members, and parametrized trajectories for the standard continuation homotopies have the compactness, transversality and gluing properties of perturbed RFH. In particular, the affine Hamiltonian comparisons and their contractions define chain maps and chain homotopies satisfying the usual action-energy estimates.
        \end{enumerate}
    \end{definition}
    
    The tomograph gives the critical-point transversality in (i) almost everywhere, whereas Floer transversality and the coherent moduli-space constructions remain analytic requirements. The analytic constants may depend on $n$. Barcode existence and local finiteness are consequences proved below.

    \begin{lemma}\label{lemma:normalized-critical-translated-pairs}
        For every $z$ sufficiently small, critical points of $\mathscr{A}_{n,z}$ correspond bijectively to translated point pairs of $\varphi_z\circ f^n$, and the critical action is $\eta$.
    \end{lemma}
    \begin{proof}
        Variation in $\eta$ enforces $\Fsig=0$ during the first time interval, when $u$ follows a reparametrized Reeb trajectory on $\Sigma$. The inverse homogeneous lift of $\varphi_z$ acts in the middle interval, followed by that of $f^n$ in the last interval. The time-one map of the unconstrained Hamiltonian part is therefore $\hat{f}^{-n}\circ\widehat{\varphi_z}^{-1} =\widehat{\varphi_z\circ f^n}^{-1}$ on the relevant trajectories. The loop-closing equation is equivalent to $(\varphi_z\circ f^n)(x)=R_\eta^\alpha x$, while the return to radius one is equivalent to $\kappa_{\varphi_z\circ f^n}(x)=1$. Conversely, each translated point pair defines the corresponding loop. Along homogeneous Hamiltonian trajectories, $\lambda(X_{K_{n,t}})=K_{n,t}$ and $\lambda(X_{L_{z,t}})=L_{z,t}$, whereas on $\Sigma$ the constraint $\Fsig$ is zero and $\lambda(R_\alpha)=1$. Substitution into \eqref{eq:normalized-rabinowitz-functional} gives action $\eta$. This calculation follows \cite[Lemma~2.2]{AlbersMerry2013}, with the separated supports arranged to match the order of composition.
    \end{proof}

    \medskip\noindent\textbf{The completed action filtration}

    Fix an admissible family. For a regular member, write $\mca{G}_{n,z}=\Crit\mathscr{A}_{n,z}$ and $A_{n,z}(x)=\mathscr{A}_{n,z}(x)$ for $x\in\mca{G}_{n,z}$. The parity grading suffices throughout, and every barcode count below sums multiplicities over all degrees. In each degree, we complete in the direction of decreasing action by setting
    \begin{equation}\label{eq:decreasing-action-chain-completion}
        \begin{aligned}
            C_{n,z}=\biggl\{&\sum_{x\in\mca{G}_{n,z}}a_xx:\;a_x\in\F,
            \#\{x:a_x\ne0,\ A_{n,z}(x)\ge r\}<\infty
            \text{ for every }r\in\mathbb{R}\biggr\}.
        \end{aligned}
    \end{equation}
    Let $C_{n,z}^{\le t}$ consist of sums supported in action at most $t$. The Floer differential extends to these sums because, above any fixed lower action level, its evaluation on a given chain involves only finitely many generators and a finite-window differential. Compatibility of the window differentials gives $\partial^2=0$. For a regular upper level $b$, $C_{n,z}^{\le b}$ is equivalently the chain-level inverse limit of $\CF_{n,z}^{(a,b)}$ as regular $a\to-\infty$. Taking homology \emph{after} this completion, we set
    \begin{equation}\label{eq:completed-action-persistence}
        V_{n,z}^t=H_*(C_{n,z}^{\le t}),\qquad t\in\mathbb{R},
    \end{equation}
    with inclusion-induced structure maps. Lemma~\ref{lemma:completion-homology-inverse-limit} justifies the required comparison with inverse limits of finite-window homology.

    A \emph{persistence module} is a family of vector spaces $V^t$ with maps $V^a\to V^b$ for $a\le b$, whose compositions agree with the corresponding structure maps and whose maps at equal parameters are identities \cite[Section~1.1]{ChazalDeSilvaGlisseOudot2016}. Its persistence diagram is the multiset associated to the module's rectangle-valued persistence measure \cite[Sections~2.1--2.7]{ChazalDeSilvaGlisseOudot2016}. The numerical endpoints indicate the parameter values at which classes appear and disappear. The \emph{bottleneck distance} is the infimum of the costs of partial matchings between diagrams, where matched endpoints move by at most the cost and an unmatched interval must have length at most twice the cost \cite[Section~4.2]{ChazalDeSilvaGlisseOudot2016}. We omit endpoint decorations.

    \begin{definition}[Associated barcode]\label{def:associated-finite-birth-barcode}
        For a regular member, $\mca{B}_{n,z}$ is the part of the \emph{undecorated persistence diagram} of $V_{n,z}$ with finite birth coordinate, represented as a multiset of bars $[p,q)$ with $p\in\mathbb{R}$ and $p<q\le+\infty$. Its existence follows from Lemma~\ref{lemma:regular-barcode-endpoint-counts} and the persistence-measure construction of~\cite[\S\S2.4--2.8]{ChazalDeSilvaGlisseOudot2016}. For a nonregular parameter, define $\mca{B}_{n,z}$ as the bottleneck limit of the barcodes at regular parameters approaching $z$. Proposition~\ref{prop:barcode-regularization-local-finiteness} proves that this limit exists and is independent of the regularizing sequence. Set $\mca{B}_n=\mca{B}_{n,0}$.
    \end{definition}

    The notation $[p,q)$ specifies numerical endpoints, with endpoint decorations omitted in the bottleneck metric. We include $q=+\infty$ and omit zero-length bars and births at $-\infty$, since the latter never contribute to the bounded birth windows used below. This convention does not assert that the entire completed module is a direct sum of interval modules. The construction depends on the specified completion, filling and $K_{n,\bigcdot}$, without asserting independence of unrelated Hamiltonian extensions.

    Hereafter, \emph{locally finite barcode} means that for every compact interval $I\subset\mathbb{R}$ and every $\varepsilon>0$,
    \begin{equation}\label{eq:bounded-birth-window-local-finiteness}
        \#\{[p,q)\in\mca{B}:p\in I,\ q-p\ge\varepsilon\}<\infty,
    \end{equation}
    counting multiplicity and allowing $q=+\infty$. In persistence-measure terminology, \emph{$v$-tameness} and \emph{$h$-tameness} mean finiteness on every vertically and horizontally infinite strip, respectively, separated from the diagonal \cite[Section~2.8]{ChazalDeSilvaGlisseOudot2016}. The displayed bounded-birth condition is the barcode form of $v$-tameness, which differs from local finiteness of a \emph{module} in their Proposition~2.21.

    \begin{lemma}\label{lemma:regular-barcode-endpoint-counts}
        For every regular member of an admissible family, $V_{n,z}$ is $v$-tame and $h$-tame. Its finite-birth barcode satisfies, for every bounded interval $I$,
        \begin{equation}\label{eq:regular-barcode-endpoint-census}
            \sum_{[b,d)\in\mca{B}_{n,z}}
            \bigl(\mathbf{1}_I(b)+\mathbf{1}_I(d)\bigr)
            =N_I(n,z),\qquad \mathbf{1}_I(+\infty)=0,
        \end{equation}
        where $N_I(n,z)=\#\{x\in\mca{G}_{n,z}:A_{n,z}(x)\in I\}$. In particular,
        \begin{equation}\label{eq:regular-barcode-birth-bound}
            \#\{[b,d)\in\mca{B}_{n,z}:b\in I\}\le N_I(n,z)<\infty.
        \end{equation}
        Thus admissibility implies local finiteness for regular members without any additional assumption.
    \end{lemma}
    \begin{proof}
        Admissibility makes the action-compatible basis locally finite, and the energy identity makes the differential strictly action-decreasing. The completed filtration therefore satisfies the hypotheses of Proposition~\ref{prop:endpoint-census}, which proves the endpoint identity and both tameness assertions. Omitting the nonnegative death count gives the birth bound. The proof in Appendix~\ref{app:completed-filtrations-barcode-endpoints} also verifies that omitting births at $-\infty$ loses no finite death endpoints for these regular complexes.
    \end{proof}

    \begin{remark}
        Since action is unbounded below, finite-dimensional window quotients need not make $V_{n,z}^t$ finite-dimensional or \emph{$q$-tame} (having finite rank for every structure map $V^a\to V^b$ with $a<b$). The preceding argument proves local finiteness of the barcode through finite cokernels and persistence measures, without invoking an interval-decomposition theorem for pointwise finite-dimensional modules or requiring an additional analytic assumption.
    \end{remark}

    \begin{lemma}[Filtered stability for regular members]\label{lemma:regular-barcode-filtered-stability}
        There is a constant $C_*$ independent of $n$ such that, for regular parameters $z,w$ of an admissible family,
        \begin{equation}\label{eq:regular-barcode-bottleneck-stability}
            d_{\rm bot}(\mca{B}_{n,z},\mca{B}_{n,w})
            \le C_*|z-w|.
        \end{equation}
    \end{lemma}
    \begin{proof}
        Consider the affine homotopy $K_{n,t}+L_{z,t}+\rho(s)(L_{w,t}-L_{z,t})$, with $\rho(-\infty)=0$, $\rho(+\infty)=1$ and $\rho'\ge0$. The action-energy identity in Lemma~\ref{lemma:continuation-action-energy-identity} bounds the filtration shift in either continuation direction by
        \[
            \Delta(z,w):=\int_0^1\|L_{w,t}-L_{z,t}\|_\infty\,\D t
            \le C_L|z-w|.
        \]
        The standard contraction of the concatenated homotopy gives chain homotopies to the identity, whose filtration cost is at most $2\Delta(z,w)$ by the same action-energy estimate. The existence and compactification of the relevant moduli spaces are among the analytic requirements in Definition~\ref{def:analytically-admissible-family}. The maps extend to~\eqref{eq:decreasing-action-chain-completion} because, for any output action threshold, a bounded filtration shift allows contributions only from the finitely many input terms above the corresponding shifted threshold. The same reasoning applies to the chain homotopies. Consequently, the modules $V_{n,z}$ and $V_{n,w}$ are $\Delta(z,w)$-\emph{interleaved}, meaning that structure-compatible maps in both directions shift the parameter by $\Delta(z,w)$ and their composites are the structure maps of shift $2\Delta(z,w)$ \cite[Section~3.2]{ChazalDeSilvaGlisseOudot2016}.

        Combining the interpolation lemma, the box lemma and measure stability in \cite[Lemmas~3.4 and 4.22, Theorem~4.28]{ChazalDeSilvaGlisseOudot2016} yields stability for persistence \emph{measures} without requiring $q$-tameness. By Lemma~\ref{lemma:regular-barcode-endpoint-counts}, both finite interiors contain the entire finite-birth region, including the line $q=+\infty$. The resulting matching restricts to this region because finite-distance matches cannot change a finite birth into $-\infty$, and any unmatched bar there has length at most $2\Delta(z,w)$. This proves~\eqref{eq:regular-barcode-bottleneck-stability}, with $C_*=C_L$ (or any larger fixed constant). Although compactness constants may depend on $n$, the filtration displacement is independent of it.
    \end{proof}

    \begin{proposition}\label{prop:barcode-regularization-local-finiteness}
        For each fixed $n$, the regular barcodes extend uniquely by bottleneck limits to every $z\in\Bpar$. These limits satisfy \eqref{eq:bounded-birth-window-local-finiteness} and
        \begin{equation}\label{eq:regularized-barcode-bottleneck-stability}
            d_{\rm bot}(\mca{B}_{n,z},\mca{B}_{n,w})\le C_*|z-w|
            \quad(z,w\in\Bpar).
        \end{equation}
        In particular, $\mca{B}_n=\mca{B}_{n,0}$ is well-defined.
    \end{proposition}
    \begin{proof}
        Fix $n,z_0$ and regular $z_j\to z_0$. By Lemma~\ref{lemma:regular-barcode-filtered-stability}, the barcodes $\mca{B}_{n,z_j}$ form a bottleneck Cauchy sequence. To construct its limit within the class \eqref{eq:bounded-birth-window-local-finiteness}, pass to a subsequence with successive matchings of costs $\delta_j$ such that $\sum_j\delta_j<\infty$, and follow the matched bars. Along every track continuing indefinitely, the finite endpoints converge, while an infinite death endpoint remains infinite. Keep tracks with positive limiting length and omit those converging to the diagonal. A bar whose length exceeds twice the remaining total matching cost cannot disappear into the diagonal, so this construction gives a limit multiset and bottleneck convergence. Tracks born after a given stage have limiting length at most twice the remaining total cost. Hence, in a bounded birth window, every limit bar of length at least $\varepsilon$ can be tracked from one sufficiently late stage, where only finitely many such bars exist. The limit therefore satisfies~\eqref{eq:bounded-birth-window-local-finiteness}.

        More explicitly, write $I^\delta=\{s:\operatorname{dist}(s,I)\le\delta\}$ and let $b_\varepsilon^I$ denote the count in~\eqref{eq:bounded-birth-window-local-finiteness}. For $0<\delta<\varepsilon/4$ and a sufficiently close regular approximation $\mca{B}_{n,z_j}$, a matching of cost less than $\delta$ gives
        \begin{equation}\label{eq:limit-barcode-local-finiteness}
            b_\varepsilon^I(\mca{B}_{n,z_0})
            \le b_{\varepsilon-2\delta}^{I^\delta}(\mca{B}_{n,z_j})<\infty.
        \end{equation}
        Interleaving two regularizing sequences and applying \eqref{eq:regular-barcode-bottleneck-stability} shows that the limit is independent of the sequence. Passing this inequality to the limit proves~\eqref{eq:regularized-barcode-bottleneck-stability}. If $z_0$ is already regular, the limit equals its original barcode. The argument takes place at the level of diagrams, without assuming a Floer complex with isolated generators at a degenerate parameter.
    \end{proof}
\subsection{The fixed-window entropy upper bound}\label{section:fixed-window-barcode-entropy}
    \begin{definition}[Fixed-window TP barcode entropy]\label{def:fixed-window-translated-point-barcode-entropy}
        For a barcode $\mca{B}$ and a compact interval $I$, write
        \[
            b_\varepsilon^I(\mca{B})
            =\#\{[b,d)\in\mca{B}:b\in I,\ d-b\ge\varepsilon\}.
        \]
        With $\log^+x=\log\max\{1,x\}$, the \emph{fixed-window TP barcode entropy} of the chosen normalized family is
        \[
            h_{\rm TP}^I(f;\alpha;\{K_{n,\bigcdot}\})
            =\lim_{\varepsilon\downarrow0}\limsup_{n\to\infty}
            \frac{1}{n}\log^+b_\varepsilon^I(\mca{B}_n).
        \]
        The outer limit exists by monotonicity. We also put
        \[
            L_{\varepsilon,\mathrm{end}}^I(n)
            =\#\{[b,d)\in\mca{B}_n:d-b\ge\varepsilon,
            \ b\in I\text{ or }d\in I\}.
        \]
        The latter is called a \emph{two-endpoint long-bar count}.
    \end{definition}
    
    Both counts use the length conventions of \cite{CineliGinzburgGurel2024Entropy,CineliGinzburgGurel2024Growth}, with the indicated finite action-window restrictions. Infinite bars are included with length $\infty$.

    \begin{theorem}\label{thm:fixed-window-barcode-entropy-bound}
        Let $(W,\lambda)$ be a Liouville domain with boundary $(Y,\alpha)$, let $f\in\Cont_0^\infty(Y,\ker\alpha)$, and choose an analytically admissible family as in Definition~\ref{def:analytically-admissible-family}. For every compact Reeb-shift interval $I$ and every $\varepsilon>0$, there is a constant $C_{I,\varepsilon}$ independent of $n$ such that
        \begin{equation}\label{eq:fixed-window-long-bar-volume-bound}
            b_\varepsilon^I(\mca{B}_n)
            \le L_{\varepsilon,\mathrm{end}}^I(n)
            \le C_{I,\varepsilon}\operatorname{Vol}_{\dim Y}(\Gamma_n).
        \end{equation}
        Consequently,
        \begin{equation}\label{eq:fixed-window-barcode-entropy-bound}
            h_{\rm TP}^I(f;\alpha;\{K_{n,\bigcdot}\})\le h_{\rm top}(f).
        \end{equation}
    \end{theorem}
    \begin{proof}
        Fix $0<\delta<\varepsilon/4$ and construct a contact tomograph for $I^\delta=\{s:\dist(s,I)\le\delta\}$. Using the constant in Proposition~\ref{prop:barcode-regularization-local-finiteness}, shrink the fixed parameter ball $\Bpar$ so that $C_*|z|<\delta$ throughout the ball. A bar in $\mca{B}_n$ with length at least $\varepsilon$ cannot be matched to the diagonal at a cost less than $\delta$. Instead, it matches a distinct bar in $\mca{B}_{n,z}$ with length at least $\varepsilon-2\delta$, and each finite endpoint moves by less than $\delta$. Thus every bar counted by $L_{\varepsilon,\mathrm{end}}^I(n)$ matches a bar with an endpoint in $I^\delta$.

        Theorem~\ref{thm:translated-point-crofton-bound} gives transversality of the translated-pair equation for almost every parameter, while admissibility provides the regular Floer choices. Lemmas~\ref{lemma:regular-barcode-endpoint-counts} and \ref{lemma:normalized-critical-translated-pairs} then yield
        \[
            L_{\varepsilon,\mathrm{end}}^I(n)
            \le \#\TP(\varphi_z\circ f^n;I^\delta).
        \]
        Indeed, each matched bar has at least one endpoint in this window, so the endpoint census bounds the number of these bars by the number of generators. Integration over $\Bpar$, followed by the contact Crofton estimate, gives
        \[
            \operatorname{Vol}(\Bpar)L_{\varepsilon,\mathrm{end}}^I(n)
            \le C_{I^\delta}\operatorname{Vol}(\Gamma_n).
        \]
        Since the parameter ball and constants were fixed before $n$, this proves \eqref{eq:fixed-window-long-bar-volume-bound}. Taking the upper exponential growth rate and applying Theorem~\ref{thm:cocycle-graph-volume-growth}, followed by the limit $\varepsilon\downarrow0$, establishes \eqref{eq:fixed-window-barcode-entropy-bound}.
    \end{proof}

    \begin{proof}[Proof of Theorem~\ref{thm:intro-barcode}]
        This is the birth-window part of Theorem~\ref{thm:fixed-window-barcode-entropy-bound}.
    \end{proof}

    The window is fixed independently of the iterate. Appendix~\ref{subsection:continuation-and-hamiltonian-cutoffs} presents the relevant continuation estimates and their scope, without asserting exact invariance under unrelated Hamiltonian cutoffs. The completed action filtration also differs from the contact-parameter filtration of \cite{Cant2025}.

\section{Super-exponential growth of short bars}\label{section:superexponential-short-bar-growth}

    Under strong translated nondegeneracy, the unperturbed translated-pair counts equal the generator counts in the completed filtration of Subsection~\ref{section:analytic-admissibility-and-barcodes}. Regularization alone does not establish this equality for degenerate solutions.

    \begin{definition}[Short-bar counts]\label{def:short-bar-window-counts}
        For a bounded interval $I$ and $\varepsilon>0$, define
        \[
            \begin{aligned}
                S_\varepsilon^I(n)
                &=\#\{[b,d)\in\mca{B}_n:b\in I,\ 0<d-b<\varepsilon\},\\
                S_{\varepsilon,\mathrm{end}}^I(n)
                &=\#\{[b,d)\in\mca{B}_n:0<d-b<\varepsilon,
                \ b\in I\text{ or }d\in I\}.
            \end{aligned}
        \]
        These are, respectively, the counts of \emph{short bars born in $I$} and \emph{short bars with an endpoint in $I$}.
    \end{definition}
    
    Both counts refer to the barcode of the full completed filtration rather than that of a quotient complex truncated to $I$. The length convention follows \cite{CineliGinzburgGurel2024Growth}, with the endpoint restriction specified here as part of the definition.

    \begin{lemma}\label{lemma:short-bar-subtraction-bound}
        Assume the hypotheses of Theorem~\ref{thm:intro-barcode} and strong translated nondegeneracy. For every bounded interval $K$, every $\varepsilon>0$ and every $H>h_{\rm top}(f)$, there is a constant $C=C(K,\varepsilon,H,f)$ such that
        \begin{equation}\label{eq:short-subtraction}
            S_{\varepsilon,\mathrm{end}}^K(n)
            \ge \tfrac{1}{2}\#\TP(f^n;K)-Ce^{Hn}
            \qquad(n\ge1).
        \end{equation}
    \end{lemma}
    \begin{proof}
        For regular unperturbed members, Lemma~\ref{lemma:regular-barcode-endpoint-counts} identifies the number of generators in $K$ with the number of finite barcode endpoints in $K$, counted with multiplicity. Each bar contributes at most two such endpoints. Separating the bars of length less than $\varepsilon$ from those of length at least $\varepsilon$, with infinite bars included in the latter class, gives
        \[
            \#\TP(f^n;K)
            \le 2S_{\varepsilon,\mathrm{end}}^K(n)
            +2L_{\varepsilon,\mathrm{end}}^K(n).
        \]
        Applying the two-endpoint estimate of Theorem~\ref{thm:fixed-window-barcode-entropy-bound} to $\overline{K}$ bounds the last count by $C_{K,\varepsilon}\operatorname{Vol}(\Gamma_n)$. Theorem~\ref{thm:cocycle-graph-volume-growth} gives $\operatorname{Vol}(\Gamma_n)\le C_H e^{Hn}$ after increasing $C_H$ to cover finitely many initial iterates, which proves the claim.
    \end{proof}

    \begin{proof}[Proof of Theorem~\ref{thm:intro-short}]
        Suppose $\#\TP(f^n;K)$ has infinite upper exponential growth rate. The entropy of a smooth map on a compact manifold is finite \cite[Chapter~3]{KatokHasselblatt1995}; choose $H>h_{\rm top}(f)$. For every $A>H$, arbitrarily large $n$ satisfy $\#\TP(f^n;K)\ge e^{(A+1)n}$. Lemma~\ref{lemma:short-bar-subtraction-bound} then gives
        \[
            S_{\varepsilon,\mathrm{end}}^K(n)\ge\tfrac14e^{(A+1)n}
        \]
        along an unbounded subsequence. Letting $A\to\infty$ proves the endpoint-window assertion.

        Now let $J\subset I$ be compact, with $I$ open, and set $\rho=\dist(J,\R\setminus I)>0$. For a given $\varepsilon>0$, fix $0<\delta<\min\{\varepsilon,\rho\}$. Any bar of length less than $\delta$ with an endpoint in $J$ has its birth in $I$. Indeed, the assertion is immediate when the birth lies in $J$; when its death $d$ lies in $J$, we have $d-\delta<b<d$. Consequently,
        \begin{equation}\label{eq:inner-window-short}
            S_\varepsilon^I(n)\ge S_{\delta,\mathrm{end}}^J(n).
        \end{equation}
        Applying the first part of the proof to $J$ establishes the claimed growth of bars born in $I$.

        For the blender application, indicate the shift window in the notation for the residual sets and choose
        \[
            \widetilde{\mca{R}}_\infty
            =\mca{R}_\infty(s_*)\cap\mca{R}_\infty(s_*/2).
        \]
        By Theorem~\ref{thm:intro-translated}, both sets are residual in the same $C^1$-open set $\mca{U}$. The second factor provides the growth assumption in the compact interior interval $J=[-s_*/2,s_*/2]\subset I=(-s_*,s_*)$, and Lemma~\ref{lem:radial-filling} provides an explicit Liouville filling for $\alpha=u\alpha_{\rm can}$. The preceding argument therefore proves the final assertion for each map that also satisfies the stated isotopy, nondegeneracy and admissibility hypotheses. This application uses no density or existence statement for the additional intersection.
    \end{proof}

    \begin{corollary}\label{cor:prescribed-short-bar-growth}
        In the setting of Theorem~\ref{thm:intro-translated}, fix $s_*>0$ and a sequence of positive real numbers $(a_n)$. There is a $C^1$-residual set $\mca{R}_a^{\rm short}\subset\mca{U}$ such that every member satisfying the additional hypotheses of Theorem~\ref{thm:intro-short} obeys, for each fixed $\varepsilon>0$,
        \[
            \limsup_{n\to\infty}
            \frac{S_\varepsilon^{(-s_*,s_*)}(n)}{a_n}=+\infty.
        \]
    \end{corollary}
    \begin{proof}
        Apply the explicit Baire construction \eqref{eq:translated-point-open-dense-sets} in the smaller window $(-s_*/2,s_*/2)$ with positive integers $b_n=\lceil a_n+e^{n^2}\rceil$. This construction gives a residual set and an unbounded sequence $n_j$ with at least $n_j b_{n_j}$ translated pairs in that window. For each fixed $\varepsilon$, choose an interior threshold $\delta$ as in \eqref{eq:inner-window-short}. Lemma~\ref{lemma:short-bar-subtraction-bound} then yields
        \[
            S_\varepsilon^{(-s_*,s_*)}(n_j)
            \ge \tfrac{1}{2} n_j b_{n_j}-Ce^{Hn_j}
            \ge \tfrac{1}{4} n_j b_{n_j}
        \]
        for all large $j$. Dividing by $a_{n_j}\le b_{n_j}$ proves the claim. The residual set is chosen before $\varepsilon$, although the threshold and constant in the estimate may depend on $\varepsilon$.
    \end{proof}

    \begin{remark}\label{remark:short-bar-windows-and-quantifiers}
        Growth of generators in $I$ alone does not force growth of births in $I$. For example, take $M_n=\lceil e^{n^2}\rceil$ independent pairs $\partial y_i=x_i$ with actions $A(x_i)=-s_* -\delta_n$ and $A(y_i)=-s_*+\delta_n$, where $0<\delta_n<s_*$ and $\delta_n\to0$. There are $M_n$ generators in $(-s_*,s_*)$, although every barcode interval $[-s_*-\delta_n,-s_*+\delta_n)$ is born outside that window. This algebraic example explains the interior-window hypothesis in Theorem~\ref{thm:intro-short}, without providing a contact counterexample.

        Infinite upper exponential growth rate does not imply that a single map dominates every prescribed super-exponential sequence. Accordingly, the residual sets in Corollary~\ref{cor:prescribed-short-bar-growth} and Theorem~\ref{thm:intro-translated} depend on the prescribed sequence.
    \end{remark}

\appendix
\section{Hamiltonian normalization and continuation}\label{app:hamiltonian-lifts-rabinowitz-functional}

    This appendix verifies the lift and nondegeneracy statements needed for the normalized functional~\eqref{eq:normalized-rabinowitz-functional} and specifies the fixed-iterate analytic input drawn from perturbed Rabinowitz Floer theory.

\subsection{Lifts and critical-point nondegeneracy}

    \begin{lemma}[Compactly supported homogeneous lifts]\label{lemma:compactly-supported-hamiltonian-lift}
        Let $f\in\Cont_0^\infty(Y,\ker\alpha)$ and fix $n\ge1$. The homogeneous lift
        \begin{equation}\label{eq:homogeneous-contact-lift}
            \hat{f}^{\,n}(r,x)=(re^{-G_n(x)},f^n(x))
        \end{equation}
        preserves $\lambda=r\alpha$. For any chosen contact isotopy from $\id$ to $f^n$ and any compact annulus in the symplectization, there is a compactly supported Hamiltonian isotopy on $\widehat{W}$ agreeing with its homogeneous lift on a neighborhood of all trajectories starting in that annulus. Along these trajectories, its Hamiltonian $H_t$ satisfies $\lambda(X_{H_t})=H_t$. The isotopy can be supported in any prescribed open time interval.
    \end{lemma}
    \begin{proof}
        The pullback identity is $re^{-G_n}(f^n)^*\alpha=r\alpha$. If $k_t$ generates the chosen contact isotopy, its homogeneous lift is generated by $rk_t(x)$; see \cite[\S2, preceding Eq.~(2.3)]{AlbersMerry2013}. The Liouville flow embeds the entire symplectization in $\widehat{W}$ \cite[\S1, preceding Lemma~1.9]{AlbersMerry2013}. Images of a compact annulus under the lifted isotopy form a compact subset of this embedded symplectization. Choose a radial cutoff equal to one near that subset and vanishing near the compact core and sufficiently far out. Uniqueness of Hamiltonian trajectories then gives the claimed agreement, throughout which
        \[
            \lambda(X_{H_t})=\D{H_t}(r\partial_r)=r\partial_r(rk_t)=rk_t=H_t.
        \]
        A smooth time reparametrization places the support in the prescribed interval. The cutoff is chosen after $n$, and its size may depend on $n$.
    \end{proof}

    \begin{lemma}\label{lemma:radial-defining-hamiltonian}
        There is a smooth function $\Fsig:\widehat{W}\to\R$ such that
        \begin{equation}\label{eq:radial-defining-hamiltonian}
            \Fsig^{-1}(0)=\Sigma,\qquad
            \Fsig=r-1\text{ near }\Sigma,\qquad
            \D{\Fsig}\text{ has compact support}.
        \end{equation}
        It is constant and negative in the compact core, constant and positive sufficiently far out, and $X_{\Fsig}|_\Sigma=R_\alpha$.
    \end{lemma}
    \begin{proof}
        Extend $r-1$ across a two-sided collar by a smooth monotone radial profile that preserves its sign on either side and is constant outside a slightly larger collar. Extend these constants over the core and the positive end. The Hamiltonian convention~\eqref{eq:hamconv} gives $X_{r-1}=R_\alpha$ near $\Sigma$. This is the defining-function construction used in \cite[\S2.3, Eq.~(2.14)]{AlbersMerry2018}.
    \end{proof}

    \begin{proposition}[Nondegeneracy condition]\label{prop:normalized-hessian-translated-pair-kernel}
        At the critical point of $\mathscr{A}_{n,z}$ corresponding to a translated pair $(x,\eta)$ of $f_{n,z}:=\varphi_z\circ f^n$, the Hessian kernel is naturally identified with the solutions $(\zeta,\sigma)\in T_xY\times\R$ of
        \begin{equation}\label{eq:normalized-hessian-kernel}
            \begin{cases}
                Df_{n,z}(x)\zeta-DR_\eta^\alpha(x)\zeta-\sigma R_\alpha(f_{n,z}(x))=0,\\
                \D{g_{f_{n,z}}}(x)\zeta=0.
            \end{cases}
        \end{equation}
        Consequently, the critical point is Morse nondegenerate exactly when the translated pair is nondegenerate.
    \end{proposition}
    \begin{proof}
        The first variation of~\eqref{eq:normalized-rabinowitz-functional} gives
        \[
            \dot{u}=X_{K_{n,t}+L_{z,t}}(u)+\eta\chi(t)X_{\Fsig}(u),
            \qquad \int_0^1\chi(t)\Fsig(u(t))\,\D{t}=0.
        \]
        Its linearization characterizes the Hessian kernel. Testing the differentiated first variation against arbitrary loop variations gives the linearized differential equation, while variation of the multiplier gives the linearized integral constraint. On the first time interval, $\D{\Fsig}(u(t))\nu(t)$ is constant because $\Fsig$ is preserved by its own Hamiltonian flow. The integral constraint and $\int\chi=1$ therefore imply $\nu(0)=(0,\zeta)$ for some $\zeta\in T_xY$. At the end of that interval, the variation is
        \[
            (0,DR_\eta^\alpha(x)\zeta+\sigma R_\alpha(R_\eta^\alpha x)).
        \]
        Near the critical trajectory, the remaining flow is the inverse homogeneous lift of $f_{n,z}$, so the linearized loop-closing condition is
        \[
            D\hat{f}_{n,z}(1,x)(0,\zeta)
            =(0,DR_\eta^\alpha(x)\zeta+\sigma R_\alpha(R_\eta^\alpha x)).
        \]
        Since $g_{f_{n,z}}(x)=0$, the left side equals $(-\D{g_{f_{n,z}}}(x)\zeta,Df_{n,z}(x)\zeta)$, proving~\eqref{eq:normalized-hessian-kernel}. Conversely, uniqueness for the linear ordinary differential equation associates each solution of that system with a unique solution satisfying the integral and loop-closing constraints. This establishes the asserted isomorphism.
    \end{proof}

\subsection{Analytic input and the scope of continuation}\label{subsection:continuation-and-hamiltonian-cutoffs}

    For fixed $n,z$, the separated time supports place the functional~\eqref{eq:normalized-rabinowitz-functional} in the setting of perturbed Rabinowitz Floer theory. An $\omega$-compatible almost complex structure is an endomorphism $J$ with $J^2=-\id$ for which $\omega(\bigcdot,J\bigcdot)$ is a Riemannian metric. On the positive end, choose $J$ invariant under radial rescaling, preserving $\xi$, and satisfying $J(r\partial_r)=R_\alpha$. These are the \emph{cylindrical contact-type} conditions of \cite[Proposition~2.8]{AlbersMerry2013}, whose convention for $J$ has the opposite sign.

    The multiplier estimate and compactness for action-bounded trajectories follow from \cite[Lemma~2.11, Proposition~2.10 and Theorem~2.9]{AlbersFrauenfelder2010Leafwise}. Their hypotheses hold here because $\widehat{W}$ is exact and convex at infinity, $\lambda|_\Sigma$ is a contact form, $\D{\Fsig}$ and the perturbation have compact support, and $\lambda(X_{\Fsig})>0$ near $\Sigma$. Exactness prevents sphere bubbling. All resulting constants may depend on $n$.

    The finite-window differential and its square-zero identity follow the construction in \cite[\S2.3, Eq.~(2.43)]{AlbersFrauenfelder2010Leafwise}. The energy identity
    \[
        E=\mathscr{A}_{n,z}(c_-)-\mathscr{A}_{n,z}(c_+)
    \]
    implies strict action decrease along every nonconstant trajectory. Floer regularity, gluing, and the parametrized constructions remain explicit requirements of Definition~\ref{def:analytically-admissible-family}, and multiplier-dependent almost complex structures may be needed \cite[Remark~2.17]{AlbersMerry2018}.

    The action is defined on all free homotopy classes because $\omega=\D{\lambda}$; compare \cite[Remark~2.6]{AlbersMerry2018}. A Floer cylinder preserves the free homotopy class of its ends. Since the preceding estimates use no capping disk, the finite-window construction applies to each class. For a regular member, a bounded closed window contains finitely many generators across all classes. Indeed, the translated-pair solution set is closed in $Y\times[a,b]$ and hence compact, while nondegeneracy makes its points isolated. This argument does not extend any particular computation of contractible RFH to noncontractible classes.

    A \emph{continuation map} counts solutions with an $s$-dependent Hamiltonian interpolating between two choices \cite[proof of Theorem~2.16]{AlbersFrauenfelder2010Leafwise}. The following identity controls its filtration displacement.

    \begin{lemma}[Continuation energy]\label{lemma:continuation-action-energy-identity}
        Replace $K_{n,t}+L_{z,t}$ in~\eqref{eq:normalized-rabinowitz-functional} by a smooth compactly supported homotopy $H_{s,t}$, constant for $|s|$ large, and denote the resulting functional by $\mathscr{A}_s$. A finite-energy negative-gradient solution with asymptotics $c_\pm$ satisfies
        \begin{equation}\label{eq:continuation-action-energy-identity}
            E=\mathscr{A}_-(c_-)-\mathscr{A}_+(c_+)
            -\int_{\R\times S^1}\partial_sH_{s,t}(u(s,t))\,\D{s}\,\D{t}.
        \end{equation}
        Regular compact continuation shifts the filtration by at most $\int_{\R\times S^1}\|\partial_sH_{s,t}\|_\infty\,\D{s}\,\D{t}$.
    \end{lemma}
    \begin{proof}
        Differentiation along the negative-gradient solution gives
        \[
            \frac{\D}{\D{s}}\mathscr{A}_s=-\|\grad\mathscr{A}_s\|^2
            -\int_0^1\partial_sH_{s,t}(u)\,\D{t}.
        \]
        Integration and the supremum-norm bound prove both assertions; compare \cite[Lemma~2.7]{AlbersFrauenfelder2010Leafwise}.
    \end{proof}

    The radial-cutoff independence proved in \cite[Lemma~2.15 and the discussion after Eq.~(2.23)]{AlbersMerry2018} applies to their chosen contact-path functional. For unrelated Hamiltonian extensions with the same endpoint contactomorphism, we use the stated continuation estimates without claiming an exact identification of barcodes.

\section{Completed filtrations and barcode endpoints}\label{app:completed-filtrations-barcode-endpoints}

    We justify the endpoint census for the completed filtration of Subsection~\ref{section:analytic-admissibility-and-barcodes}. Coefficients lie in a field, and persistence diagrams are understood in the measure-theoretic sense of \cite[Sections~2.1--2.8]{ChazalDeSilvaGlisseOudot2016}.

    \begin{lemma}[Homology of the decreasing-action completion]\label{lemma:completion-homology-inverse-limit}
        Let $C$ be a chain complex completed in the direction of decreasing action, as in \eqref{eq:decreasing-action-chain-completion}, with a locally finite action-compatible generating set and strictly action-decreasing differential. Put $V^t=H_*(C^{\le t})$. For regular levels $a_j\downarrow-\infty$ and $t$, there is a natural isomorphism
        \[
            V^t\cong\varprojlim_j H_*(C^{(a_j,t]}).
        \]
        Moreover,
        \begin{equation}\label{eq:separated-persistence}
            \bigcap_j\operatorname{im}(V^{a_j}\longrightarrow V^t)=0.
        \end{equation}
    \end{lemma}
    \begin{proof}
        Write $C_j=C^{(a_j,t]}=C^{\le t}/C^{\le a_j}$. The completion gives $C^{\le t}=\varprojlim_j C_j$ degreewise, and each $C_j$ is finite-dimensional. In an inverse sequence of finite-dimensional vector spaces, the images in each fixed term form a decreasing sequence of subspaces that eventually stabilizes. Thus the systems of chain groups, cycles $Z_k(C_j)$, boundaries $B_k(C_j)$ and homology groups satisfy the \emph{Mittag--Leffler condition}, which is precisely this stabilization of images \cite[Definition~12.31.2]{StacksProjectInverseSystems}.

        Apply the exactness of inverse limits under this condition \cite[Lemma~12.31.3]{StacksProjectInverseSystems} to
        \[
            0\to Z_k(C_j)\to (C_j)_k\to B_{k-1}(C_j)\to0,
            \qquad
            0\to B_k(C_j)\to Z_k(C_j)\to H_k(C_j)\to0.
        \]
        Kernels commute with inverse limits, while the first sequence identifies the image of the limiting differential with $\varprojlim B_{k-1}(C_j)$. The second sequence then establishes the homology isomorphism. Finally, the long exact sequence of $0\to C^{\le a_j}\to C^{\le t}\to C_j\to0$ identifies $\operatorname{im}(V^{a_j}\to V^t)$ with the kernel of the projection $V^t\to H_*(C_j)$. Since these projections jointly separate points of their inverse limit, \eqref{eq:separated-persistence} follows.
    \end{proof}
    
    For each completed complex $C$ below, let $\mca{G}$ be its chosen action-compatible generating set, and let $A:\mca{G}\to\R$ assign to each generator its action, or filtration value. Then $C^{\le t}$ consists of completed chains supported on generators $x$ with $A(x)\le t$. In Subsection~\ref{section:analytic-admissibility-and-barcodes}, these are $\mca{G}=\mca{G}_{n,z}$ and $A=A_{n,z}$, with $A_{n,z}(x)=\mathscr{A}_{n,z}(x)$. A generator $x$ corresponding to a translated-point-shift pair $(y,\eta)$ has action $A(x)=\eta$ by Lemma~\ref{lemma:normalized-critical-translated-pairs}.

    \begin{proposition}\label{prop:endpoint-census}
        Under the hypotheses of Lemma~\ref{lemma:completion-homology-inverse-limit}, assume that every bounded action interval contains finitely many generators in total. The persistence diagram has finite birth and finite death counts on bounded intervals. Let $\mca{B}$ be its finite-birth part, including infinite deaths and omitting zero-length intervals. If $K\subset\R$ is bounded, then
        \begin{equation}\label{eq:endpoint-census}
            \#\{x\in\mca{G}:A(x)\in K\}
            =\sum_{[b,d)\in\mca{B}}
            \bigl(\mathbf{1}_K(b)+\mathbf{1}_K(d)\bigr),
            \qquad \mathbf{1}_K(+\infty)=0.
        \end{equation}
        Counts include multiplicity and all homological degrees.
    \end{proposition}
    \begin{proof}
        For regular $a<b$, the window exact sequence gives
        \[
            \dim\operatorname{coker}(V^a\to V^b)\le\dim H_*(C^{(a,b]}),
            \qquad
            \dim\ker(V^a\to V^b)\le\dim H_*(C^{(a,b]}).
        \]
        Summing over degrees removes the homological degree shift in the second inequality. The two finite bounds control, respectively, vertical and horizontal strips of the persistence measure. Thus the module is $v$-tame and $h$-tame in the terminology of \cite[Section~2.8]{ChazalDeSilvaGlisseOudot2016}, and its diagram is defined on the corresponding finite regions. Since the filtration is constant between consecutive action values, finite diagram endpoints occur at those values.

        Fix one such value $c$ and choose regular $c_-<c<c_+$ with no other action values between them. The window complex has zero differential, since all its generators have action $c$. Its long exact sequence gives
        \begin{equation}\label{eq:local-endpoint-census}
            m_c=\sum_k\dim\operatorname{coker}(V_k^{c_-}\to V_k^{c_+})
            +\sum_k\dim\ker(V_k^{c_-}\to V_k^{c_+}),
        \end{equation}
        where $m_c$ is the number of generators at $c$.

        To identify the first term with all births at $c$, consider, for $t\ge c_+$, the finite-dimensional quotient
        \[
            Q^t=\frac{\operatorname{im}(V^{c_+}\to V^t)}{\operatorname{im}(V^{c_-}\to V^t)}.
        \]
        The structure maps of this quotient are surjective, and $Q^{c_+}$ is the indicated cokernel. Its dimension can drop only finitely many times, with the surviving quotient accounting for deaths at $+\infty$. Equivalently, the vertical strip formula for the persistence measure counts precisely $\dim Q^{c_+}$ births at $c$, all of positive length.

        For deaths, put $K_c=\ker(V^{c_-}\to V^{c_+})$, which is finite-dimensional. By Lemma~\ref{lemma:completion-homology-inverse-limit},
        \[
            \bigcap_j\bigl(K_c\cap\operatorname{im}(V^{a_j}\to V^{c_-})\bigr)=0.
        \]
        This decreasing sequence of subspaces is eventually zero. Consequently, the entire death multiplicity $\dim K_c$, measured by horizontal strips, has birth bounded below by a finite $a_j$. Every death at $c$ is therefore counted when births at $-\infty$ are omitted. The absence of other action values in $(c_-,c_+)$ also rules out zero-length bars here, so \eqref{eq:local-endpoint-census} counts exactly all endpoints at $c$. Summing over the finitely many such values in a bounded interval proves \eqref{eq:endpoint-census}.
    \end{proof}

    The finite-dimensional identity behind this proposition is the usual barcode pairing of generators, with two generators for a finite bar and one for an infinite bar; compare \cite[Remark~2.1]{CineliGinzburgGurel2024Growth}. Lemma~\ref{lemma:completion-homology-inverse-limit} provides the additional argument required by the action completion used here.

\end{document}